\documentclass{amsart}
\usepackage{amssymb,amsmath,latexsym,epsfig,amsthm,enumerate}

\newtheorem{thm}{Theorem}[section]
\newtheorem{dfn}[thm]{Definition}
\newtheorem{cor}[thm]{Corollary}
\newtheorem{prop}[thm]{Proposition}
\newtheorem{lem}[thm]{Lemma}
\newtheorem{problem}[thm]{\bf Problem}%

\usepackage{pst-grad} 
\usepackage{pst-plot} 
\usepackage{mathrsfs}

\newcommand{\be}{\begin{equation}}
\newcommand{\ee}{\end{equation}}

\newcommand{\comment}[1]{}

\newcommand{\uor}[2]{#1\stackrel{u}{\longrightarrow} #2}
\newcommand{\vor}[2]{#1\stackrel{v}{\longrightarrow} #2}

\begin{document}

 \title{Spanning trees of 3-uniform hypergraphs}
 
\author[Andrew Goodall]{Andrew Goodall$^1$} 
\address{Unaffiliated.} 
\email{goodall.aj@googlemail.com}
\author[Anna de Mier]{Anna de Mier$^2$}


\address{Departament de Matem\`atica Aplicada II, Universitat Polit\`ecnica de Catalunya, Jordi Girona
  1-3, 08034 Barcelona, Spain.} \email{anna.de.mier@upc.edu}

\begin{abstract}

Masbaum and Vaintrob's ``Pfaffian matrix
tree theorem'' implies that counting spanning trees of a 3-uniform
hypergraph (abbreviated to  3-graph) can be done in polynomial time for a class of ``3-Pfaffian''
3-graphs,  comparable to and related to the class of Pfaffian
graphs. We prove a complexity result for recognizing a 3-Pfaffian
3-graph and describe two large classes of 3-Pfaffian 3-graphs --- one of these is given by a forbidden
subgraph characterization analogous to Little's for bipartite Pfaffian
graphs, and the other consists of a class of partial Steiner triple
systems for which the property of being 3-Pfaffian can be reduced to
the property of an associated graph being Pfaffian.
We exhibit an infinite set of partial Steiner triple systems that are not
3-Pfaffian, none of which can be reduced to any other
by deletion or contraction of triples. 

We also find some necessary or sufficient conditions for the
existence of a spanning tree of a 3-graph (much more succinct than can be obtained by
the currently fastest polynomial-time algorithm of Gabow and Stallmann for
finding a spanning tree) and a superexponential lower bound on the
number of spanning trees of a Steiner triple system.

\end{abstract}


\maketitle



\footnotetext[1]{Research supported in part by the hosting department
  while visiting the other author and Marc Noy at various times in 2008-10.} 
\footnotetext[2]{Research supported in part by projects MTM2008-03020 and DGR2009-SGR1040.}

\section{Introduction}

\subsection{Spanning trees of 3-uniform hypergraphs}

In this paper we investigate the problem of the existence, finding and
counting of spanning trees of 3-uniform hypergraphs (henceforth called
{\it 3-graphs} for short). 
The initial motivation for our work was Masbaum and Vaintrob's Pfaffian matrix
tree theorem~\cite{MV02}. They introduce the notion of an orientation (or
equivalently a sign)  of a spanning tree of a 3-graph. The Pfaffian matrix
tree theorem gives a generating function for signed spanning trees of a 3-graph.
We shall be particularly interested in how this spanning tree orientation
can be used to identify a large class of 3-graphs for which the
problem of counting the number of spanning trees can be done in
polynomial time. This class is comparable to that of Pfaffian
graphs, for which there is a polynomial-time algorithm for counting the number of perfect
matchings. A classical theorem of Kasteleyn~\cite{K67} is that
planar graphs are Pfaffian: can we find a similar class of 3-graphs
for which counting the number of spanning trees can be done in
polynomial time? 

We should be clear at the outset about how we are defining a spanning
tree of a 3-graph, for there are various natural alternatives. (More
detailed definitions of these and other terms from the theory of
hypergraphs are given in Section~\ref{sec:notation} below.)  A
spanning tree of a 3-graph $H$ is an inclusion-maximal subset $T$ of the hyperedges of $H$ that covers all
the vertices subject to the condition that $T$ does not contain a cycle of hyperedges. If $B_H$ is the usual bipartite vertex-hyperedge incidence graph
associated with $H$, then a spanning tree of $H$ in this sense corresponds precisely to a spanning tree of $B_H$ with the property that either all three edges of $B_H$ incident with a given hyperedge belong to the tree or none of them do. 
Alternatively, if each hyperedge $\{a,b,c\}$ of $H$ is represented as
a triangle of edges $ab, bc, ca$ in a graph $G_H$ on the same vertex set as $H$, then a
spanning tree of $H$ corresponds to a cactus subgraph of $G_H$ covering all
vertices. 
See~\cite{Abd04} for a generalization of the Masbaum--Vaintrob theorem
to arbitrary hypergraphs in which spanning trees are now cacti with
cycles of any odd length and not just triangles.

Spanning trees of 3-graphs differ in fundamental ways from spanning
trees of ordinary graphs: a closer correspondence is to be found with
perfect matchings. Whereas for spanning trees of graphs the
problems of the existence, finding and counting of spanning trees each have a
straightforward polynomial-time algorithm, the same is not true for
spanning trees of 3-graphs. 

The augmenting path algorithm finds a maximum matching of a bipartite
graph in polynomial time. Consequently, both the problem of whether
there is a perfect matching of a bipartite graph and the problem of
finding one can be solved in polynomial time. 
Edmonds' maximum matching algorithm~\cite{E65} solves in polynomial time the existence and search
problems for whether an arbitrary graph has a perfect
matching. 
   
Lov\'asz's matroid matching algorithm~\cite{L78, L80} provides a
polynomial-time algorithm solving the problem of the existence and
finding of a spanning tree of a 3-graph. 
However, since it solves such a general and complicated problem, the algorithm is involved, has running time a polynomial
of high degree and is not optimal when restricting attention from
linear matroids to the graphic matroids underlying the case of
3-graphs.  
The augmenting path algorithm for
linear matroids of Gabow and Stallmann~\cite{GS86} has running time
$O(mn^2)$ with $O(mn)$ space for graphic matroids of rank $n$ and size
$m$, improved to using $O(m)$ space (alternatively $O(mn\log^6n)$ time using $O(m \log^4n)$ space) by the
same authors in~\cite{GS85}. In this paper we give some
straightforward necessary or sufficient conditions that give simple
criteria for the existence of a spanning tree of a 3-graph and in the
case of Steiner triple systems a superexponential lower bound on the
number of spanning trees.  

Our focus then turns to the problem of counting spanning trees of 3-graphs.
This problem is $\#\mathsf{P}$-complete even for a very restricted
  class of 3-graphs, which is a consequence of the fact that counting perfect matchings is $\#\mathsf{P}$-complete for general graphs~\cite{Val79}.
Masbaum and Vaintrob define an orientation or sign of a spanning tree
of a 3-graph using orientations of hyperedges in a way that closely
follows the definition of the sign of a perfect matching, as
elucidated by Hirschman and Reiner~\cite{HR04}. Just as the existence
of a Pfaffian orientation of the edges of a graph enables the number
of perfect matchings of a graph to be computed in polynomial time, so
the existence of what we shall call a ``3-Pfaffian orientation'' of a 3-graph allows the number of spanning trees to be calculated in polynomial time.  
This observation was made by Caracciolo et al. in the conclusion of their paper~\cite{CMSS08}.  

Having identified a property of 3-graphs that enables counting of
spanning trees to be done in polynomial time, how quickly can we
verify that a graph has this property? Compare the case of Pfaffian
graphs:  it is not known whether there is a polynomial-time checkable certificate for a
graph to have a Pfaffian orientation.  
Vazirani and Yannakakis~\cite{VY89} show that the problem of determining whether a
graph $G$ has a Pfaffian orientation and that of determining whether a given
orientation of $G$ is Pfaffian are polynomial-time equivalent. They appeal to 
Lov\'asz' polynomial-time algorithm \cite{Lovasz87} for computing the
binary rank and finding a basis of the vector space of matchings of a graph. 
They also show that the problem of deciding whether a graph has a Pfaffian orientation is in {\sf co-NP}.
We show that the problem of deciding the existence of a 3-Pfaffian
orientation is also in {\sf co-NP}, but we do not know if it is equivalent to deciding if a given orientation of hyperedges is 3-Pfaffian.

Although checking whether a graph is Pfaffian is not known to be
polynomial time, Little~\cite{Little75} gave a structural
characterization of 
Pfaffian bipartite graphs as those that do not contain an even subdivision of $K_{3,3}$
with a perfect matching in the complement. A natural question is
whether there is any similar characterization of 3-Pfaffian 3-graphs:
we prove such a characterization for a special subclass of tripartite
3-graphs. Whether tripartite 3-Pfaffian 3-graphs in general have a
similar description in terms of forbidden subgraphs remains open.  

\subsection{Outline of the paper}

In Section~\ref{sec:notation} we introduce some of the basic notions and
notation required in the paper.
We refer to~\cite{RT06} for a recent survey of the topic of Pfaffian
orientations, and~\cite{LP86} for matching theory. 

In Section~\ref{sec:existence_counting} we present some elementary
results about the problem of deciding if there is a spanning
tree of a 3-graph and about the problem of 
counting them. We begin in
Subsection~\ref{sec:polytime} with a summary of what is known about the
complexity of these problems in general.
In Subsection~\ref{sec:sp_trees_complete} we consider the case of the
complete 3-graph, for which we can enumerate the number of spanning
trees, and, more importantly, thereby establish in Lemma~\ref{lem.sp trees pm + function} 
a correspondence
between spanning trees of a 3-graph and perfect matchings of a graph that is basic to the
rest of the paper. In Subsection~\ref{sec:nec_and_suff_sp_trees} we describe
some straightforward necessary or sufficient conditions for the
existence of a spanning tree of a 3-graph.
Theorem~\ref{thm:lower_bound_STS} gives a lower bound
on the number of spanning trees of a Steiner triple system. 
 
In Section~\ref{sec:oriented} we initiate our study of orientations of spanning
trees of 3-graphs and the property of a 3-graph having a ``3-Pfaffian
orientation,'' which by
Masbaum and Vaintrob's Pfaffian matrix-tree theorem~\cite{MV02}
implies a polynomial-time algorithm for counting spanning trees. We
begin in Subsection~\ref{sec:orns_sp_trees} by defining orientations of
spanning trees, which are defined relative to an orientation of triples. 
Theorem~\ref{thm:positive_negative} gives an explicit formula for the number of positively
and negatively oriented spanning trees of the complete 3-graph under a canonical orientation of its triples.
In Subsection~\ref{sec:tree_generating_polynomials} we introduce the
notion of a ``3-Pfaffian orientation'',
analogous to a Pfaffian orientation of a graph. In fact in
Theorem~\ref{thm:1-susp} we see that if we make a 3-graph $H$ by adding an extra vertex to
every edge of a graph $G$ then a 3-Pfaffian orientation of $H$
corresponds exactly to a
Pfaffian orientation of $G$. In Subsection~\ref{sec:complexity} we
prove that deciding if a 3-graph has a 3-Pfaffian orientation is in
${\sf co}$-${\sf NP}$.

In Section~\ref{sec:suspensions} we consider a family of 3-graphs for which we can
characterize the property of having a 3-Pfaffian orientation in terms
of forbidden subgraphs, similar to Little's characterization of
Pfaffian bipartite graphs (Theorem~\ref{thm:2-susp_3-Pfaffian} and Corollary~\ref{cor:2-susp_obstructions}).
  
In Section~\ref{sec:STS} we find a large class of
partial Steiner triple systems that have 3-Pfaffian orientations
(Theorem~\ref{thm:PSTS_3-Pfaff}) and also describe an infinite family of partial Steiner triple
systems that do not have a 3-Pfaffian orientation
(Theorem~\ref{thm:interlaced}). Furthermore, we prove that this second family
cannot be reduced by deletion and contraction of triples to a finite set of
non-3-Pfaffian 3-graphs.

Finally, in Section~\ref{sec:open_problems}, we highlight some open
problems.

\medskip


\section{Notation and terminology}\label{sec:notation}


A {\em 3-graph\/} is a 3-uniform hypergraph $H=(V, \Delta)$, where
$\Delta\subseteq\binom{V}{3}$. There are no repeated hyperedges and no
hyperedges of size 2 or 1. We shall use the name {\em
  triple\/} for a hyperedge of $H$. 
The {\em underlying graph\/} of a 3-graph $H=(V,\Delta)$ is the
multigraph $G=(V,E)$ with edge set $E=\{\{a,b\}:\exists \ {c\in V}\;
\{a,b,c\}\in \Delta\}$, an edge $\{a,b\}$ occuring with multiplicity $|\{c\in
V: \{a,b,c\}\in \Delta\}|$. We identify a triple of $H$ with its corresponding
triangle in the underlying graph $G$. We write $abc$ for the triple
$\{a,b,c\}$ of $H$ or triangle of $G$ and $ab$ for the edge
$\{a,b\}$ of $G$.

{\em Deleting} a triple
$abc\in\Delta$ gives the 3-graph $H\backslash abc=(V,\Delta\setminus
 abc)$.
A {\em sub-3-graph} of $H$ is a 3-graph obtained from $H$ by deleting
some subset of triples. 
{\em Contracting} a triple $abc$ gives the 3-graph $H/abc=(V\setminus\{b,c\},\Delta')$
where $\Delta'$ is defined as follows. A triple $ijk$ belongs to
$\Delta'$ if (i) $ijk$ and $abc$ are
  disjoint,  or (ii)  $ijk$ is
obtained from a triple that meets
$abc$ in one vertex by relabelling this common vertex by $a$ if it
is equal to $b$ or $c$.  
In other words, to form $H/abc$ from $H$ we set $a=b=c$ and remove all triples that have
decreased in size to a pair or singleton and also any repeated
triples. 
 
In terms of the underlying graph $G$ of $H$, deleting a triple $abc$
of $H$ corresponds to deleting the edges $ab, bc, ca$ of
$G$. Contracting $abc$ corresponds to contracting $ab,bc, ca$ and
removing any edges that are no longer an edge of a triangle.

The {\em degree} of a vertex $a\in V$ in $H$ is defined by
$d(a)=\#\{t\in \Delta:a\in t\}$, equal to half the degree of $a$ in the underlying graph $G$. 
The {\em multiplicity} of a pair $ab\in\binom{V}{2}$ in $H$ is defined
by $m(ab)=\#\{t\in\Delta:\{a,b\}\subseteq t\}$, equal to the multiplicity of the edge $ab$ in the underlying graph $G$.

A {\em path\/} in a 3-graph $H=(V,\Delta)$ is an alternating sequence of
$\ell+1$ distinct vertices and $\ell$ distinct triples, $a_0, t_1, a_1,\ldots, a_{\ell-1}, t_\ell,
a_\ell$, with the property that $a_{i-1}\in t_{i}\ni a_{i}$ for
$i\in[\ell]$. A path is usually identified with its set of triples
$\{t_1,\ldots, t_\ell\}$. 
 Observe that it is not required that a path with $\ell$ triples
spans $2\ell +1$ vertices, although most of the paths that appear in the paper
have this property.

The 3-graph $H$ is {\em connected} if for each pair
of vertices $u,v\in V$ there is a path $u,t_1,\ldots , t_\ell, v$ in $H$ that
joins them. $H$ is connected if and only if its underlying graph is connected.

A {\em cycle\/} in $H$ is a closed path, i.e.,\ an alternating sequence
of $\ell$ distinct 
vertices and $\ell$ distinct triples $a_0, t_1,\ldots, a_{\ell-1}, t_\ell$
terminated by the starting vertex $a_\ell=a_0$, with the
property that $a_{i-1}\in t_{i}\ni a_{i}$.
A cycle is usually identified with its set of triples $\{t_1,\ldots,
t_\ell\}$. Two triples sharing two vertices form a cycle.

A {\em forest} of $H$ is a set of triples $T\subseteq\Delta$ with the
property that there is no cycle $C\subseteq T$. Between any two
vertices in a forest there is at most one path. A {\em spanning tree}
of $H$ is a sub-3-graph $T$ containing no cycles such that $\bigcup T=V$,
i.e.,\ a
connected forest spanning $V$. If
$H$ has a spanning tree then $|V|$ is necessarily odd and $T$ contains
$\frac{|V|-1}{2}$ triples. 
The connected 3-graph on $\{u,v,a,b,c\}$
with triples $uva, uvb, uvc$ has no spanning tree.
A {\em leaf} of a tree $T$ is a triple with two vertices of degree $1$ (belonging to no other triple of $T$). A spanning tree of $H$ has at
least one leaf $abc$, and at least two leaves if $|V|\geq 5$.  The 3-graph
$T-\{b,c\}$ obtained by deleting vertices $b,c$  is a spanning tree of $H-\{b,c\}$ if
and only if $abc$ is a leaf of $T$ for some $a$ and where $b,c$ have
degree $1$.


\section{Elementary results on the existence and counting of spanning trees}\label{sec:existence_counting}


%
\subsection{Complexity of existence, finding and counting of spanning trees
    of 3-graphs}\label{sec:polytime}

Given a 3-graph $H=(V,\Delta)$ and triples $abc$ put in arbitrary linear order $a<b<c$, define the subgraph $G'$ of its underlying graph $G=(V,E)$ on edge set
$E'=\{ab, ac: abc\in\Delta, a<b<c\}$ of size $2|\Delta|$. 
Partition $E'$ into pairs $ab,ac$ with $abc\in\Delta, a<b<c$. A matching of
the graphic matroid defined by $G'$ is a forest of $G'$ such that for
each $abc\in\Delta$ with $a<b<c$ if
$ab$ belongs to the forest then so does $ac$. A maximum matching
has the greatest number of pairs possible. The 3-graph $H$ has a
spanning tree if and only if the maximum matching has size
$\frac{|V|-1}{2}$. Thus the problem of determining whether a 3-graph
has a spanning tree is a special case of the matroid matching problem.
As mentioned in the introduction, this gives a
polynomial-time algorithm for finding a spanning tree of a 3-graph.

For $k\geq 4$ the problem of deciding if a $k$-uniform hypergraph has
a spanning tree is {\sf NP}-complete \cite{AF95}.

Counting spanning trees of a 3-graph is \#{\sf P}-complete. This follows since
counting perfect matchings of a graph is a $\#{\sf P}$-complete problem in general \cite{Val79} and this reduces to the problem of counting spanning trees for the class of 3-graphs with the property that there is a vertex that is contained in all triples. 

On the other hand, counting perfect matchings is polynomial time for
the class of graphs that have a Pfaffian orientation.
One of the aims of this paper is to develop the analogous notion of
 a Pfaffian orientation for 3-graphs and thereby characterize a
 class of 3-graphs with the property that counting spanning trees has a polynomial-time algorithm.

\subsection{Spanning trees of complete $3$-graphs}\label{sec:sp_trees_complete}

For a 3-graph $H=(V,\Delta)$ let $\mathcal{T}(H)=\{T\subseteq \Delta:
\mbox{\rm $T$ is a spanning tree of $H$}\}$.
Note that $\mathcal{T}(H\backslash abc)=\{T\in\mathcal{T}(H): \, abc\not\in T\}$
and there is a bijection between
$\mathcal{T}(H/abc)$ and $\{T\in\mathcal{T}(H):\, abc\in T\}.$
If $abc$ is in no spanning tree of $H$ then
$\mathcal{T}(H)=\mathcal{T}(H\backslash abc)$. If $abc$ is in every
spanning tree of $H$ then contracting the triple $abc$ defines a
bijection from $\mathcal{T}(H)$ to $\mathcal{T}(H/abc)$.

Fix $n\in\mathbb{N}$ and denote by $K_{2n+1}^{(3)}$ the complete $3$-graph with vertex set $[2n+1]$ and triple set $\binom{[2n+1]}{3}$, the set all $3$-subsets of $[2n+1]$.
For short we write $\mathcal{T}$ for the set of spanning trees of
$K_{2n+1}^{(3)}$. 

The following result can be found for example in~\cite{SS06}, but we
include a proof as it prepares the ground for the next lemma and for
Theorem~\ref{thm:positive_negative} later.
\begin{thm}\label{thm:pruefer} The number of spanning trees of $K_{2n+1}^{(3)}$ is given by
$$|\mathcal{T}|=(2n-1)!!(2n+1)^{n-1}.$$ 
\end{thm}
\begin{proof}
The proof uses a similar construction to the Pr\"{u}fer code for
spanning trees of ordinary graphs.

A tree spanning at least five vertices always has at least
two leaves; a rooted tree spanning five or more vertices has at least
one leaf not containing the root as a vertex of degree 1. 

Suppose we are given a spanning tree $T$ on $[2n+1]$. 
 We remove triples from
  $T$ leaf by leaf in a
  canonical way until we are left with a tree consisting of just one
  triple.
At the end of the algorithm described below we obtain a sequence
$\gamma=\gamma_{n}\in[2n+1]^{n-1}$ and a perfect matching $M=M_n$ of
$[2n]$. If $n=1$, we take $\gamma$ to be the empty sequence and $M=\{ 12\}$. 
For $n\geq 2$, the algorithm proceeds as follows.

\begin{enumerate}[(1)]
\item Intialize $\gamma_1$ as the
  empty sequence, $M_0$ as the empty matching and $T_1=T$ as the
  spanning tree of $K_{2n+1}^{(3)}$ that is to be encoded. Root $T$ at
  vertex $2n+1$.

Start with $i=1$.
\item At step $i$ consider the rooted tree $T_i$.
Remove the leaf containing the smallest vertex label in $T_i$ while
 not containing the root $2n+1$ as a vertex of degree $1$,
 thereby obtaining the next rooted tree $T_{i+1}$. (If a leaf
  contains $2n+1$ as a vertex of degree $1$ it is ignored and the leaf
  with the next smallest vertex is taken.) Record
as $c_i$ the vertex of degree greater than $1$ in this leaf and set
$\gamma_{i+1}=\gamma_i c_i$. The other two
vertices of degree $1$ in the leaf $a_ib_ic_i$ are paired in the matching $M_{i}=M_{i-1}\cup\{a_ib_i\}$.  
\item If the remaining 
tree $T_{i+1}$ has only one triple (i.e., $i=n-1$) then this triple
takes the form $a_nb_n(2n\!+\!1)$; in this case set
$M=M_n=M_{n-1}\cup\{a_nb_n\}$, $\gamma=\gamma_{n}$, 
  and stop.  Otherwise
increment $i$ to $i+1$ and go to (2).
\end{enumerate}
Conversely, given a sequence $\gamma=c_1c_2\ldots c_{n-1}\in[2n+1]^{n-1}$ and a
perfect matching $M=\{a_1b_1,\ldots, a_{n}b_n\}$ of $[2n]$ a unique spanning
tree of $K_{2n+1}^{(3)}$ is constructed as follows.
\begin{enumerate}[(1)]
\item Initialize $i=1$, $\gamma_1=\gamma$, $M_1=M$, $T_1$ the empty tree (no
  triples or vertices).
\item   Find the vertex $a_i$ with smallest label that does not occur as
  an element of the sequence $\gamma_i$ and that occurs in the
  matching $M_i$, but is not paired with $c_i$.  Let  $b_i$ be the vertex such
that $a_ib_i\in M_i$. Set $T_{i+1}=T_i\cup\{a_ib_ic_i\}$,
$M_{i+1}=M_i\setminus \{a_ib_i\}$ and $\gamma_{i+1}=c_{i+1}\ldots c_{n-1}$.

\item After step $i=n-1$ the sequence $\gamma_n$ is empty and
  $M_n=\{a_nb_n\}$. Set $T=T_n\cup \{a_nb_n(2n+1)\}$ and stop.  Otherwise, increment $i$ to $i+1$ and
  go to (2). 
\end{enumerate}

Spanning trees of $K_{2n+1}^{(3)}$ are thus in bijection with pairs
$(\gamma,M)$, where $\gamma\in [2n+1]^{n-1}$ and $M$ is a perfect matching of
$[2n]$. Since there are $(2n-1)!!$ such perfect matchings, the result follows.
\hspace*{\fill} \end{proof}

The first part of the proof of Theorem~\ref{thm:pruefer} can be applied to any
$3$-graph $H$, yielding a correspondence between spanning trees of $H$ and
pairs $(M,f)$, where $M$ is a perfect matching of $H-v$ and $f:M\rightarrow V$
is a function satisfying a certain condition.

\begin{lem}\label{lem.sp trees pm + function}
Let $H=(V,\Delta)$ be a $3$-graph with underlying graph $G$, and let $v\in V$. 
Given a spanning tree $T$ of $H$, there is a unique perfect matching $M$
of $G-v$ and a function $f:M\rightarrow V$ such that the set of triples of $T$
is equal to $\{ijf(ij):ij\in M\}$. Conversely, a perfect
matching $M$ of $G-v$ and a function $f:M\rightarrow V$ determine a spanning
tree of $H$ if
$\{ijf(ij):ij\in M\}\subseteq \Delta$ and there is no set
of edges
$\{i_0j_0,\ldots, i_{\ell-1}j_{\ell-1},\, i_\ell j_\ell=i_0j_0\}\subseteq M$ such that
$f(i_{m-1} j_{m-1})\in\{i_m,j_m\}$ for $m\in[\ell]$. 
\end{lem} 
\begin{proof}
Rooting a spanning tree $T$ of $H$  at the vertex $v$, we construct a unique
perfect matching $M$ of $G-v$ and associated function $f:M\rightarrow
V$ as follows. 

If
$|V|=3$ then $T=\{vij\}$ and set $M=\{ij\}$ and $f(ij)=v$. Assume now
that $|V|>3$. Then every leaf of $T$ has one vertex of degree greater
than $1$, by which it is attached to the rest of the tree, and the
remaining two vertices are of degree $1$. Let $ijk$
be a leaf
of $T$ with vertices $i,j$ of degree $1$.  Remove this
leaf from $T$. Inductively the remaining tree $T\backslash ijk$ determines a unique perfect
matching $M'$ of $G-\{v,i,j\}$ and function $f:M'\rightarrow V\backslash\{i,j\}$. Extend
$M'$ to a perfect matching $M$ of $G-v$ by adding the edge $ij$ and the
function $f$ by setting $f(ij)=k$.

Conversely, given a perfect matching $M$ of $G-v$ and a function
$f:M\rightarrow V$, the 3-graph on $V$ having as  set of triples
$T=\{ij f(ij):ij\in M\}$ is a spanning tree of $H$ if $T\subseteq
\Delta$ and there is no cycle  of triples. It is easy to see that this amounts
to the condition on $f$ in the statement of the theorem. For such an $f$, the $3$-graph $T$ is a tree
with $(|V|-1)/2$ triples, and therefore it  spans the $|V|$ vertices of $H$.
\hspace*{\fill} \end{proof}

\subsection{Necessary or sufficient conditions for the existence of spanning trees}\label{sec:nec_and_suff_sp_trees}

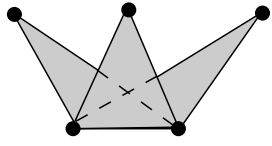
\begin{figure}[ht]
\caption{Smallest connected 3-graph on an odd number of vertices without a
  spanning tree. (Shaded triangles are triples.) }\label{fig:twin}
\begin{center}

\scalebox{1} 
{
\begin{pspicture}(0,-0.9)(3.51,0.9)
\definecolor{color1089b}{rgb}{0.8,0.8,0.8}
\pspolygon[linewidth=0.02,fillstyle=solid,fillcolor=color1089b](3.38,0.75)(0.88,-0.79)(2.3,-0.79)
\pspolygon[linewidth=0.02,fillstyle=solid,fillcolor=color1089b](0.92,-0.75)(2.28,-0.73)(0.08,0.75)
\psdots[dotsize=0.2](0.1,0.73)
\psdots[dotsize=0.2](3.4,0.75)
\pspolygon[linewidth=0.02,fillstyle=solid,fillcolor=color1089b](1.6,0.83)(0.86,-0.79)(2.28,-0.77)
\psdots[dotsize=0.2](1.62,0.79)
\psdots[dotsize=0.2](0.88,-0.79)
\psdots[dotsize=0.2](2.28,-0.79)
\psline[linewidth=0.02cm,fillcolor=color1089b,linestyle=dashed,dash=0.16cm 0.16cm](0.92,-0.69)(2.0,-0.09)
\psline[linewidth=0.02cm,fillcolor=color1089b,linestyle=dashed,dash=0.16cm 0.16cm](1.2,-0.01)(2.22,-0.75)
\end{pspicture} 
}
\end{center}
\end{figure}

The most straightforward necessary conditions for the existence of a
spanning tree of a 3-graph $H=(V,\Delta)$ is that $H$ is connected and that
$|V|$ is odd. The $3$-graph in Figure~\ref{fig:twin} shows that these
conditions are not sufficient. 

Our first non-trivial condition for the existence of
spanning trees is a sufficient one and is as follows.
\begin{thm}\label{thm.sp trees exist} Suppose $H=(V,\Delta)$ is a
  3-graph such that $|V|$ is odd and each pair of vertices has
  multiplicity at least 1 in $H$. 
Then $H$ has a spanning tree.
\end{thm}
\begin{proof}
Assume $T\subseteq \Delta$ is a tree of maximum size and suppose that
$|T|<\frac{|V|-1}{2}$.
Let $U\subset V$ be the set of vertices not spanned by $T$. Then
$|U|$ is even, containing at least two vertices $u,v$. 
Since there is some triple
containing $\{u,v\}$, there is $w\in V$ such that $uvw\in\Delta$ and in
fact $w\in U$ for otherwise we could add the triple $uvw$ as a leaf to $T$ and
obtain a larger tree of $H$.

Set $S=\{uvw\}$, vertex-disjoint from $T$. 
For any leaf $abc$ of $T$ with vertices $a,b$ of degree $1$ in $T$
there is a triple $uai$ containing the pair
$\{u,a\}$. By the remark in the previous paragraph $i\in V\setminus U$. If $uai$ is a triple for some $i\neq b$ then
deleting $abc$ from $T$ and adding the triples $uai$ and $uvw$ gives a
larger tree, contradicting the fact that $T$ has maximum size. 
So we may assume that the only triple that contains $u$ and at least one of
$a,b$ is $uab$, and that this is true for every leaf $abc$ of
$T$. We then remove
all the leaves $abc$ of $T$ and put the triples $uab$ in $S$. 

We repeat this argument, at each stage looking at triples containing
$u$ and vertices of degree $1$ in the leaves of what is left of the
initial tree $T$. There are just two possible outcomes: either
(i) at some stage we can join the remaining subtree of $T$
and the tree $S$ containing $uvw$ by a triple to make a larger tree than
the original tree $T$, or (ii) we remove all the leaves of $T$ and end
up with a larger tree $S$ that spans all but one
of the vertices that are spanned by $T$ and also the vertices
$u,v,w$. Both possibilities
 contradict the hypothesis that $T$ has maximum size. 

Hence the maximum tree $T$ spans all the vertices of $H$, i.e., $T$ is a spanning tree of $H$.\hspace*{\fill} \end{proof}

An extremal case of Theorem~\ref{thm.sp trees exist} is when each
pair of vertices  is contained in exactly one triple, i.e., $H$ is a Steiner triple
system. The condition on the multiplicity of pairs of vertices implies  that a
Steiner triple system on $n$ points also has the property that
every vertex is of degree $\frac{n-1}{2}$, and that $n$ is congruent with $1$
or $3$ modulo $6$.  R.M. Wilson~\cite{Wilson74}
 showed that the number  of
non-isomorphic Steiner triple systems on $n\equiv 1 \; \mbox{ or } \,
3\pmod 6$ points lies between $(e^{-5}n)^{n^2/12}$ and
$(e^{-\frac12}n)^{n^2/6}$. (Given the truth of the then conjecture of Van der Waerden
on the size of permanents, Wilson improved the lower bound, and
further conjectured that the actual number is in fact asymptotically
$(e^{-\frac12}n)^{n^2/6}$.) There is just one isomorphism class for
$n\in\{3,7, 9\}$, two for $n=13$, eighty for $n=15$.

For Steiner triple systems  we can not only assert the existence of a spanning tree but
also give a superexponential lower bound on the number of spanning trees.

\begin{thm}\label{thm:lower_bound_STS}
If $H=(V,\Delta)$ is a Steiner triple system on $|V|=n$ vertices then
$H$ has $\Omega((n/6)^{n/12})$ spanning trees.
\end{thm}

\begin{proof}
Brouwer~\cite{Brouwer81} proved that any Steiner triple system on $n$
vertices has a transversal (set of pairwise disjoint triples) covering all but $5n^{2/3}$ vertices, and
Alon, Kim and Spencer~\cite{AKS97} improved this to all but
$O(n^{1/2}{\rm ln}^{3/2}n)$ vertices.  
Let $P\subseteq\Delta$ be such a
set of pairwise disjoint triples that together cover
$U\subseteq V$, with $|U|=n-k$ and $k=o(n)$. Let $r=(n-1)/2$.

We give a procedure that generates $\prod_{i=0}^s  (r-k-1-6i)$
spanning trees, where $s$ is the largest integer such that $r-k-1-6s>0$ ($s$ is $n/12-o(n)$). 
Unfortunately, this procedure may
give repeated trees; we then show that each tree cannot appear
more than $n/6$ times.  Recall that in a Steiner triple system every
vertex belongs to $r$ triples. In $H_{|U}$ every vertex belongs to
at least $r-k$ triples. Let $u_0$ be a vertex of $U$.  The
construction of a spanning tree consists in first using $P$ to construct a ``comb-like''  tree of
 $H_{|U}$  and then extending this tree to a spanning tree of $H$. So let us begin by considering
the restriction $H_{|U}$.
 Let $t_0$ be any
 triple containing $u_0$, subject only to the condition that
 $t_0\not\in P$. 
 Say  $t_0=\{u_0,u'_0,u''_0\}$. Let
 $p_0,p_{1},p_{2}$ be the triples in $P$ that contain
 $u_0,u'_0,u''_0$, respectively. Clearly the triples $t_0,p_{0},p_{1},p_{2}$
 form a tree $T_0$. Let $u_1$ be any of the four vertices in $(p_{1}\cup
 p_{2}) \backslash t_{0}$. There are at least $r-k-7$ triples that
 contain $u_1$ but no other
 vertex of $T_0$. Let $t_1=\{u_1,u'_1,u''_1\}$ be any one of them. Let $p_3$ and
 $p_4$ be the triples in $P$ that contain $u'_1$ and $u''_1$,
 respectively. Let $T_1=T_0\cup\{t_1,p_3,p_4\}$. We proceed recursively in this way as long as $r-k-1-6i$ is positive: we choose $u_i$ to be any of the four
 vertices in $p_{2i-1}\cup p_{2i}$ that are not in $t_{i-1}$ and we choose $t_i=\{u_i,u'_i,u''_i\}$ a
 triple containing $u_i$ and no other vertex in $T_{i-1}$. Then we take the two
 triples $p_{2i+1},p_{2i+2}$ in $P$ that contain $u'_i,u''_i$ and set $T_i=T_{i-1}\cup\{t_i,p_{2i+1},p_{2i+2}\}$. 

Once we have a tree $T_s$, covering $6s+9$ vertices, we need to complete it to a
spanning tree of $H$. We repeatedly use the following claim.

\medskip

{\sc Claim.} Let $T$ be a tree of $H$ and let $W$ be the set of vertices not
spanned by $T$. Then there are vertices $a,b$ of $W$ such that the triple that
contains them has its third vertex in $T$.

\medskip

\textit{Proof of the claim.} Suppose it were not the case. Then the
triples on $W$ would form a
Steiner triple system. But since $W$ has even cardinality this is
impossible. \hspace*{\fill} $\Box$

Therefore, by adding a leaf at a time, we can complete $T_s$ to a
spanning tree of $H$. There may be many ways of completing
$T_s$, but we just take one of them arbitrarily.

If we fix the starting vertex $u_0$, by applying the procedure just described we
obtain  $4^s\prod_{i=0}^s (r-k-1-6i)$ spanning trees of $H$. Indeed, at step
$i$ we need to choose one of four vertices and then we know that this vertex
belongs to at least $r-k-1-6i$ triples that are contained in $H|_U$ but do not
contain any vertex already in the tree.  It could be, however, that the same
tree is produced several times. For instance, the tree in
Figure~\ref{fig:skeleton} could appear in two different ways.

\begin{figure}[ht]
\begin{center}
\includegraphics[width=11cm, scale=1.7]{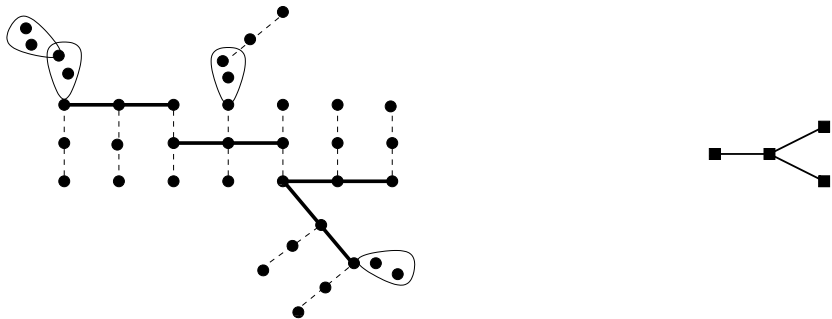}
\caption{A spanning tree and its skeleton (here $s=2$). Red triples are shown as thick long
  lines, blue triples are thin dashed lines and green triples are depicted as bags. }\label{fig:skeleton}
\end{center}
\end{figure}

Next we bound the number of possible repetitions of a given spanning
 tree $T$.
Let us first colour the triples of $T$ in the following way. The triples from $P$ are coloured blue; the triples entirely
contained in $U$  and that intersect three triples of $P$ are coloured red, and the remaining triples
are coloured green. Observe that green triples are only included in the final
stage of the construction of a tree (when the claim is used), whereas blue and red triples can appear both during
the first steps of the construction and also at the end.

The \emph{skeleton} of $T$ is the graph whose vertices are the
red triples and where two vertices are adjacent if there is a blue edge
in the tree intersecting the corresponding red triples in different
vertices.
The skeleton is a forest (it will be a tree if there is only
one red triple containing $u_0$); root each component of the forest at the
vertex corresponding to the triple that contains $u_0$. Observe that the
skeleton contains at least one rooted path of length $s$. 
If by the above procedure the same tree is produced more than once,
the corresponding skeleton has at least two different rooted paths of length $s$. 
The skeleton of a spanning tree
contains at most $n/6$ vertices, since in the tree every red triple
has two blue triples attached. There are at most $n/6-s$ vertices that can be the
end of a rooted path of length $s$. Since we are only interested in a lower
bound for the number of trees, certainly there are no more than $n/6$ rooted
paths of length $s$ in the skeleton, so each tree is produced at most $n/6$ times.

Therefore the number of spanning trees of a Steiner triple system is a least
$$\frac{4^s\,\prod_{i=0}^s (n-k-1-6i)}{n/6}. $$
This is $\Omega((n-k-1)!^{1/6})$ and since $k=o(n)$ we thus have $\Omega((n/2)!^{1/6})$ spanning trees, which by Stirling's approximation gives the statement of the theorem. 
\hspace*{\fill} \end{proof}

We now return to the question of the existence of spanning trees and
will this time present a necessary condition. Consider again a 3-graph $H=(V,\Delta)$ with underlying graph $G=(V,E)$. The
hypergraph obtained from $H$ by deleting vertices in $S\subseteq V$ is
denoted by $H-S$. This may contain hyperedges
of size $1$, $2$ or $3$. The underlying graph $G-S$ consists
of triangles for each triple, edges for each pair, and isolated
vertices for each singleton of $H-S$. A connected component of $H-S$
corresponds exactly to a connected component of the graph $G-S$.  
Let $q(H-S)$ denote the number of connected
components of $H-S$ spanning an odd number of vertices, which is also
equal to $q(G-S)$, the number of odd connected components of $G-S$. We shall use
$q(H-S)$ and $q(G-S)$ interchangeably. 
\begin{thm}\label{thm:Tutte-like} If $H=(V,\Delta)$ has a spanning
  tree then $q(H-S)\leq |S|-1$ for each
  non-empty $S\subseteq V$.
\end{thm} 
\begin{proof} Given that $H$ has a spanning tree $T$, $|V|$
  is odd. Since $q(H-S)\leq q(T-S)$ it suffices to
  prove that $q(T-S)\leq |S|-1$ for each non-empty $S\subseteq
  V$. Beginning with $|S|=1$, take $S=\{v\}$ and root $T$ at $v$. 
To each triple $abv$ of $T$ rooted at $v$ there corresponds a branch
of $T$ comprising all triples that lie on a path from
$v$ that starts with the triple $abv$.  Denote this branch by
$T_{ab}$. The 3-graph $T_{ab}$ is a tree.
Removing $v$ from $T$ creates a connected component
$T_{ab}-v$ for each $abv\in T$. Each hypergraph $T_{ab}-v$ spans an
even number of vertices. Hence the statement of the
  theorem is true for any 3-graph $H$ when $|S|=1$. Assume as induction hypothesis that the
  statement 
  is true for all 3-graphs $H$ and sets $S$ with $|S|\leq k$, where $1\leq k\leq
  |V|-1$. Suppose $|S|=k$ and take $v\in V\setminus S$. By hypothesis
  $q(H-S)\leq q(T-S)\leq k-1$ and we wish to
  prove that $q(T-S-v)\leq k$.  
  
Root $T$ at $v$ as before. 
For each $abv\in T$ define $S_{ab}=V(T_{ab})\cap S$. Possibly
$S_{ab}=\emptyset$, in which case $q(T_{ab}-S_{ab})=1$ and upon
removing $v$ we obtain one even component $T_{ab}-v$.
The non-empty sets $S_{ab}$ partition $S$. By induction hypothesis,
if $S_{ab}\neq\emptyset$ then $q(T_{ab}-S_{ab})\leq |S_{ab}|-1$ for each $abv\in T$. The vertex $v$
belongs to a unique component of $T_{ab}-S_{ab}$ for each $abv\in T$
and furthermore has degree $1$ in $T_{ab}$. Removing $v$ from $T_{ab}$
therefore creates no new components in $T_{ab}-S_{ab}$ and switches the parity of the size of the component of $T_{ab}-S_{ab}$
that contains $v$. Hence 
\begin{align*}q(T-S-v) & \leq \mathop{\sum_{ab: \; abv\in T}}_{S_{ab}\neq\emptyset}(|S_{ab}|-1)\; +
  \#\{ab: abv\in T, \, S_{ab}\neq\emptyset\}\\
 & =\mathop{\sum_{ab: \; abv\in T}}_{S_{ab}\neq\emptyset}|S_{ab}|=|S|=k.\end{align*}
This completes the inductive step. 
\hspace*{\fill} \end{proof}

 The condition of Theorem~\ref{thm:Tutte-like} although necessary for
 the existence of a spanning tree of a 3-graph is not sufficient, unlike its
 counterpart for perfect matchings of graphs (Tutte's 1-factor theorem).
The following lemma implies that if we can find a 3-graph $H$ whose
 underlying graph $G$ is  Hamiltonian then $H$ satisfies the
 conclusion of Theorem~\ref{thm:Tutte-like}.

\begin{lem}\label{lem:Hamiltonian}
Let $G=(V,E)$ be a graph with an odd number of vertices. If $G$ is Hamiltonian then $q(G-S)\leq
|S|-1$ for each non-empty $S\subseteq V$.\end{lem}
\begin{proof}  Removing $S$
  from $G$ creates at most $|S|$ connected components since this is true of the
  Hamiltonian cycle of $G$. Therefore if the condition $q(G-S)\geq |S|$ holds
  for some $S$ then there must be equality. Since $q(G-S)$ has the same
  parity as $|V|-|S|$ and $|V|$ is odd, $q(G-S)$ has parity opposite to $|S|$,
  and  hence equality is impossible.
\hspace*{\fill} \end{proof}

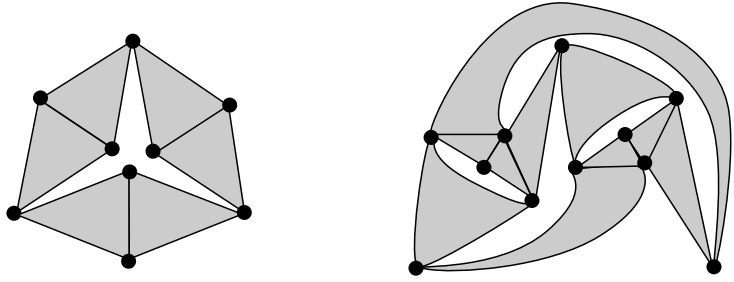
\begin{figure}[ht]
\begin{center}
\caption{Examples of 3-graphs with no spanning tree but satisfying the necessary condition
  of Theorem~\ref{thm:Tutte-like} for a spanning tree to
  exist. The 3-graphs are given by their underlying graph: each shaded triangle in the graph is a triple of
  the 3-graph.}\label{fig:Hamiltonian}

\scalebox{1} 
{
\begin{pspicture}(0,-1.9147654)(9.737807,1.9047654)
\definecolor{color1028b}{rgb}{0.8,0.8,0.8}
\psbezier[linewidth=0.02,fillstyle=solid,fillcolor=color1028b](6.6229424,0.009224757)(6.503393,-0.07064854)(5.7867265,-0.010733009)(5.667267,-0.030690774)(5.547807,-0.05064854)(6.355538,1.8947654)(7.5387983,1.7255926)(8.722059,1.5564198)(9.53087,1.0271131)(9.629339,0.2686757)(9.727807,-0.48976162)(9.459112,-1.8347654)(9.43024,-1.7869742)(9.401367,-1.739183)(9.627808,-0.3347654)(9.347807,0.2852346)(9.067807,0.9052346)(8.335195,1.3264372)(7.7578073,1.3264372)(7.18042,1.3264372)(6.8021317,1.1867329)(6.6229424,0.58799994)(6.4437532,-0.010733009)(6.742492,0.08909805)(6.6229424,0.009224757)
\psbezier[linewidth=0.02,fillstyle=solid,fillcolor=color1028b](5.596259,-1.8038379)(5.984711,-1.8765799)(7.221785,-1.791096)(7.904796,-1.3673861)(8.587808,-0.9436763)(8.547808,-0.5947654)(8.47716,-0.47629702)(8.406514,-0.35782865)(7.7278075,-0.4747654)(7.5643673,-0.47058958)(7.4009275,-0.46641377)(7.90562,-0.6945229)(7.2143593,-1.2946441)(6.523098,-1.8947654)(5.2078075,-1.731096)(5.596259,-1.8038379)
\pspolygon[linewidth=0.02,fillstyle=solid,fillcolor=color1028b](8.507808,-0.3747654)(9.407807,-1.7747654)(8.927808,0.4852346)
\pspolygon[linewidth=0.02,fillstyle=solid,fillcolor=color1028b](8.967808,0.4852346)(8.507808,-0.4547654)(8.2478075,-0.025634965)
\pspolygon[linewidth=0.02,fillstyle=solid,fillcolor=color1028b](7.607807,-0.4547654)(8.487807,-0.4347654)(8.2478075,0.0052345996)
\psbezier[linewidth=0.02,fillstyle=solid,fillcolor=color1028b](7.444163,1.1652346)(7.600518,1.2652346)(8.8552065,0.72426903)(8.92954,0.5052346)(9.003874,0.28620014)(8.917495,0.746444)(8.127808,0.2252346)(7.33812,-0.29597476)(7.6559625,-0.5747654)(7.580974,-0.4147654)(7.5059853,-0.2547654)(7.2878075,1.0652345)(7.444163,1.1652346)
\pspolygon[linewidth=0.02,fillstyle=solid,fillcolor=color1028b](7.3878074,1.1652346)(6.647807,-0.0547654)(7.047807,-0.9347654)
\pspolygon[linewidth=0.02,fillstyle=solid,fillcolor=color1028b](6.3678074,-0.4342526)(7.0078073,-0.8547654)(6.6238074,-0.0347654)
\pspolygon[linewidth=0.02,fillstyle=solid,fillcolor=color1028b](5.6678076,-0.0147654)(6.627807,-0.0147654)(6.3678074,-0.4347654)
\psbezier[linewidth=0.02,fillstyle=solid,fillcolor=color1028b](5.502442,-1.7747654)(5.6678076,-1.8147655)(7.127807,-0.8947654)(7.033494,-0.9147654)(6.939181,-0.9347654)(6.9278073,-1.0147654)(6.2296915,-0.6547654)(5.531576,-0.2947654)(5.7669287,-0.0747654)(5.6746855,-0.0547654)(5.582442,-0.0347654)(5.3370767,-1.7347654)(5.502442,-1.7747654)
\rput{-147.48312}(5.4955153,0.59782827){\pstriangle[linewidth=0.02,dimen=outer,fillstyle=solid,fillcolor=color1028b](2.8349311,-1.1829877)(1.2113466,1.3611292)}
\rput{33.060146}(0.6908296,-1.0322391){\pstriangle[linewidth=0.02,dimen=outer,fillstyle=solid,fillcolor=color1028b](2.0844529,-0.078008726)(1.2113466,1.4514928)}
\rput{146.19704}(0.7391727,-1.1678711){\pstriangle[linewidth=0.02,dimen=outer,fillstyle=solid,fillcolor=color1028b](0.54701596,-1.1973826)(1.2113466,1.4514928)}
\rput{-33.746876}(-0.15495498,0.8498361){\pstriangle[linewidth=0.02,dimen=outer,fillstyle=solid,fillcolor=color1028b](1.3234271,-.00021242457)(1.2113466,1.3611292)}
\psdots[dotsize=0.2,dotangle=2.5439634](1.7021602,1.227653)
\psdots[dotsize=0.2,dotangle=2.5439634](0.122352816,-1.0647274)
\psdots[dotsize=0.2,dotangle=2.5439634](0.47447222,0.47241664)
\psdots[dotsize=0.2,dotangle=2.5439634](1.4255234,-0.2060387)
\rput{90.00736}(-0.18361169,-1.9914763){\pstriangle[linewidth=0.02,dimen=outer,fillstyle=solid,fillcolor=color1028b](0.9038043,-1.8481152)(1.2113466,1.5211661)}
\rput{-89.78866}(3.4817326,1.3378044){\pstriangle[linewidth=0.02,dimen=outer,fillstyle=solid,fillcolor=color1028b](2.4122405,-1.8478504)(1.2113466,1.5389062)}
\psdots[dotsize=0.2,dotangle=-55.889362](3.1830263,-1.0534269)
\psdots[dotsize=0.2,dotangle=-55.889362](1.6615332,-0.513311)
\psdots[dotsize=0.2,dotangle=-55.889362](1.6440027,-1.7016736)
\psdots[dotsize=0.2,dotangle=-104.82428](1.9706051,-0.23691759)
\psdots[dotsize=0.2,dotangle=-104.82428](2.9898005,0.37433425)
\psdots[dotsize=0.2](7.4078074,1.1652346)
\psdots[dotsize=0.2](5.4678073,-1.7947654)
\psdots[dotsize=0.2](9.427808,-1.7747654)
\psdots[dotsize=0.2](6.647807,-0.0347654)
\psdots[dotsize=0.2](8.2478075,-0.0147654)
\psdots[dotsize=0.2](8.507808,-0.3947654)
\psdots[dotsize=0.2](6.3678074,-0.4547654)
\psdots[dotsize=0.2](5.6678076,-0.0547654)
\psdots[dotsize=0.2](7.0078073,-0.8947654)
\psdots[dotsize=0.2](7.587807,-0.4547654)
\psdots[dotsize=0.2](8.927808,0.4652346)
\psline[linewidth=0.02cm](8.227807,0.0052345996)(8.487807,-0.3747654)
\end{pspicture} 
}

\end{center}
\end{figure}

Figure~\ref{fig:Hamiltonian} 
gives examples of Hamiltonian graphs that underlie 
3-graphs without a spanning tree,  thereby  showing that the condition of
Theorem~\ref{thm:Tutte-like} is not sufficient. We will see in a moment why these 3-graphs
have no spanning trees.

Recall that for $H$ to have a spannning tree its underlying graph $G$ must
be connected. 
A {\em block\/} of a connected graph $G$ is a maximal 2-connected
subgraph.

\begin{prop}\label{prop:even_blocks}
Suppose the underlying graph $G$ of a 3-graph $H$
has a block that spans an even number of vertices. Then $H$ has no spanning tree.
\end{prop}
\begin{proof}
Given a spanning tree $T$ of $H$ and block $B$ of $G$, the restriction
of $T$ to the block $B$ is a tree spanning the vertices of $B$. Therefore $B$ has
an odd number of vertices.
\hspace*{\fill} \end{proof}

The parity observation behind Proposition~\ref{prop:even_blocks} can be extended to give a more general necessary
condition for the existence of a spanning tree. 

Let $H=(V,\Delta)$ be a $3$-graph.
Given subsets $V_1,\ldots, V_k$ of $V$, consider the induced
sub-3-graphs $H_i=(V_i,\Delta_i)$, where $\Delta_i=\{ abc \in \Delta:
a,b,c \in V_i\}$. Suppose moreover that the $\Delta_i$ form a partition of $\Delta$.  If
$H$ has a spanning tree $T$, this spanning tree restricted to $H_i$ yields a
spanning forest $F_i=T\cap \Delta_i$ of $H_i$. Moreover, the number of components of $F_i$ is
of the same parity as $|V_i|$.

\begin{figure}[ht]
\caption{A $3$-graph and a spanning tree (triples are triangles, the ones
belonging to the tree are shaded). Taking $V_1=\{1,2,3,4,5\}$,
$V_2=\{2,3,4,5,6,7,8\}$ and $V_3=\{6,7,8,9\}$ gives $U_1=\{2,3,4,5\}$,
$U_2=\{2,3,4,5,6,7,8\}$ and $U_3=\{6,7,8\}$. The star-partition and the
corresponding tree on $U_1\cup U_2 \cup U_3$ are shown on the right.  }\label{fig:petitex}
\begin{center}
\includegraphics[width=12cm]{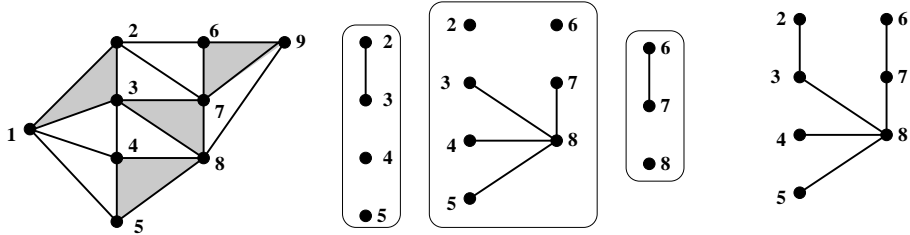}
\end{center}
\end{figure}

Let $U_i=V_i\cap (\cup_{j\neq i} V_j)$. A \emph{star-partition} of $U_i$ is a graph on $U_i$ such that each connected
component is a star, and the number of components is of the same parity as
$|V_i|$. 
For each $i$, take a star partition of $U_i$ such that vertices in a star belong to the same connected component of the underlying graph of $F_i$. Now define a graph
on $U=\bigcup U_i$ by taking the union of these star-partitions. Then this graph is a tree
on $U$. (See Figure~\ref{fig:petitex} for an example.)

The converse can be used to determine whether a $3$-graph has no spanning tree. Given
$U_i$ as above, if there is no tree on $U$ that is a union of star-partitions
then $H$ has no spanning tree. Even if there is such a tree, sometimes
the non-existence of a spanning tree can be inferred by showing that the structure
of the required star-partitions cannot be obtained from the $3$-graph.

\begin{figure}
\caption{Graph occurring as two edge-disjoint induced subgraphs of the graph on the right of Figure~\ref{fig:Hamiltonian}.}\label{fig:G_1}
\begin{center}

\scalebox{1} 
{
\begin{pspicture}(0,-1.64)(2.6336095,1.64)
\definecolor{color1028b}{rgb}{0.8,0.8,0.8}
\psbezier[linewidth=0.02,fillstyle=solid,fillcolor=color1028b](0.3991038,-1.4867061)(0.5645981,-1.5613412)(1.0916355,-1.4736325)(1.3826225,-1.0388952)(1.6736095,-0.60415804)(1.596568,-0.16616647)(1.5736095,-0.1)(1.5506508,-0.033833534)(2.4672189,-0.22304337)(1.2336094,-0.08)(0.0,0.06304337)(1.3829734,-0.34852025)(1.0884718,-0.96426016)(0.79397,-1.58)(0.23360947,-1.412071)(0.3991038,-1.4867061)
\pspolygon[linewidth=0.02,fillstyle=solid,fillcolor=color1028b](1.5936095,-0.02)(2.4936094,-1.42)(2.0136094,0.84)
\pspolygon[linewidth=0.02,fillstyle=solid,fillcolor=color1028b](2.0536094,0.84)(1.5936095,-0.1)(1.3336095,0.32913044)
\pspolygon[linewidth=0.02,fillstyle=solid,fillcolor=color1028b](0.6936095,-0.1)(1.5736095,-0.08)(1.3336095,0.36)
\psbezier[linewidth=0.02,fillstyle=solid,fillcolor=color1028b](0.529965,1.52)(0.6863204,1.62)(1.9410084,1.0790344)(2.0153422,0.86)(2.089676,0.6409655)(2.0032969,1.1012094)(1.2136095,0.58)(0.42392218,0.058790635)(0.7417648,-0.22)(0.666776,-0.06)(0.5917872,0.1)(0.37360945,1.42)(0.529965,1.52)
\psdots[dotsize=0.2,fillstyle=solid,dotstyle=o](0.49360946,1.52)
\psdots[dotsize=0.2,fillstyle=solid,dotstyle=o](0.35360947,-1.52)
\psdots[dotsize=0.2,fillstyle=solid,dotstyle=o](2.5136094,-1.42)
\psdots[dotsize=0.2](1.3336095,0.34)
\psdots[dotsize=0.2](1.5936095,-0.04)
\psdots[dotsize=0.2](0.6736095,-0.1)
\psdots[dotsize=0.2](2.0136094,0.82)
\psline[linewidth=0.02cm](1.3136095,0.36)(1.5736095,-0.02)
\end{pspicture} 
}
\end{center}
\end{figure}
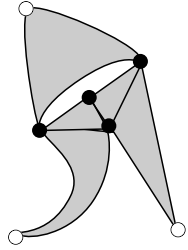

For instance, let us use this method to show that the $3$-graph on the right
of Figure~\ref{fig:Hamiltonian} has no spanning tree. Let $H_1$ and $H_2$ be
the two sub-$3$-graphs isomorphic to the $3$-graph illustrated on
Figure~\ref{fig:G_1}. Then the sets $U_1$ and $U_2$ are equal, and consist of
the three vertices that are common to both sub-$3$-graphs. Since both $H_1$
and $H_2$ have an odd number of vertices, the only possible star-partitions
for $U_1$ and $U_2$ are a star $K_{1,2}$ or three isolated vertices.  For the
union of two such star-partitions to be a tree, the only possibility is to
take one of each. Hence there must be a spanning forest of the $3$-graph in
Figure~\ref{fig:G_1} in which
the three white vertices  belong to the same component. But this forces the
spanning forest to contain three triples that form a cycle.

A particularly simple case is when  all the sets $U_i$ have 
size two.
Since the  only star-partitions of a set of two vertices are either two isolated
vertices or an edge, depending on the parity,  it is straightforward to check whether there is a tree that is a union
of star-partitions. (This argument applies to the $3$-graph on the left of
Figure~\ref{fig:Hamiltonian}.)


\section{Oriented spanning trees and 3-Pfaffian orientations}\label{sec:oriented}


\subsection{Orientations of spanning trees}\label{sec:orns_sp_trees}
An {\em orientation\/} of a finite subset of $\mathbb{N}$ is an order up to even permutation. The {\em canonical orientation\/} takes elements in the order consistent with the order $1<2<3<\cdots$ on $\mathbb{N}$.

A $(2n\!+\!1)$-cycle $(\,s(1)\;s(2)\;\cdots s(2n\!+\!1)\,)$ determines an orientation $s(1),s(2),$$\ldots,$ $s(2n\!+\!1)$ of $[2n\!+\!1]$ given by the permutation $s$. The permutation $s$ is an even or odd permutation according as it determines the same or opposite orientation of $[2n\!+\!1]$ to the canonical orientation $1<2<\cdots <2n\!+\!1$.

Suppose $H=([2n\!+\!1],\Delta)$ is a 3-graph. A triple $ijk\in\Delta$ can be assigned one of two orientations (order up to even permutation), or, what is the same thing here, a cyclic order, either $(\,i\; j\; k\,)$ or $(\, j\; i\; k\,)$. If $i<j<k$ then the canonical orientation is defined by taking the cyclic order $(i\; j\; k)$.
In other words, given $i<j$, the triple orientation $(i\;j\; k)$ is the canonical one if $k\not\in\{i+1,\ldots, j-1\}$, while if $k$ lies between $i$ and $j$ then the orientation $(\,i\; j\; k\,)$ is opposite to the canonical orientation of $ijk$.

\begin{dfn} \label{def:tree_sign} 
Suppose that we are given an orientation $\omega$ of the triples of a 3-graph
$H=([2n\!+\!1],\Delta)$. To each $t=ijk\in \Delta$ is associated a cyclic
permutation $\sigma(t,\omega)$ given by $(\,i\;j\;k\,)$ or
$(\,j\;i\;k\,)$, whichever is consistent with $\omega$. As shown by Masbaum
and Vaintrob~\cite{MV02}, if $T$ is a spanning tree of $H$, the
product $$\prod_{t\in T}\sigma(t,\omega)$$
is a $(2n+1)$-cycle $(\,s(1)\;s(2)\;\cdots \; s(2n+1)\,)$. 
The {\em
  orientation} of $T\in\mathcal{T}(H)$ associated with the triple
orientation $\omega$ is the order up to even permutation of vertices
taken in the order $s(1), s(2),
\ldots, s(2n\!+\!1)$  given by the cycle.
The {\em sign of the
 spanning tree $T$\/},
${\rm sgn}(T,\omega)$ is the sign of the permutation
 $s$.
\end{dfn}

It is also shown by Masbaum and Vaintrob~\cite{MV02} that the permutation $s$ in Definition~\ref{def:tree_sign} is determined up to conjugation by even
  permutations: not only does it not matter which of the $2n+1$ ways
  the cycle $(s(1)\; s(2)\; \cdots \; s(2n\!+\!1)\,)$ is written, but
  its sign is also independent of the order in which the factors are
  taken in the product over triples of $T$. 

We fix the notation $\omega_0$ for the canonical orientation on each triple
$ijk$ given by the cycle $(\,i\;j\;k\,)$ consistent with the natural
order $i<j<k$.

For two triple orientations $\omega_1$ and $\omega_2$ of
$H$ and $T\in\mathcal{T}(H)$ we have  
\be\label{eqn:sign_relation}{\rm sgn}(T,\omega_2)=(-1)^{\#\{t\in T: \sigma(t,\omega_1)\neq\sigma(t,\omega_2)\}}{\rm sgn}(T,\omega_1).\ee

\begin{figure}[ht]
\caption{Embeddings of two spanning trees of $K_7^{(3)}$ on vertex set
  $[7]=\{1,2,3,4,5,6,7\}$. The
  left-hand tree has oriented triples $(\, 1\; 2\; 4\,)$, $(\, 2\; 7\;
  6\,)$ and $(\, 3 \; 6\; 5\,)$. The right-hand tree has oriented
  triples $(\, 1\; 2\; 4\,)$, $(\, 3\; 7\; 4\,)$ and $(\, 4\; 6\;
  5\,)$. The linear order on $[7]$ given below each tree is obtained by taking the
  vertex labels the first time we encounter them, but other orders are
  possible by taking vertices later than at their first appearance
  (there is an even number of intermediate vertices between any two
  appearances of a given vertex). This order of appearance is then
  written as a permutation of $[7]$, whose sign gives the sign of the
  tree under the given triple orientation. }\label{fig:plane_tree_sign}
\begin{center}

\scalebox{1}
{
\begin{pspicture}(0,-3.4697502)(9.443833,3.4397502)
\pstriangle[linewidth=0.02,dimen=outer](2.4338334,1.5636874)(1.52,1.26)
\psdots[dotsize=0.2](2.4338334,2.8436873)
\psdots[dotsize=0.2](1.6738334,1.6036874)
\psdots[dotsize=0.2](3.1538334,1.6036874)
\pstriangle[linewidth=0.02,dimen=outer](1.6738334,0.40368736)(1.6,1.2)
\pstriangle[linewidth=0.02,dimen=outer](2.4738333,-0.6763126)(1.56,1.12)
\psdots[dotsize=0.2](0.87383336,0.40368736)
\psdots[dotsize=0.2](2.4338334,0.42368737)
\psdots[dotsize=0.2](1.6938334,-0.65631264)
\psdots[dotsize=0.2](3.2338333,-0.65631264)
\usefont{T1}{ptm}{m}{n}
\rput(2.098052,2.8236873){\small 1}
\usefont{T1}{ptm}{m}{n}
\rput(1.3882084,1.6236874){\small 2}
\usefont{T1}{ptm}{m}{n}
\rput(3.4102397,1.6236874){\small 4}
\usefont{T1}{ptm}{m}{n}
\rput(0.58508337,0.42368737){\small 7}
\usefont{T1}{ptm}{m}{n}
\rput(2.664927,0.44368735){\small 6}
\usefont{T1}{ptm}{m}{n}
\rput(3.4778957,-0.65631264){\small 3}
\usefont{T1}{ptm}{m}{n}
\rput(1.3794583,-0.6763126){\small 5}
\rput{179.93015}(14.58807,3.2593863){\pstriangle[linewidth=0.02,dimen=outer](7.293041,0.9835908)(1.6399988,1.3010967)}
\rput{-58.309822}(2.6210387,5.973274){\pstriangle[linewidth=0.02,dimen=outer](6.6643662,-0.047061786)(1.7060809,1.3689222)}
\rput{59.661194}(4.477149,-6.515352){\pstriangle[linewidth=0.02,dimen=outer](7.91976,-0.048398484)(1.6413996,1.3893167)}
\psdots[dotsize=0.2](6.4738336,2.2836874)
\psdots[dotsize=0.2](7.2938333,0.96368736)
\psdots[dotsize=0.2](8.113833,2.2836874)
\psdots[dotsize=0.2](5.6538334,0.98368734)
\psdots[dotsize=0.2](6.5138335,-0.43631265)
\psdots[dotsize=0.2](8.093833,-0.41631263)
\psdots[dotsize=0.2](8.893833,0.96368736)
\usefont{T1}{ptm}{m}{n}
\rput(6.2482085,2.2836874){\small 2}
\usefont{T1}{ptm}{m}{n}
\rput(5.4049273,0.9436874){\small 6}
\usefont{T1}{ptm}{m}{n}
\rput(6.2594585,-0.45631263){\small 5}
\usefont{T1}{ptm}{m}{n}
\rput(8.297895,-0.45631263){\small 3}
\usefont{T1}{ptm}{m}{n}
\rput(9.105083,0.9436874){\small 7}
\usefont{T1}{ptm}{m}{n}
\rput(8.298052,2.3036873){\small 1}
\usefont{T1}{ptm}{m}{n}
\rput(7.29024,0.64368737){\small 4}
\psbezier[linewidth=0.02,linestyle=dotted,dotsep=0.09cm,arrowsize=0.1cm 2.0,arrowlength=1.5,arrowinset=0.6]{->}(1.8538333,2.7236874)(1.6338334,2.0636873)(0.0,0.9342101)(0.39383337,0.26368737)(0.78766674,-0.4068354)(2.0149386,0.49062827)(2.0338333,0.26368737)(2.0527282,0.03674647)(2.0538332,0.0036873643)(1.8338333,-0.16339783)(1.6138333,-0.33048305)(1.0938333,-0.52985525)(1.1538334,-0.83631265)(1.2138333,-1.14277)(3.7144477,-1.3381523)(3.7138333,-0.7563126)(3.713219,-0.174473)(3.2138333,0.28368735)(2.7938333,0.6841358)(2.3738334,1.0845842)(1.7538334,1.2693375)(2.1938334,1.4105932)(2.6338334,1.5518488)(3.6151714,1.1163814)(3.6338334,1.5636873)(3.6524951,2.0109935)(3.2138333,2.3236873)(2.7338333,2.7836874)
\usefont{T1}{ptm}{m}{n}
\rput(2.152896,-2.2163126){\small 1276534}
\psbezier[linewidth=0.02,linestyle=dotted,dotsep=0.09cm,arrowsize=0.1cm 2.0,arrowlength=1.5,arrowinset=0.6]{->}(6.873833,1.2036873)(6.7138333,1.1236874)(5.333833,1.4036874)(5.1138334,1.2236874)(4.893833,1.0436873)(6.2338333,-1.0163126)(6.5738335,-0.83631265)(6.913833,-0.65631264)(7.0138335,0.44368735)(7.1138334,0.44368735)(7.2138333,0.44368735)(7.4938335,0.48368737)(7.6338334,0.16368736)(7.7738333,-0.15631263)(7.815108,-0.8311378)(8.153833,-0.7563126)(8.4925585,-0.68148744)(9.433833,0.6236874)(9.333834,1.1636873)(9.233833,1.7036873)(8.273833,0.96368736)(7.6738334,1.1436874)(7.0738335,1.3236873)(8.926917,1.7376246)(8.353833,2.5836873)(7.7807503,3.4297502)(5.8052564,2.9461474)(6.0338335,2.3636873)(6.26241,1.7812272)(6.18276,2.145078)(6.7938333,1.4436873)
\usefont{T1}{ptm}{m}{n}
\rput(7.169146,-2.1963127){\small 4653712}
\usefont{T1}{ptm}{m}{n}
\rput(2.1764896,-2.7363126){\small ( 3 7 4 6 )}
\usefont{T1}{ptm}{m}{n}
\rput(2.2571146,-3.2563126){\small negative sign}
\usefont{T1}{ptm}{m}{n}
\rput(7.0664897,-2.6963127){\small ( 1 4 3 5 7 2 6 )}
\usefont{T1}{ptm}{m}{n}
\rput(7.118677,-3.2163126){\small positive sign}
\end{pspicture} 
}

\end{center}
\end{figure}
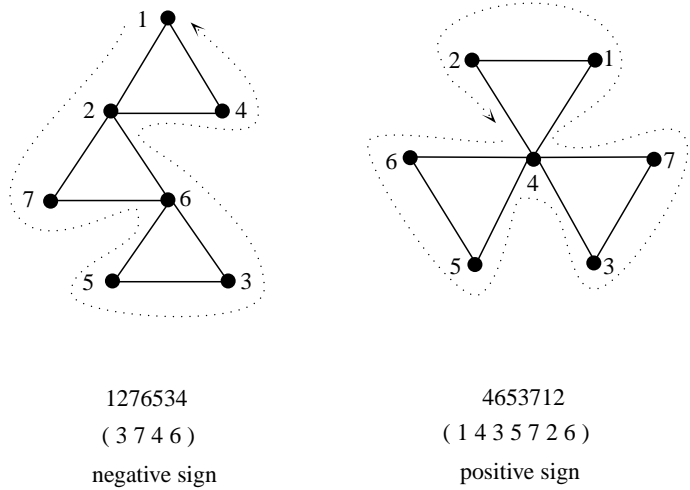

A convenient way to calculate the sign of a spanning tree is as
follows and illustrated in
Figure~\ref{fig:plane_tree_sign}. Given a spanning tree $T$ of $H$ and a triple orientation $\omega$, embed the underlying graph
of $T$ in the plane so that boundaries of the interior faces are the triples
of $T$ and so that the vertices of a triple $ijk$  appear in anticlockwise
order consistent with the triple orientation $\omega$. 
Starting at an arbitrary vertex, tour the tree in an anticlockwise
sense, reading off a cyclic string of $3n$ vertex labels. Remove
repeated vertex labels until a cyclic string of length $2n+1$ remains,
equal to $(\,s(1)\;s(2)\;\cdots\;s(2n+1)\,)$ for some permutation $s$
of $[2n+1]$. Then the sign of $s$ as a permutation is equal to ${\rm
  sgn}(T,\omega)$. 

For a given orientation $\omega$ of triples of $H=([2n\!+\!1],\Delta)$, let $\mathcal{T}^+(H)=\{T\in\mathcal{T}(H):{\rm
  sgn}(T,\omega)=+1\}$ and
$\mathcal{T}^-(H)=\{T\in\mathcal{T}(H):{\rm
  sgn}(T,\omega)=-1\}$. These sets will be denoted by $\mathcal{T}^+$
and $\mathcal{T}^-$ respectively when $H=K_{2n+1}^{(3)}$ is complete. 

 For $S\subseteq\Delta$ define $(\,i\;j\,)S$ to be the
 set obtained from $S$ by switching  $i$ and $j$ in
 triples containing either of these two vertices.
If $i$
 and $j$
 have the property that $\{t-\{i\}:t\in\Delta, i\in
 t\}=\{t-\{j\}:t\in\Delta, j\in t\}$ then this action set-stabilizes
 $\Delta$. Furthermore, if $T$ is a spanning tree of $H$ then in this
 situation $(\,i\; j\,)T$ is also a spanning tree of $H$. Under the
 canonical orientation, the sign of
 $(\, i\; j\,)T$ is related to the sign of $T$ in a particularly
 straightforward way when $j=i\!+\!1$:  

\begin{lem}\label{lem: sign reversing}
Let $H=([2n\!+\!1],\Delta)$ be a 3-graph with canonical orientation of
its triples. Suppose that
$i\in[2n]$ has the property that $\{t-\{i\}:t\in \Delta, i\in
t\}=\{t-\{i\!+\!1\}: t\in \Delta, i\!+\!1\in t\}$. Then $(i\;
i\!+\!1)T$ is a spanning tree of $H$ with opposite sign to that of $T$ if
$\{i,i\!+\!1\}$ is not contained in any triple of $T$, while $(i\;
i\!+\!1)T$ has the same sign as $T$ if some triple of $T$ contains $\{i,i\!+\!1\}$.
\end{lem}
\begin{proof} 
Let us start with a fixed embedding of $T$ in the plane, as described
  above. Interchanging the labels $i$ and $i+1$ in the embedding gives an
  embedding of $(i\ i+1)T$ in which all triples appear in anticlockwise order
 if and only if there is no triple containing both $i$ and $i+1$. If this is the case,
  when touring the embedding of $(i\ i+1)T$  in anticlockwise order, we obtain
  the same permutation of the vertices as when touring $T$, except that
  elements $i$ and $i+1$ are transposed. Hence clearly the sign of $(i\ i+1)T$
  is opposite to that of $T$. 
If $T$ contains a (necessarily unique) triple $\{i,i+1,j\}$, consider the orientation  $\omega$  that
agrees with $\omega_0$ everywhere except in the triple $\{i,i+1,j\}$. Then
${\rm sgn}((i\ i+1)T,\omega_0)=-{\rm sgn}((i\ i+1)T,\omega)$
by equation~\eqref{eqn:sign_relation} and ${\rm sgn}((i\
i+1)T,\omega)=-{\rm sgn}(T,\omega_0)$ by the same argument about touring
the embedding as before.\hspace*{\fill} \end{proof}

The involution $T\mapsto (\,i\;i\!+\!1\,)T$ of
  Lemma~\ref{lem: sign reversing} specializes to the sign-reversing
  involution on perfect matchings on $[2n]$ of \cite[Lemma
    2.1]{Stem90} when applied to 3-graphs in which every triple
  contains the vertex $2n+1$ (where spanning trees of the 3-graph
  correspond precisely to perfect matchings on $[2n]$).
 
\begin{thm}\label{thm:positive_negative}
The distribution of positive and negative spanning trees of $K_{2n+1}^{(3)}$ under the canonical orientation is given by
$$|\mathcal{T}^+|-|\mathcal{T}^-|=(2n+1)^{n-1}.$$
\end{thm}
\begin{proof}  
Let $\mathcal{T}_i$ denote the set of trees that have a triple
containing $\{2i\!-\!1, 2i\}$.
 By Lemma~\ref{lem: sign reversing}, 
the involution $\tau_i:T\mapsto(\,2i\!-\!1\;2i\,)T$ reverses the sign
of trees in $\mathcal{T}\setminus \mathcal{T}_i$. 
If $T\in\mathcal{T}_j$ then $(\,2i\!-\!1\; 2i\,)T\in\mathcal{T}_j$, since the pairs $\{\{2i-1,2i\}:i\in[n]\}$ are pairwise disjoint.  So for each $j\in[n]$ the restriction of $\tau_i$ to
$\mathcal{T}_j$ is a map $\mathcal{T}_j\rightarrow\mathcal{T}_j$
reversing the sign of trees in $\mathcal{T}_j\setminus
\mathcal{T}_i$. (On $\mathcal{T}_i$
itself $\tau_i$ fixes the sign of every tree.) 
Hence $$|\mathcal{T}^+|-|\mathcal{T}^-|=\sum_{T\in\mathcal{T}}{\rm sgn}(T,\omega_0)=\sum_{T\in
  \mathcal{T}_1\cap\mathcal{T}_2\cap\cdots\cap\mathcal{T}_n}{\rm
  sgn}(T,\omega_0).$$
A tree belonging to $\bigcap_{i\in[n]}\mathcal{T}_i$ has set of
  triples equal to $\{\{2i-1,2i, f(i)\}:i\!\in\![n]\}$ for some
  function $f:[n]\rightarrow[2n+1],$ (with $f(i)=2n+1$
  for at least one value of $i$). The canonical
  orientation of a triple $\{2i\!-\!1,2i, f(i)\}$ is $(\,2i\!-\!1\; 2i\; f(i)\,)$, no
  matter whether $f(i)>2i$ or $f(i)<2i\!-\!1$. To show that a tree 
 $T\in\bigcap_{i\in[n]}\mathcal{T}_i$ is positively oriented
    under the canonical orientation $\omega_0$,
we embed it $T$ in the plane so
    that the vertices of a triple appear in anticlockwise order consistent
    with the orientation $\omega_0$. Traversing the tree anticlockwise starting   at vertex $2n+1$, the vertices appear, for some permutation $\pi$ of $[n]$, in the order
    $2n+1,2\pi(1)\!-\!1, 2\pi(1),\, 2\pi(2)\!-\!1, 2\pi(2),\,\ldots,\, 2\pi(n)\!-\!1, 2\pi(n)$. This is an even permutation of $1,2,3,4,\ldots, 2n-1,
    2n, 2n+1$.

 To evaluate $|\bigcap_{i\in[n]}\mathcal{T}_i|$, use the ``Pr\"{u}fer code'' described in the proof of
   Theorem~\ref{thm:pruefer}, in which the perfect
 matching $M$ is fixed equal to $\{\{2i\!-\!1, 2i\}:i\in[n]\}$. Trees in
 $\bigcap_{i\in[n]}\mathcal{T}_i$ are in bijective
   correspondence with sequences $\gamma\in[2n+1]^{n-1}$.
\hspace*{\fill} \end{proof}

\subsection{Tree generating polynomials}\label{sec:tree_generating_polynomials}
Let $y=(y_t:t\in \Delta)$ be a set of commuting indeterminates indexed
by triples of the sub-3-graph $H=([2n\!+\!1],\Delta)$ of
$K_{2n+1}^{(3)}$. (Here we depart from Masbaum and Vaintrob~\cite{MV02},
but follow for example Caracciolo et el.~\cite{CMSS08}, by indexing
the indeterminates by
triples rather than oriented triples. In other words,
$y_{ijk}=y_{\{i,j,k\}}=y_{jik}$, and so on.)
The \emph{tree generating polynomial} of $H$ is defined by
$$\mathcal{P}(H,y)=\sum_{T\in\mathcal{T}(H)}\prod_{t\in T}y_t.$$
The \emph{signed tree generating polynomial} associated with an orientation $\omega$ of the edges  is defined by
$$\mathcal{P}^\omega(H,y)=\sum_{T\in\mathcal{T}(H)}{\rm sgn}(T,\omega)\prod_{t\in T}y_t.$$
   
By equation~\eqref{eqn:sign_relation} in the previous subsection, the polynomial $\mathcal{P}^{\omega}(H,y)$ is related to
$\mathcal{P}^{\omega_0}(H,y)$ by substituting $-y_t$ for $y_t$ for triples
$t$ on which $\omega$ is opposite to $\omega_0$.

Define the antisymmetric $(2n\!+\!1)\times (2n\!+\!1)$ matrix $\Lambda$ with $(i,j)$ entry given by 
$$\Lambda_{i,j}=\sum_{k\neq i,j}\epsilon_{i,j,k}y_{ijk},$$
where $\epsilon_{i,j,k}=+1$ if $(\, i\; j\; k\, )$ is a cyclic permutation of
$i<j<k$, $\epsilon_{i,j,k}=-1$ if $(\, i\; j\; k\, )$ is opposite to this canonical orientation, and $\epsilon_{i,j,k}=0$ if two of the indices are equal.
 Let $\Lambda^{(k)}$ denote the matrix obtained from $\Lambda$ by deleting row
 $k$ and column $k$. The following is the Pfaffian matrix-tree theorem of
 Masbaum and Vaintrob.
\begin{thm}\label{thm: MV} {\rm \cite{MV02}}
For any $k\in[2n\!+\!1]$ the signed tree polynomial associated with the canonical orientation $\omega_0$ is given by 
$$\mathcal{P}^{\omega_0}(K_{2n+1}^{(3)},y)=(-1)^{k-1}{\rm Pf}(\Lambda^{(k)}).$$
\end{thm}
An orientation $\omega$ of the triples of $K_{2n+1}^{(3)}$ restricted to $\Delta\subseteq\binom{[2n+1]}{3})$ gives an orientation of the sub-3-graph $H=([2n+1],\Delta)$; the signed tree polynomial $\mathcal{P}^{\omega}(H,y)$ is obtained from $\mathcal{P}^{\omega}(K_{2n+1}^{(3)},y)$ upon setting $y_t=0$ if $t\not\in \Delta$.

\begin{dfn}
An orientation $\omega$ of the triples of a 3-graph
$H$ is {\em 3-Pfaffian\/} if ${\rm sgn}(T,\omega)$ is
constant for $T\in\mathcal{T}(H)$. 
A 3-graph is said to be \emph{3-Pfaffian} if there exists some 3-Pfaffian
orientation of its triples.
\end{dfn}
See Subsection~\ref{ssec:minimal} for some examples of 3-Pfaffian and
non-3-Pfaffian 3-graphs.

For a 3-Pfaffian orientation $\omega$ of $H$,
$\mathcal{P}^\omega(H,y)=\pm \mathcal{P}(H,y)$; in particular, in this
case by Theorem~\ref{thm: MV} the number of spanning trees of $H$ will
be computable in polynomial time by the evaluation of
$\mathcal{P}^\omega(H;1)$ (setting $y_t=1$ for each
$t\in\Delta$). To evaluate $\mathcal{P}^{\omega}(H,1)$ from
$\mathcal{P}^{\omega_0}(H,y)$ set $y_t=+1$
for triples $t$ on which $\omega$ is the same as $\omega_0$ and
$y_t=-1$ when $\omega$ is opposite to $\omega_0$ on $t$

The correspondence from Lemma~\ref{lem.sp trees pm + function} between spanning trees and perfect matchings $M$ of $G-\{2n+1\}$
together with an 
``apex-choosing'' function $f:M\rightarrow [2n]$ is used by Hirschman and
Reiner~\cite{HR04}  to prove a useful alternative formulation  of the Masbaum--Vaintrob theorem.

\begin{thm} \label{thm: HR} {\rm \cite{HR04}} For a 3-graph $H=([2n+1],\Delta)$,
$$\mathcal{P}^{\omega_0}(H,y)=\mathop{\sum_{\mbox{\rm \tiny perfect matchings $M$ of $[2n]$}}}_{f:M\rightarrow[2n+1]}{\rm sgn}(M)\mathop{\prod_{ij\in M}}_{i<j}\epsilon_{i,j,f(ij)}y_{ijf(ij)},$$
where ${\rm sgn}(M)$ is the sign of the perfect matching $M$, given by 
$${\rm sgn}(M)=(-1)^{\mbox{\rm \tiny cross}(M)},$$
$${\rm cross}(M)=\#\{i<j<k<l:\{i,k\},\{j,l\}\in M\}.$$
\end{thm}

An orientation of the edges of a graph
  $G=([2n],E)$ is Pfaffian if for all perfect matchings $M$ of $G$ the
quantity $${\rm sgn}(M) \cdot (-1)^{\#\{i<j:j\longrightarrow i\}}$$ is constant,
  where ${\rm sgn}(M)$ is defined as in the previous theorem and
  $j\longrightarrow i$ denotes an oriented edge with $j$ directed towards $i$.  

As a straightforward application, Theorem~\ref{thm: HR} yields a simple criterion for
3-graph $H$ to be 3-Pfaffian when $H$ has the property that all its triples
contain a common vertex --- in the terminology of
Section~\ref{sec:suspensions} below, that is to say when the 3-graph $H$ is the 1-suspension of an ordinary graph.

\begin{thm}\label{thm:1-susp}
Let $G$ be a graph on vertex set $[2n]$ and edge set $E\subseteq\binom{[2n]}{2}$. Let $H$ be the $3$-graph on vertex set $[2n+1]$ and triple set $\Delta=\{ij(2n\!+\!1):ij\in E\}$.

Then $H$ has a 3-Pfaffian orientation if and only if $G$ has a
Pfaffian orientation: if edge $ij$ has orientation $i\longrightarrow j$ in the
Pfaffian orientation of $G$ then the triple orientation given by
$(\,i\; j\; 2n\!+\!1\,)$ defines a 3-Pfaffian orientation of $H$, and conversely.
\end{thm}
\begin{proof}
 Spanning trees of $H$ are in one-one correspondence with perfect
matchings of $G$. As described in the statement of the theorem, orientations
of the edges of $G$ are also in one-one correspondence with orientations of
the triples of $H$.  Let us denote by $\stackrel{\omega}{\longrightarrow}$ the
orientation of $G$ corresponding to orientation $\omega$ of $H$, i.e., $i\stackrel{\omega}{\longrightarrow} j$ if and only if the triple $ij(2n+1)$
is oriented $(\,i\; j\; 2n\!+\!1\,)$ in $\omega$.  Recall that
$\mathcal{P}^{\omega}(H,y)=\mathcal{P}^{\omega_0}(H,\bar{y}), $
where $\bar{y}_t$ equals or is opposite to $y_t$ depending on whether
$\omega$ and $\omega_0$ agree or not in $t$. 

In the expansion of Theorem~\ref{thm: HR} the term $\epsilon_{i,j,2n+1}$ is
constant equal to $1$ if $i<j$, therefore,
\begin{align*}\mathcal{P}^{\omega_0}(H,\bar{y}) & = \mathop{\sum_{\mbox{\rm \tiny perfect
      matchings $M$ of $[2n]$}}} {\rm sgn}(M)\mathop{\prod_{ij\in
      M}}_{i<j}\bar{y}_{ij(2n+1)}\\
& =\mathop{\sum_{\mbox{\rm \tiny perfect
      matchings $M$ of $[2n]$}}} {\rm sgn}(M) (-1)^{\# \{i<j: j\rightarrow i\}} .\end{align*}

It is clear then that $\omega$ is a $3$-Pfaffian orientation of $H$ if and
only if $\stackrel{\omega}{\longrightarrow}$ is a Pfaffian orientation of $G$.\hspace*{\fill} \end{proof}

Recall that $\mathcal{T}(H)$ denotes the set of spanning trees of a
3-graph $H$ and that 
$\mathcal{T}(H\backslash abc)=\{T\in\mathcal{T}(H): \, abc\not\in T\}$
and there is a bijection between
$\mathcal{T}(H/abc)$ and $\{T\in\mathcal{T}(H):\, abc\in T\}.$
If $abc$ is in no spanning tree of $H$ then
$\mathcal{T}(H)=\mathcal{T}(H\backslash abc)$. If $abc$ is in every
spanning tree of $H$ then contracting the triple $abc$ defines a
bijection from $\mathcal{T}(H)$ to $\mathcal{T}(H/abc)$.

 \begin{lem}\label{lem.delete contract}
If a triple $abc$ occurs in no spanning tree of 
$H$ then $H$ is 3-Pfaffian if and only if $H\backslash abc$ is
3-Pfaffian. Similarly,  if $abc$ occurs in every spanning tree of
$H$ then $H$ is 3-Pfaffian if and only if $H/abc$ is 3-Pfaffian.
\end{lem}
\begin{proof} The only thing to prove is that contracting a triple
  $abc$ either preserves the sign of all spanning trees
  $\{T\in\mathcal{T}(H):abc\in T\}$ or reverses all their signs.

Let $V=[2n+1]$ and $abc$ the triple to be contracted. By labelling the vertices suitably we may assume
$\{b,c\}=\{2n,2n+1\}$: the property that all spanning trees have the
same sign is unaffected by a permutation of vertex labels.
 
 Embed a given tree $T\in\mathcal{T}(H)$ in the plane so that the
 orientation of triples corresponds to the anticlockwise order of its vertices. The anticlockwise appearance of vertices
around $T$ up to cyclic permutation takes the form $AaBbCc$, where $A,B,C$ are each an even length
sequence of vertices, and $A\cup B\cup C\cup\{a,b,c\}$ is a partition of
$[2n+1]$.  Upon contracting $abc$ to a vertex with label $a$, a
spanning tree of $H/abc$ on vertex set $[2n-1]$ is obtained with
vertices around the tree appearing in the order $AaBC$ up to cyclic permutation. The parity of
$AaBbCc$ as a permutation of $[2n+1]$ is the same as the parity of $AaBCbc$
since $B$ has even length. Since $b,c$ are greater than all the other
vertex labels the parity of $AaBCbc$ is equal to that of $AaBC$ plus
that of $bc$. Hence all spanning trees of $H$ have their sign multiplied by
the sign of $bc$ as a permutation of $\{2n,2n+1\}$ when contracting
the triple $abc$.
\hspace*{\fill}\end{proof}

\begin{dfn} A 3-graph $H=(V,\Delta)$ is {\em minimally
non-3-Pfaffian\/} with respect to triple deletion and contraction if
$H$ is non-3-Pfaffian and there is no $t\in\Delta$ such that $H\backslash t$
or $H/t$ is non-3-Pfaffian. (A 3-graph that has no spanning trees is
vacuously 3-Pfaffian.)
\end{dfn}

Lemma~\ref{lem.delete contract} implies that in a minimal
non-3-Pfaffian 3-graph (with respect to triple deletion and contraction) 
each triple occurs in at least one
spanning tree and no triple occurs in all spanning trees.

Since the property of being 3-Pfaffian is preserved by deletion
and contraction, if a 3-graph $H$ after deletion and contraction of
triples gives a non-3-Pfaffian graph then $H$ must be
non-3-Pfaffian. This is the same as restricting attention to spanning
trees of $H$ that contain a given subset of triples (those that are contracted) and disjoint
from another subset of triples (those deleted).  More generally, if some
subset of the class $\mathcal{T}(H)$ of all spanning trees of $H$ can be shown to be impossible to make all the same
sign then the same is true of the whole class $\mathcal{T}(H)$, i.e.,
$H$ is non-3-Pfaffian.

\subsection{Complexity results for orientations}\label{sec:complexity}

As observed in the previous subsection, $\mathcal{P}(H;1)=|\mathcal{T}(H)|$, and
$\mathcal{P}^\omega(H;1)=|\mathcal{T}^+(H)|-|\mathcal{T}^-(H)|$. 
Consider the distribution
of $\mathcal{P}^\omega(H;1)$ when $\omega$ is a triple orientation
chosen uniformly at random (u.a.r.). 
By equation~\eqref{eqn:sign_relation}, if we let
$y_t$ take values in $\{-1,+1\}$ u.a.r. for $t\in\Delta$
then the random variable $\mathcal{P}^{\omega_0}(H;y)$ is equal to the
random variable
$\mathcal{P}^\omega(H;1)$
under an orientation $\omega$ of triples taken u.a.r. The following lemma is analogous to the well-known
result~\cite{LP86} that the expected value of the determinant of the skew
adjacency matrix of a graph $G$ (under all possible orientations of its
edges) is equal to the number of perfect matchings of $G$. 

\begin{lem}\label{prop:variance}
Suppose $H=([2n\!+\!1],\Delta)$ is  a 3-graph with a fixed
orientation $\omega$ of its triples. For each $t\in\Delta$ let $y_{t}$ take values in $\{-1,+1\}$
independently uniformly at random, while $y_{t}=0$ when $t\not\in\Delta$. Then
$$\mathbb{E}[\mathcal{P}^{\omega}(H;y)]=0,$$
$$\mathbb{E}[\mathcal{P}^{\omega}(H;y)^2]=|\mathcal{T}(H)|.$$ 
\end{lem} 
\begin{proof} 
The random variables $y_t$ for $t\in\Delta$ are independent, each with
expected value $\mathbb{E}(y_t)=0$. For $S\subseteq\Delta$ let
$y_S=\prod_{s\in S}y_s$. Then $\mathbb{E}(y_T)=0$ for each spanning
tree $T$ and
$$\mathbb{E}\left[\sum_{T\in\mathcal{T}(H)}{\rm sgn}(T,\omega)y_T\right]=\sum_{T\in\mathcal{T}(H)}{\rm
  sgn}(T,\omega)\mathbb{E}(y_T)=0.$$
Also,
$$\mathbb{E}\left[\big(\sum_{T\in\mathcal{T}(H)}{\rm
    sgn}(T,\omega)y_T\big)^2\right]=\sum_{S,T\in\mathcal{T}(H)}{\rm sgn}(S,\omega){\rm
  sgn}(T,\omega)\mathbb{E}(y_{S\bigtriangleup T}),$$
where $\mathbb{E}(y_{S}y_T)=\mathbb{E}(y_{S\bigtriangleup T})$,
for if $t\in S\cap T$ then $y_t^2=1$.
Since $\mathbb{E}(y_{S\bigtriangleup T})=0$ unless $S\bigtriangleup
T=\emptyset$ in which case $\mathbb{E}(y_\emptyset)=1$ this yields
$$\mathbb{E}\left[\big(\sum_{T\in\mathcal{T}(H)}{\rm
    sgn}(T,\omega)y_T\big)^2\right]=\sum_{T\in\mathcal{T}(H)}{\rm
  sgn}(T,\omega)^2=|\mathcal{T}(H)|.$$
\hspace*{\fill} \end{proof}

Whereas counting (unsigned) spanning trees of 3-graphs (evaluating
$|\mathcal{T}(H)|$) is \#{\sf P}-complete in general, the problem
of evaluating $|\mathcal{T}^+(H)|-|\mathcal{T}^-(H)|$ under
any given triple orientation turns
out to be polynomial time by Theorem~\ref{thm: MV} above, as it is the
evaluation of the Pfaffian of a polynomial-size matrix with integer
entries (each bounded in absolute value by $2n\!-\!1$).

\begin{cor}\label{cor:pos_neg_orn} A 3-graph $H$ has a spanning tree, i.e.,
  $P(H;1)=|\mathcal{T}(H)|\neq 0$, if and only if there is some 
  triple orientation $\omega$ of $H$ such that
$P^{\omega}(H;1)=|\mathcal{T}^+(H)|-|\mathcal{T}^-(H)|\neq 0$. 
\end{cor} 
\begin{proof} Clearly $|\mathcal{T}^+(H)|-|\mathcal{T}^-(H)|\neq 0$
implies the existence of a spanning tree. 
By Lemma~\ref{prop:variance} the variance of
$|\mathcal{T}^+(H)|-|\mathcal{T}^-(H)|$ is positive if and only if
$|\mathcal{T}(H)|\neq \emptyset$. 
\hspace*{\fill} \end{proof}

If there is a point $y$ such that $\mathcal{P}^\omega(H;y)\neq 0$ then
$H$ has a spanning tree.  
Caracciolo et al.~\cite{CMSS08} give an algorithm that runs in expected polynomial time for deciding the
existence of a spanning tree. Since the polynomial $\mathcal{P}^\omega(H;y)$ has $|\Delta|\leq\binom{2n+1}{3}$ variables and total degree $n$ the problem of
deciding if it is non-zero can be solved in expected polynomial time
by evaluating it at random points in a field $\mathbb{F}_q$ of
sufficiently large order $q\geq 2n$. 

We turn from the problem of deciding if there is a triple
orientation for which the difference between positively and negatively
oriented spanning trees is non-zero to
the problem of whether there is
a 3-Pfaffian orientation (for which all spanning trees have the same
sign). The former problem is polynomial time by
Corollary~\ref{cor:pos_neg_orn} and the fact that deciding if there is
a spanning tree is polynomial time. We do not know whether the problem of whether a 3-graph is
3-Pfaffian can be solved in polynomial time. (It is also unknown whether the problem of deciding if a graph is Pfaffian can be solved in polynomial time.)

However, a similar method of proof to that of Vazarani and Yannakakis~\cite{VY89} for Pfaffian orientations of graphs shows that the problem of deciding
the existence of a
$3$-Pfaffian orientation is in {\sf co-NP}.
The main idea is to write a system of linear equations whose solutions are the
$3$-Pfaffian orientations of a $3$-graph. We start by explaining this
construction, which will be also be used later in the paper. 

Let $H=(V,\Delta)$ be a $3$-graph and let $\mathcal{T}(H)$ be its collection
of spanning trees. 
Consider the triple--spanning tree incidence matrix $M\in\mathbb{F}_2^{\mathcal{T}(H)\times \Delta}$ with $(T,t)$ entry equal to $1$ if $t\in T$ and $0$ otherwise.
The rows of $M$ are the indicator vectors in $\mathbb{F}_2^\Delta$ of
the triple sets of spanning trees $T\in \mathcal{T}(H)$.
The columns  of $M$ are the indicator vectors in $\mathbb{F}_2^{\mathcal{T}(H)}$ of those trees that change orientation when the orientation of triple $t$ is reversed (i.e., those trees containing $t$).
Let $\mathbf{c}\in\mathbb{F}_2^{\mathcal{T}(H)}$ denote the indicator vector
of tree orientations under the canonical orientation of edges, that is, for
$T\in \mathcal{T}(H)$ the $T$-component of $\mathbf{c}$ is $0$ if
$\mathrm{sgn}(T,\omega_0)=1$ and is $1$ if $\mathrm{sgn}(T,\omega_0)=-1$. There is some
orientation of edges that leads to all trees $T\in\mathcal{T}(H)$ having the
same sign if and only if either of the equations
$$M\mathbf{x}=\mathbf{c}, \qquad M\mathbf{x}=\mathbf{c}+\mathbf{1}$$
has a solution ($\mathbf{x}$ is the indicator vector of a subset of
triples which when flipped in orientation change the tree orientations
to have all positive signs or all negative signs, respectively).

\begin{thm}\label{thm: existence coNP}
The problem of deciding whether a $3$-graph has a $3$-Pfaffian orientation is in {\sf co-NP}.
\end{thm} 

\begin{proof}
According to the previous discussion, deciding whether a $3$-graph $H$ has a
$3$-Pfaffian orientation is equivalent to finding a solution of either of the equations
$M\mathbf{x}=\mathbf{c}, \qquad M\mathbf{x}=\mathbf{c}+\mathbf{1}.$
The length of the vector $\mathbf{c}$ and the number of rows of $M$ is
$|\mathcal{T}(H)|$, typically exponential in the number of vertices, say $2n+1$. However,
the rank of $M$ is polynomial on $n$, since $M$ has $O(n^3)$ columns (one for
each triple). If the system is inconsistent, basic linear algebra implies that
there is a subset of rows of $M$,
say $M'$, such that $\mathrm{rank}(M')<\mathrm{rank}(M'|\mathbf{c'})$, where
$\mathbf{c'}$ is the restriction of $\mathbf{c}$ to the rows of $M'$. Since
$\mathrm{rank}(M'|\mathbf{c'})$ cannot be more than the number of columns of
$M$ plus one, there is a polynomial time certificate that the equation
$M\mathbf{x}=\mathbf{c}$ is inconsistent. Doing the same for the equation
$M\mathbf{x}=\mathbf{c+1}$, one can verify in polynomial time that $H$ has no
$3$-Pfaffian orientation.
\hspace*{\fill} \end{proof}

In the next two sections we consider two special families of 3-graphs
for which we can say more about the existence of 3-Pfaffian orientations.

\section{Suspensions of graphs and 3-Pfaffian orientations}\label{sec:suspensions}

\begin{dfn}
Let $G=(V,E)$ be a graph and $U$ a finite set disjoint from $V$. Then the
{\em suspension} of $G$ from $U$ is the 3-graph $G^U=(U\cup V,\Delta)$ with set of triples 
$\Delta=\{iju:ij\in E, u\in U\}.$
If $U$ has $k$ elements then $G^U$ is called a {\em $k$-suspension} of
$G$. (All $k$-suspensions of $G$ are isomorphic.)  
\end{dfn}

For the 3-graph $G^U$ to have a spanning tree it is necessary that $G$ has no
isolated vertices and for $|U|$ to have opposite parity to $|V|$. 

In this section we characterize those graphs $G$ whose $k$-suspension
has a 3-Pfaffian orientation.

The case of 1-suspensions has already been dealt with at the end of Subsection~\ref{sec:tree_generating_polynomials}. 
A spanning tree $T$ of a 1-suspension $G^{\{u\}}$ consists 
 of triples $\{iju:ij\in M\},$
where $M$ is a perfect matching of $G$. In particular, if $G$ has no
perfect matching then $G^{\{u\}}$ has no spanning trees. 
There is a bijective correspondence between
orientations of triples of the 1-suspension $G^{\{u\}}$ and orientations of
edges of $G$. If an edge $ij$ of $G$ is oriented $i\longrightarrow j$,
then the triple $iju$ has orientation given by the cyclic order $(i\;j\;u)$. 
By Theorem~\ref{thm:1-susp}, the $1$-suspension $G^{\{u\}}$ has a 3-Pfaffian
  orientation if and only if $G$ has a Pfaffian orientation.

For $3$-suspensions and upwards, there is no orientation that makes all
spanning trees have the same sign, unless of course there is no spanning
tree. 
\begin{thm}
Let $G$ be a graph and $u,v,w\not\in V(G)$. If the $3$-suspension
$G^{\{u,v,w\}}$ has a spanning tree then it has no 3-Pfaffian
orientation. For $k\geq 4$, the analogous result holds for the  $k$-suspension of $G$.
\end{thm}
\begin{proof}

Up to symmetry in $u,v,w$, a spanning tree of $G^{\{u,v,w\}}$ takes one of the following two
forms:
\begin{enumerate}[(i)]
\item $\{uxa, vxb, wxc\}\cup \{iju:ij\in M_1\}\cup\{ijv:ij\in
  M_2\}\cup\{ijw:ij\in M_3\}$, where $M_1, M_2, M_3$ are matchings
  together spanning $G-\{a,b,c,x\}$, or 
\item $\{uxa, vxb, vyc, wyd\}\cup \{iju:ij\in M_1\}\cup\{ijv:ij\in
  M_2\}\cup\{ijw:ij\in M_3\}$,  where $M_1, M_2, M_3$ are matchings
  together spanning $G-\{a,b,c,d,x,y\}$.
\end{enumerate}

Recall that a $3$-Pfaffian $3$-graph remains $3$-Pfaffian after the deletion and
contraction of triples; therefore, it is enough to show that after suitable
contractions and deletions of $G^{\{u,v,w\}}$ we obtain a $3$-graph that is not $3$-Pfaffian.

Suppose first that $G^{\{u,v,w\}}$ has a spanning tree as in case (i). Fix the
matchings $M_1, M_2, M_3$. Let $G_1$
be the graph on $\{a,b,c,x\}$ with edges $\{ax,bx,cx\}$ and let
$H_1=G_1^{\{u,v,w\}}$. The $3$-graph $H_1$ is obtained from $G^{\{u,v,w\}}$ by
contracting the triples $\{iju:ij\in M_1\}\cup\{ijv:ij\in
  M_2\}\cup\{ijw:ij\in M_3\}$ and deleting the triples that remain and do not
  belong to $G_1^{\{u,v,w\}}$.

The spanning trees of $G_1^{\{u,v,w\}}$ are $S_\pi=\{a\pi(u)x,
b\pi(v)x,c\pi(w) x\}$, where $\pi$ ranges over the six permutations of $\{u,v,w\}$.
 Consider the order $a<b<c<u<v<w<x$ on the vertices of $G_1^{\{u,v,w\}}$
 and let $\omega_0$ be the canonical orientation associated with this order. It is easy to see that the  sign of
 $S_\pi$ under this orientation is the sign of the permutation
$$a \, \pi(u)\, b\, \pi(v) \, c \pi(w)\, x. $$
Therefore, three of the trees $S_{\pi}$ are positive and three are
negative. Since each of the triples appears in exactly two of the trees,
changing the orientation of any of the triples keeps the parity of the
number of positive and negative  trees. Thus, there is no orientation that
makes all six trees the same sign, as needed.

If $G^{\{u,v,w\}}$ has a spanning tree as in case (ii), one argues
analogously by considering the $3$-graph $G_2^{\{u,v,w\}}$, where $G_2$ has
  edges $\{ax,bx,cy,dy\}$.

Finally, that a $k$-suspension $G^{U}$ of a graph is non-3-Pfaffian for $k\geq 3$
follows by a similar argument by permuting $3$ of the vertices in $U$
while fixing the rest of the spanning tree. 

\hspace*{\fill} \end{proof}

The case of 2-suspensions is the richer one and occupies the rest of this
section. The main result (Theorem~\ref{thm:2-susp_3-Pfaffian}) is a characterization of
those graphs for which the 2-suspension is $3$-Pfaffian in terms of forbidden
subgraphs. This is similar in spirit to the result by Little~\cite{Little75} characterizing 
Pfaffian bipartite graphs as those without an even subdivision of $K_{3,3}$
with a perfect matching in the complement. Before stating and proving our
characterization, we need  another result akin to the theory of Pfaffian
orientations. Recall that all perfect matchings of a graph have the same sign in
a given orientation if
and only if any cycle of even length whose complement has a perfect matching
has an odd number of edges in each direction. Our goal is to establish a
similar characterization of 3-Pfaffian orientations of 2-suspensions (Theorem~\ref{thm:carac_3pfaffian_or}). 
For this we need first to describe spanning trees of $H=G^{\{u,v\}}$ and their unions in terms
of the graph $G$, so that conditions arise for all trees to have the same sign
under a given orientation. From now on the $2$-suspension is denoted
by $G^{u,v}$.
 
Spanning trees of $G^{u,v}$ correspond to matchings $M_u$ and $M_v$ of $G=(V,E)$
with the property that $V(M_u)\cap V(M_v)=\{i\}$ for some single
vertex $i$ and $V(M_u)\cup V(M_v)=V$. We call the subgraph $M_u\cup
M_v$ a {\em quasi-perfect matching} of $G$. A quasi-perfect matching
consists of a collection of independent edges and a single path on two edges, which together partition the vertices of $G$.
A spanning tree $T$ of $G^{u,v}$ has 
triple set 
$$\{iju:ij\in M_u\}\cup\{ijv:ij\in M_v\},$$
for some quasi-perfect matching $M_u\cup M_v$ of $G$.
 
Having described the spanning trees of $G^{u,v}$, we 
calculate their sign under a given orientation of
$G^{u,v}$. Triple orientations of $G^{u,v}$ can be obtained
from orientations of $G$, and vice versa. To do this we assume that the vertex set of $G$
is $[2n+1]$ and that the vertex set of $G^{u,v}$ is ordered
$1<2<\cdots< 2n+1<u<v$. Recall that the canonical orientation of a triple
$ijk$ takes $i,j,k$ in linear order up to even permutation. 

Suppose we are given an orientation of triples of $G^{u,v}$. This
orientation of triples is determined by its sign relative to the
canonical orientation, positive or negative according as it has the same or
opposite sense respectively. For each edge $ij$ of $G$ there are two triples of $G^{u,v}$,
namely $iju$ and $ijv$.  We call the edge $ij$ \emph{agreeing\/} if the
orientations of $iju$ and $ijv$ are both equal or both contrary to the  canonical orientation; we call it
\emph{opposite\/} otherwise. The $u$-orientation of $G$ is the orientation of
 $G$ that orients the edge $ij$ with $i<j$ as $i\longrightarrow j$ if  the
 triple $iju$ has orientation $(i\, j\, u)$ and $j\longrightarrow i$ otherwise. 
 The notation $i\stackrel{u}{\longrightarrow}j$ means that the
 $u$-orientation of the 
 edge $ij$ is $i\longrightarrow j$.
 Analogous notions are defined with
 respect to $v$.

\begin{lem}\label{lem:sign_tree}
Let $T$ be a tree of $G^{u,v}$ with associated quasi-perfect matching $M_u\cup
M_v$. Let $xy$ and $yz$ be the edges of the path of length 2 in
$M_u\cup M_v$, with $xy\in
M_u$, and let  $i_1j_1,\ldots,i_{n\!-\!1}j_{n\!-\!1}$ be the
other edges of $M_u\cup M_v$,
written such that $\uor{i_\ell}{j_\ell}$ or $\vor{i_\ell}{j_\ell}$ depending on whether
$i_\ell j_\ell$ belongs to $M_u$ or to $M_v$. 

Then the sign of $T$ is the product of the sign of the permutation
$$\left(\begin{array}{ccccccccccc} 1&2& 3& 4 &\ldots & 2(n-1) & 2n-1&  2n &2n+1 &u&v\\ i_1 &j_1 &i_2&
    j_2 &\ldots & j_{n\!-\!1} & u& x &y &z &v   \end{array} \right)  $$
and $(-1)^{\alpha_u(xy)+\alpha_v(yz)}$ where $\alpha_u(xy)=0$ if $\uor{x}{y}$ and $\alpha_u(xy)=1$ otherwise, and
similarly for $\alpha_v(yz)$.  
\end{lem}

\begin{proof}
The formula follows from the definition of the sign of a tree in terms of the traversal
of a planar embedding together with the fact that switching the orientation of
one edge switches the sign of the tree.

More concretely, if we draw the planar embedding of the tree assuming that the
triples $xyu$ and $yzv$ are oriented $(x\; y\; u)$ and $(y \; z\; v)$,
respectively, and then we traverse the tree in anticlockwise sense, the permutation
whose sign we need is
$$\left(\begin{array}{cccccccccccccc} 1&2& \ldots & 2\ell & 2\ell+1 & & &\ldots && &\ldots&
     2n+1 &u&v\\ i_1 &j_1 &\ldots & j_{\ell} & u& x &y &z &v
    &i_{\ell+1}&\ldots&j_{n\!-\!2} &i_{n\!-\!1} &j_{n\!-\!1}   \end{array} \right),$$
where we assume that the edges $i_1j_1,\ldots,i_\ell j_\ell$ are the ones
in $M_u$. This permutation and the one in the statement differ in an even
number of transpositions hence they have the same sign. The term
$(-1)^{\alpha_u(xy)+\alpha_v(yz)}$ collects the change of sign if the triples
$xyu$ and $yzv$ are oriented differently.
\hspace*{\fill} \end{proof}

The following lemma is an easy consequence, but it will be used often in the sequel. Given a subgraph $G'$ of $G$, its
{\em complement\/} is the graph induced by the vertices not in $G'$,
i.e., $G-V(G')$.

\begin{lem}\label{lem:edge_agrees}
If an edge $ij$ of $G$ is such that its complement has a quasi-perfect
matching and $G^{u,v}$ has a 3-Pfaffian orientation, then $ij$ is agreeing (in that orientation). 
\end{lem}

\begin{proof}
Let $Q$ be the quasi-perfect matching in the complement of $ij$. There are many
spanning trees of $G^{u,v}$ that correspond to the quasi-perfect matching
$Q\cup \{ij\}$ of $G$. Of all these trees, let $T_1$ and $T_2$ be two of them
such that they only differ in that $T_1$ contains the triple $iju$ and
$T_2$ contains the triple $ijv$. By Lemma~\ref{lem:sign_tree}, any
orientation that gives the same sign to $T_1$ and $T_2$ must agree on $ij$.
\hspace*{\fill} \end{proof}

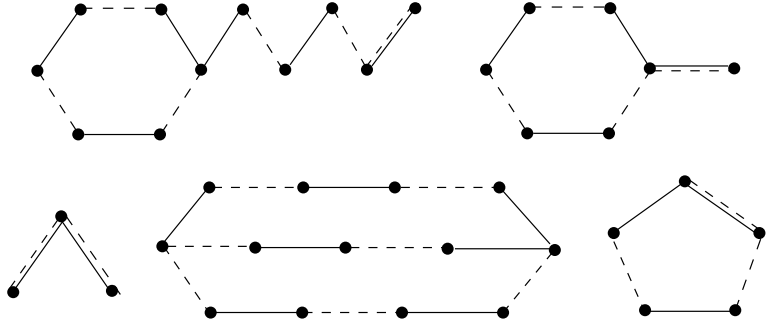
\begin{figure}[ht]
\caption{Some examples of connected components of type (H) and (T)
  given in Lemma~\ref{lem:unions_qpm}. For clarity, one
quasi-perfect matching is depicted by a solid line, the other by a
dashed line.}\label{fig:compsqpm}
\begin{center}

\scalebox{0.8} 
{
\begin{pspicture}(0,-2.65)(12.62,2.65)
\psline[linewidth=0.02cm,linestyle=dashed,dash=0.16cm 0.16cm](4.92,-2.53)(6.6,-2.53)
\psline[linewidth=0.02cm](4.14,-1.45)(5.68,-1.45)
\psdots[dotsize=0.2](1.22,2.51)
\psdots[dotsize=0.2](3.92,2.51)
\psdots[dotsize=0.2](2.56,2.51)
\psdots[dotsize=0.2](4.62,1.51)
\psdots[dotsize=0.2](5.4,2.53)
\psline[linewidth=0.02cm](4.7,1.59)(5.36,2.51)
\psline[linewidth=0.02cm,linestyle=dashed,dash=0.16cm 0.16cm](5.42,2.51)(5.96,1.55)
\psdots[dotsize=0.2](5.98,1.51)
\psline[linewidth=0.02cm,linestyle=dashed,dash=0.16cm 0.16cm](5.96,1.55)(6.74,2.59)
\psline[linewidth=0.02cm](6.06,1.57)(6.84,2.61)
\psdots[dotsize=0.2](6.78,2.53)
\psline[linewidth=0.02cm,linestyle=dashed,dash=0.16cm 0.16cm](1.24,2.53)(2.46,2.53)
\psline[linewidth=0.02cm](1.18,2.47)(0.52,1.53)
\psline[linewidth=0.02cm](2.62,2.47)(3.2,1.55)
\psdots[dotsize=0.2](0.5,1.49)
\psdots[dotsize=0.2](1.18,0.43)
\psdots[dotsize=0.2](3.22,1.51)
\psline[linewidth=0.02cm,linestyle=dashed,dash=0.16cm 0.16cm](0.52,1.47)(1.14,0.51)
\psline[linewidth=0.02cm](1.28,0.43)(2.52,0.43)
\psdots[dotsize=0.2](2.54,0.43)
\psline[linewidth=0.02cm,linestyle=dashed,dash=0.16cm 0.16cm](3.24,1.49)(2.56,0.45)
\psline[linewidth=0.02cm](3.26,1.55)(3.88,2.51)
\psline[linewidth=0.02cm,linestyle=dashed,dash=0.16cm 0.16cm](3.96,2.49)(4.58,1.51)
\psdots[dotsize=0.2](8.68,2.53)
\psdots[dotsize=0.2](10.02,2.53)
\psdots[dotsize=0.2](12.08,1.53)
\psline[linewidth=0.02cm,linestyle=dashed,dash=0.16cm 0.16cm](8.7,2.55)(9.92,2.55)
\psline[linewidth=0.02cm](8.64,2.49)(7.98,1.55)
\psline[linewidth=0.02cm](10.08,2.49)(10.66,1.57)
\psdots[dotsize=0.2](7.96,1.51)
\psdots[dotsize=0.2](8.64,0.45)
\psdots[dotsize=0.2](10.68,1.53)
\psline[linewidth=0.02cm,linestyle=dashed,dash=0.16cm 0.16cm](7.98,1.49)(8.6,0.53)
\psline[linewidth=0.02cm](8.74,0.45)(9.98,0.45)
\psdots[dotsize=0.2](10.0,0.45)
\psline[linewidth=0.02cm,linestyle=dashed,dash=0.16cm 0.16cm](10.7,1.51)(10.02,0.47)
\psline[linewidth=0.02cm](10.72,1.57)(11.98,1.57)
\psline[linewidth=0.02cm,linestyle=dashed,dash=0.16cm 0.16cm](10.72,1.49)(12.14,1.49)
\psdots[dotsize=0.2](0.9,-0.93)
\psline[linewidth=0.02cm,linestyle=dashed,dash=0.16cm 0.16cm](0.84,-0.99)(0.04,-2.13)
\psline[linewidth=0.02cm,linestyle=dashed,dash=0.16cm 0.16cm](0.98,-0.93)(1.88,-2.23)
\psline[linewidth=0.02cm](0.94,-0.99)(0.12,-2.15)
\psline[linewidth=0.02cm](0.88,-0.97)(1.68,-2.13)
\psdots[dotsize=0.2](0.1,-2.19)
\psdots[dotsize=0.2](1.74,-2.17)
\psdots[dotsize=0.2](2.58,-1.43)
\psline[linewidth=0.02cm](2.58,-1.39)(3.36,-0.43)
\psline[linewidth=0.02cm,linestyle=dashed,dash=0.16cm 0.16cm](2.66,-1.43)(3.98,-1.43)
\psline[linewidth=0.02cm,linestyle=dashed,dash=0.16cm 0.16cm](2.6,-1.43)(3.38,-2.57)
\psdots[dotsize=0.2](3.36,-0.45)
\psdots[dotsize=0.2](3.38,-2.53)
\psdots[dotsize=0.2](4.12,-1.45)
\psline[linewidth=0.02cm,linestyle=dashed,dash=0.16cm 0.16cm](3.42,-0.45)(4.76,-0.45)
\psdots[dotsize=0.2](4.92,-0.45)
\psline[linewidth=0.02cm](4.94,-0.45)(6.44,-0.45)
\psline[linewidth=0.02cm](3.44,-2.53)(4.86,-2.53)
\psdots[dotsize=0.2](4.9,-2.53)
\psdots[dotsize=0.2](5.62,-1.45)
\psdots[dotsize=0.2](6.44,-0.45)
\psline[linewidth=0.02cm,linestyle=dashed,dash=0.16cm 0.16cm](5.7,-1.45)(7.36,-1.45)
\psline[linewidth=0.02cm,linestyle=dashed,dash=0.16cm 0.16cm](6.54,-0.45)(8.2,-0.45)
\psdots[dotsize=0.2](8.18,-0.45)
\psdots[dotsize=0.2](7.32,-1.47)
\psdots[dotsize=0.2](6.56,-2.53)
\psline[linewidth=0.02cm](8.22,-0.49)(9.02,-1.39)
\psline[linewidth=0.02cm](7.44,-1.47)(9.12,-1.47)
\psline[linewidth=0.02cm](6.6,-2.53)(8.2,-2.53)
\psdots[dotsize=0.2](8.24,-2.53)
\psdots[dotsize=0.2](9.1,-1.49)
\psline[linewidth=0.02cm,linestyle=dashed,dash=0.16cm 0.16cm](9.06,-1.53)(8.24,-2.53)
\psdots[dotsize=0.2](12.5,-1.21)
\psdots[dotsize=0.2](11.26,-0.35)
\psline[linewidth=0.02cm,linestyle=dashed,dash=0.16cm 0.16cm](12.52,-1.31)(12.1,-2.49)
\psline[linewidth=0.02cm,linestyle=dashed,dash=0.16cm 0.16cm](11.34,-0.33)(12.5,-1.15)
\psline[linewidth=0.02cm](11.26,-0.35)(10.06,-1.21)
\psline[linewidth=0.02cm](11.24,-0.37)(12.44,-1.21)
\psline[linewidth=0.02cm,linestyle=dashed,dash=0.16cm 0.16cm](10.08,-1.35)(10.6,-2.53)
\psdots[dotsize=0.2](10.08,-1.21)
\psdots[dotsize=0.2](10.6,-2.49)
\psdots[dotsize=0.2](12.1,-2.49)
\psline[linewidth=0.02cm](10.64,-2.51)(12.06,-2.51)
\end{pspicture}

}

\end{center}
\end{figure}

In order to compare the sign of two spanning trees, we look at their union,
which we next describe in terms of the associated quasi-perfect matchings.
For the rest of this section, it will be convenient to consider that  an edge is a cycle of length two.
See Figure~\ref{fig:compsqpm} for an illustration of the statement of the following lemma.

\begin{lem}\label{lem:unions_qpm}
Let $Q_1$ and $Q_2$ be two quasi-perfect matchings. Then the connected
components of $Q_1\cup Q_2$ are of the following types.
\begin{itemize}
\item[(C)] A cycle of even length. 
\item[(H)] Two edge-disjoint cycles  with a path (possibly empty)  with ends in the cycles.
\item[(T)] Three internally vertex-disjoint paths  having common endpoints
  (including a cycle of odd length as a degenerate case).

\end{itemize}
Moreover, all connected components except one are of type (C), and the
component of type (H) or (T) has an odd number of vertices.
\end{lem}

\begin{proof}
Let $p_1$ and $p_2$ be the 2-paths in $Q_1$ and $Q_2$.
If a connected component of $Q_1\cup Q_2$ contains no edges from $p_1$
or $p_2$ then we are in
case (C), since the component will result from the union of two matchings. So
we focus on the component containing $p_1$. Colour the
edges $Q_1$ blue and the edges of $Q_2$ red. An edge in $Q_1\cap Q_2$ edge is both red
and blue.  Let $xy$, $yz$ be the blue edges of $p_1$. If $xy$ is also red, $xy$ is a
cycle of length two. Otherwise $x$ must be incident to some red edge
$xx_1$, since $Q_2$ is a quasi-perfect matching of $G$. Similarly  $x_1$ is incident to a blue edge $x_1x_2$, and so on, until some
vertex $x_k$ is repeated. (There may be a choice between two red edges
along the way if the path $p_2$ is encountered when forming this
cycle. If this is the case then an arbitrary choice of red edge is made.) The edge $x_{k-1}x_k$ must be red, since every
vertex is incident to at least one blue edge, and the only vertex
incident to two blue edges is $y$. If $x_k=y$, we continue to explore the
connected component from a red edge incident with  $z$ and eventually another cycle is
closed. Otherwise we continue the component from $y$. In both cases a second
cycle is closed;  as before the last edge added must be red, hence the vertex
at which the second cycle is closed is the middle vertex of the path $p_2$. 
At this point all vertices in the component are incident with one edge
of each colour, except for one or two vertices which are adjacent to two edges
of the same colour and one or two of the other. So these are all the
edges of $Q_1\cup Q_2$ in this
component. Note that in particular both paths $p_1$ and $p_2$ are always in
the same component.
We are in case (T) or (H) according to whether the two cycles meet in an edge
or not. Note that a particular example of (T) consists of a cycle of odd
length, considering that one of the edges is a cycle of length 2.

The claim on the number of vertices follows from the fact that the total
number of vertices is odd and components of type (C) have an even number of vertices.
\hspace*{\fill} \end{proof}

Note that if the graph $G$ is bipartite the paths and cycles in the
statement of Lemma~\ref{lem:unions_qpm} are all of even length. 

By inspecting the
$u$- and $v$-orientations of the edges in a component of the
union of two quasi-perfect matchings of $G$, we are able to
characterize 3-Pfaffian orientations of $G^{u,v}$ in terms of their behaviour on even
cycles and some other small subgraphs of $G$. To reach this
characterization we require some further lemmas.

Given a graph with an orientation of its edges, a cycle of even
length is said to be \emph{oddly oriented} if when traversing
it cyclically we encounter an odd number of edges oriented forward (and
hence an odd number oriented backwards). 
By allowing cycles of length two, the next lemma is a generalization of Lemma~\ref{lem:edge_agrees}. (A cycle of length two is always oddly oriented.)

\begin{lem}\label{lem:cycle_odd}
Let $C$ be a cycle of even length in $G$ such that its complement contains a
quasi-perfect matching $Q$. If a given orientation of $G^{u,v}$ is 3-Pfaffian, then all the edges of $C$ are agreeing and the cycle
is oddly oriented (with respect to the given orientation).
\end{lem}

\begin{proof}
That all the edges of $C$ are agreeing follows from
Lemma~\ref{lem:edge_agrees}, so we focus on the second claim. 
Let $a_1,b_1,\ldots,a_k,b_k$ be the vertices of $C$ in cyclic
order.  
Construct a (partial) tree $T_Q$ of $G^ {u,v}$ from $Q$ in the following way:
if $Q=M\cup N$ for some matchings $M$ and $N$ of $G-C$, let $T_Q=\{ iju : ij\in
M\} \cup \{ijv:v\in N\}$.

 Now let $T_1$ be the
  tree having as triples $T_Q$ plus the triples $\{a_ib_iu: 1\leq i \leq k\}$ and let $T_2$ be the tree whose triples
  are $T_Q$ together with  $\{a_ib_{i-1}u: 2\leq i \leq k\}\cup \{a_1b_ku\}$  . Assume
  the edges in $C$ are oriented cyclically, that is,   $a_i\longrightarrow b_i$ and $ b_i \longrightarrow a_{i+1}$. Then by
  Lemma~\ref{lem:sign_tree} the trees $T_1$ and $T_2$ have opposite
  signs. Hence if an orientation gives both of them the same sign, an odd number of the
  edges in $C$ need to be reversed.
\hspace*{\fill} \end{proof}

\begin{lem}\label{lem:mountain}
Let $xy$ and $yz$ be two edges of $G$ such that the complement of
their union
contains a perfect matching. In any 3-Pfaffian orientation of $G^{u,v}$, one
of the two edges is agreeing and the other is opposite.
\end{lem}

\begin{proof}
  Let $M$ denote the perfect matching, and let $T_M$ be the collection of
  triples obtained by adding $u$ to the edges of $M$. Let $T_1=T_M\cup
  \{xyu,yzv\}$ and $T_2=T_M \cup \{xyv,yzu\}$.  The
  conclusion follows again by comparing the expressions for the signs of $T_1$ and
  $T_2$ given in Lemma~\ref{lem:sign_tree}.
\hspace*{\fill} \end{proof}

\begin{cor}\label{cor:nop6}
If $G^{u,v}$ is 3-Pfaffian, then $G$ does not contain a path
of length 6 whose complement has a perfect matching.
\end{cor}

\begin{proof}
Suppose for a contradiction that $a_1a_2\ldots a_7$ is path of length 6 in
$G$. Take a 3-Pfaffian orientation of $G$. The complement of the 
edge $a_3a_4$ contains a quasi-perfect matching, hence this edge is
agreeing. Similarly, $a_4a_5$ is also agreeing. But
Lemma~\ref{lem:mountain} implies that only one of $a_3a_4$ and
$a_4a_5$ can be agreeing. 
\hspace*{\fill} \end{proof}

The following lemma describes how a 3-Pfaffian orientation behaves in a path
of length 4.

\begin{lem}\label{lem:path4}
Let $x_1x_2x_3x_4x_5$ be a path of length $4$ in $G$ whose complement has a
perfect matching. In any 3-Pfaffian
orientation, the edges $x_1x_2$ and $x_4x_5$ are agreeing and the
other two are opposite. Moreover,  $\vor{x_2}{x_3}$ if and only if
$\vor{x_4}{x_3}$, and analogously for the $u$-orientation.
\end{lem}

\begin{proof}
Which edges are agreeing and which ones are opposite follows from
Lemmas~\ref{lem:edge_agrees} and~\ref{lem:mountain}. Now we proceed as in the
proof of Lemma~\ref{lem:mountain}. Let $M$ denote the perfect matching in the
complement of the path, and let $T_M$ be the collection of
triples obtained by adding $u$ to the edges of $M$. Let $T_1=T_M\cup
\{x_1x_2u,x_2x_3v,x_4x_5u\}$ and $T_2=T_M \cup
\{x_1x_2u,x_3x_4v,x_4x_5u\}$. 
The conclusion follows again by comparing the expressions for the sign of $T_1$ and
$T_2$ given in Lemma~\ref{lem:sign_tree}.
\hspace*{\fill} \end{proof}

The necessary conditions for an orientation to be 3-Pfaffian given in the
previous lemmas turn out to be sufficient. Recall that an edge is considered
to be a cycle of length two.

\begin{thm}\label{thm:carac_3pfaffian_or}
The following are equivalent for an orientation of $G^{u,v}$.
\begin{enumerate}[(i)]
\item The orientation is 3-Pfaffian.
\item With respect to this orientation,\begin{enumerate}[(a)]
  \item if $C$ is an even cycle of $G$ whose complement has a quasi-perfect
    matching,  all its edges are agreeing and $C$ is oddly oriented;
\item if $xyz$ is a path of length 2 in  $G$ whose complement has a
  perfect matching, one of the edges is agreeing and the other is opposite; 
\item if $x_1x_2x_3x_4x_5$ is a path of length 4 in $G$ whose complement
  has a perfect matching, then $\uor{x_2}{x_3}$ if and only if
  $\uor{x_4}{x_3}$, and analogously for the $v$-orientation.
 \end{enumerate}
  \end{enumerate}
\end{thm}

\begin{proof}
The implication (i)$\Rightarrow$(ii) follows from
Lemmas~\ref{lem:cycle_odd},~\ref{lem:mountain} and ~\ref{lem:path4}. For the
converse, let $T_1$ and $T_2$ be two spanning trees of $G^{u,v}$. We need to
prove that they get the same sign if the orientation satisfies the conditions
in (ii).

We first show that certain subgraphs cannot appear in $G$ if there is an
orientation satisfying (ii). A $P_6$ is a path with $6$ edges  and $K_{2,3}^{-}$
denotes the graph $K_{2,3}$ with one edge removed.

\medskip

{\sc Claim 1.} If $G$ has an orientation satisfying (ii), then $G$
has no  subgraph isomorphic to an odd cycle, a $P_6$ or a $K_{2,3}^{-}$
whose complement contains a perfect matching.

\medskip

\textit{Proof of Claim 1.} Let $C$ be an odd cycle in $G$. Since $C$ is not
2-edge-colourable, there are two consecutive edges of $G$ that are either both
opposite or both agreeing. If $C$ has a perfect matching in the complement,
condition (ii).(b) applied to this pair of edges yields a contradiction.

That $G$ contains no path of length $6$ with a perfect matching in the
complement follows from the same argument as in Corollary~\ref{cor:nop6} using
(ii).(a) and (ii).(b).

Finally, suppose $x_1x_2,x_2x_3,x_3x_4,x_4x_5,x_5x_2$ are the edges of a
$K_{2,3}^{-}$ with a perfect matching in the complement. By (ii).(a), the
edges $x_3x_4$ and $x_4x_5$ are agreeing, since each contains a quasi-perfect
matching in the complement. But by (ii).(b) one of them must be opposite.
\hspace*{\fill} $\Box$

We next see which are the connected components of $Q_1\cup
Q_2$, where $Q_1$ and $Q_2$ are the quasi-perfect matchings associated to
$T_1$ and $T_2$. 

\medskip

{\sc Claim 2.} The connected components of $Q_1\cup Q_2$ are cycles of even
length and a path of length 2 or 4.

\medskip

\textit{Proof of Claim 2.}

It will be used throughout the proof that a graph whose connected components
are even cycles and paths of odd length contains a perfect matching.

Lemma~\ref{lem:unions_qpm} gives the three types of components that can
arise. They are all even cycles (including edges), except for one of the
components that is of type (H) or (T). Let us take a connected component $D$ of
type (H). It consists of two  cycles joined by a possibly empty path. Due to the restriction on the order of $D$, only the following two
combinations can arise: the two cycles have the same parity and the path has
even length, or the two cycles have different parity and the path has odd
length. In this last case, it is easy to see that $D$ contains a spanning
subgraph consisting of an  odd cycle and a perfect matching, which together
with a perfect matching in the type (C) components contradicts Claim 1.
Hence $D$ consists of
two even cycles joined by a path of even length. If one of the cycles has
length six or more, Claim 1 is again contradicted by finding a  $P_6$ with a
perfect matching in the complement. Finally, if one of the cycles has length
$4$, it is easy to find a $K_{2,3}^{-}$ with a perfect matching in the
complement. We have thus reached the conclusion that a component of type $H$ consists of
two cycles of length 2 joined by a path of even length, that is, the
component is a path of even length, and this path can only have length 2 or 4
by Claim 1.

Next we look at possible components of type (T), that is, three paths with
common endpoints. Since the total number of vertices is odd, there are two paths of the
same parity, which together form a cycle of even length, and the other
path has necessarily even length. Reasoning as in the
preceding paragraphs, we conclude that the cycle has length 2. Thus, in fact
the (T) component is an odd cycle, which is impossible, so there are is no
component of type (T).
\hspace*{\fill} $\Box$

The only thing left is to conclude that both trees have the same sign. This follows
from Lemma~\ref{lem:sign_tree}. 

More concretely, suppose that the component of
type (H) in the union of
$Q_1\cup Q_2$ is a path of length 2, say $a\, b\, c$. It could be that  both trees
contain the triple $abu$, or that both contain the triple $abv$,  or that one of them, say $T_1$, contains $abu$
and the other one $abv$. Let us focus first on the latter case. 
To compute the sign of $T_1$,  we  compute first the sign of the permutation $\pi\   u\  a\  b\  c\  v$ , where $\pi$ are the entries that
correspond to  vertices that do not belong to the path of length $2$. 
To get the sign of $T_1$ we may
need to modify the sign according to the orientation of the path of length 2.
  The
corresponding permutation for $T_2$ can be split similarly as $\pi' u\
c\ b\ a\ v$. The two permutations $\pi$ and $\pi'$ differ in an even number of
transpositions, since all even cycles in $Q_1\cup Q_2$ are oddly
oriented. The permutations $a\ b\ c$ and $c\ b\ a$ have clearly opposite
signs, so $T_1$ and $T_2$ have the same sign if and only if
$\alpha_u(ab)+\alpha_v(bc)+\alpha_u(cb)+\alpha_v(ba)$ is odd, and this is
implied by (ii).(b). The case that  both $T_1$ and $T_2$ contain $abu$
(or $abv$) is simpler and dealt with in the same way.   

We now suppose that the component of
type (H) in the union of
$Q_1\cup Q_2$ is a path of length 4, say $a\, b\, c\, d\, e$, with $a\, b\, c$ being the
path of length $2$ in $Q_1$ and $c\, d\, e$ that in $Q_2$. 
To compute the sign of $T_1$,  we  compute first the sign of the permutation $\pi\ d\  e\  u\  a\  b\  c\  v$ , where $\pi$ are the entries that
correspond to  vertices that do not belong to the path of length $4$. We
assume that $T_1$ contains triples $abu$ and $bcv$; this is no restriction
since the tree $(ac)T_1$ has the same sign as $T_1$ by the conclusion of the previous paragraph.
To get the sign of $T_1$ we may
need to modify the sign according to the orientation of the path of length 2
and to that of edge $de$.  The
corresponding permutation for $T_2$ can be split similarly as $\pi'\ a\ b\ u\
c\ d\ e\ v$. The two permutations $\pi$ and $\pi'$ differ in an even number of
transpositions, since all even cycles in $Q_1\cup Q_2$ are oddly
oriented. Note also that $d\  e\  u\  a\  b\  c\  v$ and $a\ b\ u\
c\ d\ e\ v$ have the same sign. Hence both trees have the same sign if and
only if $\alpha_u(ab)+\alpha_v(bc)+\alpha_u(cd)+\alpha_v(de)
+\alpha_u(de)+\alpha_v(ab)$ is even. The edges $ab$ and $de$ are agreeing by
(ii).(a), therefore we only need to worry about
$\alpha_v(bc)+\alpha_u(cd)$. That this is even follows by combining the fact
that both $bc$ and $cd$ are opposite and the condition in (ii).(c).
\hspace*{\fill} \end{proof}

The conditions of Theorem~\ref{thm:carac_3pfaffian_or} for an orientation of a
2-suspension to be 3-Pfaffian are quite restrictive and suggest that there
are few of them. This is confirmed by the
following characterization by forbidden subgraphs.
As usual $C_\ell$ denotes the
cycle with $\ell$ edges.

\begin{thm} \label{thm:2-susp_3-Pfaffian} Let $G$ be a graph and $u,v\not\in V(G)$. Then the
  2-suspension $G^{u,v}$ has a 3-Pfaffian orientation if and only
  if \begin{enumerate}[(i)]
\item the graph $G-\{i\}$ is Pfaffian for each
  vertex $i$,
\item $G$ has no subgraph isomorphic to $C_3, C_5, P_6$ or
  $K_{2,3}^-$ whose complement has  a perfect matching.
\end{enumerate}

\end{thm}

\begin{proof}
If $G^{u,v}$ has a $3$-Pfaffian orientation, Claim 1 in the proof of
Theorem~\ref{thm:carac_3pfaffian_or} shows that $G$ contains no subgraph
isomorphic to an odd cycle, a $P_6$ or a $K_{2,3}^-$ whose complement contains
a perfect matching, hence (ii) holds. (Observe that excluding $P_6$  
automatically excludes all odd cycles of length at least $7$.)
To show (i) holds, consider the
$u$-orientation of $G$ corresponding to the $3$-Pfaffian orientation of $G^{u,v}$. Let $C$ be a cycle of
even length $\ell\geq 4$ in $G-\{i\}$ whose complement has a perfect matching
$M$. We need to show that $C$ is oddly oriented with respect to the
$u$-orientation. If vertex
$i$ is adjacent to some vertex in $C$, then $G$ would contain a copy of
 $K_{2,3}$ or of $P_6$ with a perfect matching in the complement, so we 
conclude that $i$ is only adjacent to vertices covered by the perfect
 matching $M$. Hence, $C$ is an even cycle whose complement in $G$ contains a
 quasi-perfect matching. Since the orientation of $G^{u,v}$ is $3$-Pfaffian,
 Lemma~\ref{lem:cycle_odd} implies that $C$ is oddly oriented, hence $G-\{i\}$
 is Pfaffian. 

\medskip

For the converse, let $B$ be a minimal graph with respect to edge deletion such that the
2-suspension $B^{u,v}$ is non-3-Pfaffian. In particular, any triple
belongs to some spanning tree of $B^{u,v}$, otherwise the corresponding edge
in $B$  could have been deleted.

Choose $ab\in E(B)$ such that there is some spanning tree of $B^{u,v}$ 
containing neither $abu$ nor $abv$. (If every edge $ab$ of $B$ has the property
that each spanning tree of $B^{u,v}$ contains $abu$ or
$abv$ then $ab$ is in every quasi-perfect matching of $B$. It is
not difficult to see that this can only happen if $B$ is a set of
vertex disjoint edges and one path of length $2$. However in this case $B^{u,v}$ is
3-Pfaffian.) 

Let $G=B\backslash ab$. By minimality of $B$, the 3-graph $G^{u,v}$ is 3-Pfaffian.
Then there is a $u$-orientation and a $v$-orientation of the edges of
 $G$ with the property that all the spanning trees of $G^{u,v}$
have the same sign when triples $iju$ are oriented according to
the $u$-orientation of $ij$ and triples $ijv$ according to the
$v$-orientation of $ij$.

Extend both the $u$- and $v$-orientation of $G$ to  orientations of
$B$ by orienting the edge $ab$ in any way. Since the resulting orientation of
$B^{u,v}$
 is not 3-Pfaffian, there exist two quasi-perfect matchings $Q^+$ and $Q^-$ 
such that they both contain $ab$ and the associated spanning trees $T^+$ and
$T^-$ have opposite signs.

Let $Q=M_u\cup M_v$ be an arbitrary quasi-perfect matching of
$G$ and consider the graphs $H^+=Q\cup Q^+$ and $H^-=Q\cup Q^-$. Lemma~\ref{lem:unions_qpm} gives the possible subgraphs that can
arise as connected components of $Q\cup Q^+$ and $Q\cup Q^-$. If one of the connected
components of type (H) or (T)  is not a path of length 2 or 4, then
we can find one of the excluded subgraphs in condition (ii), just as in Claim
2 in the proof of Theorem~\ref{thm:carac_3pfaffian_or}. 
If this is the case we are
done, so suppose that all the connected components are even cycles or paths of
length 2 or 4. 

Our next goal is to show that the edge $ab$ belongs to one of these even cycles
and not to the paths. We look at $H^+$ since the argument is symmetric. If
$H^+$ contains a path of length 2, then the paths of length 2 in the quasi-perfect
matchings $Q$ and $Q^+$ coincide and, since $Q$ does not contain $ab$, it
follows that in this case $ab$ must belong to one of the cycles of $H^+$. If
$H^+$ contains a path of length four $x_1x_2x_3x_4x_5$, it means that one of
$Q$ or $Q^+$ contains the edges $\{x_1x_2,x_3x_4,x_4x_5\}$ and the
other contains the edges $\{x_1x_2,x_2x_3,x_4x_5\}$. Thus if the
edge $ab$ belongs to this path, it is either $x_2x_3$ or $x_3x_4$. We
assume it is $x_2x_3$, and hence that  $Q$ contains $x_1x_2$,
$x_3x_4$ and
$x_4x_5$. Therefore the component of type (H) in  $H^-$ is also the path
$x_1x_2x_3x_4x_5$. The other connected components in $H^+$ and $H^-$ are
cycles of even length that do not contain $ab$ and whose complement contains a
quasi-perfect matching in $G$ (having $x_3x_4x_5$ as its path of
length 2). Since the orientation is 3-Pfaffian in $G^{u,v}$, all  even cycles in
$H^+$ and $H^-$ are oddly oriented.  Then by Lemma~\ref{lem:sign_tree} it is easy to see
that  $T^+$ and $T^-$ either have both the same sign or both the opposite as the tree associated to $Q$,
which is not possible by the choice of $T^+$ and $T^-$. So we can conclude
that $H^+$ is a collection of cycles of even length and a
path of length 2 or 4, and that the edge $ab$ belongs to one of the cycles. 

Since $H^+$ is a spanning subgraph of $B$, it is only left to decide which
other edges we can have in addition to those of $H^+$. We show that if $B$
does not contain any of the subgraphs in (ii) then there is a vertex $i$ for
which $B-\{i\}$ is not a Pfaffian graph. Let us start by
analysing what happens if the path of $H^+$ has length 4. Let
$x_1y_1zy_2x_2$ be this path. Observe that vertex $x_1$ (and similarly
$x_2$) has degree one in $B$. Indeed, if $x_1$ was joined to a vertex other than
$y_1$ it would create a $P_6$, a $K_{2,3}^-$, a $C_3$ or a $C_5$, all of them
with a perfect matching in the complement. The vertex $z$ cannot be adjacent
to any of the even cycles of length at least $4$, since this would create either a $P_6$ or a
$K_{2,3}^-$ with a perfect matching in the complement. There can be edges
joining $z$ and some of the isolated edges (cycles of length 2) of $H^+$. Let $y_1,\ldots,y_i$
($i\geq 2$) be all the neighbours of $z$. There are edges $x_iy_i$, and
all the $x_i$ have degree one. Therefore the edges $x_iy_i$ belong to
every quasi-perfect matching of $B$ and, in order to cover $z$, each
quasi-perfect matching contains exactly one of
the edges $zy_i$. Therefore, if a
a vertex $y_i$ was in other edges than $zy_i$ and $x_iy_i$, then
these other edges would belong to no quasi-perfect matching. By minimality of
$B$ we conclude that $B$ has a connected component that is isomorphic to a
star with every edge subdivided; let us call this component $S$. The case where the (H) component of $H^+$ is
a path of length 2 is argued similarly and the same conclusion reached (i.e.,
that there is a component isomorphic to a star with every edge subdivided). 

It is easy to see that $S^{u,v}$ is a 3-Pfaffian graph. Indeed, take
$\uor{x_i}{y_i}$, $\uor{y_i}{z}$ and $\vor{x_i}{y_i}$, $\vor{z}{x_i}$. This
orientation satisfies the conditions described in Theorem~\ref{thm:carac_3pfaffian_or}.
If the rest of $B$, that is, $B-S$, had a Pfaffian orientation, we could use
it to extend  the orientation of $S$ just described to an orientation
satisfying the conditions of Theorem~\ref{thm:carac_3pfaffian_or}, and therefore
$B^{u,v}$ would be 3-Pfaffian. Hence, $B-S$, or $B-\{z\}$
in particular, is not a Pfaffian graph.
\hspace*{\fill} \end{proof}

By combining Theorem~\ref{thm:2-susp_3-Pfaffian} and Little's characterization of Pfaffian
bipartite graphs we obtain a characterization of $3$-Pfaffian $2$-suspensions
of bipartite graphs. 

\begin{cor}\label{cor:2-susp_obstructions}
Let $G$ be a bipartite graph and $u,v\not\in V(G)$. Then the 2-suspension
$G^{u,v}$ has a 3-Pfaffian orientation if and only if $G$ has none of
the following as subgraphs:

\begin{enumerate}[(i)]

\item  an even
subdivision of $K_{3,3}$ whose complement in $G$ has a quasi-perfect
matching;

\item a $P_6$ or $K_{2,3}^-$ whose complement in $G$ has a
perfect matching.

\end{enumerate}
\end{cor}


\section{Partial Steiner triple systems and 3-Pfaffian orientations}\label{sec:STS}


\subsection{Partial Steiner triple systems}\label{sec:intro_STS}

In this section we consider $3$-graphs $H$ with the property that the
multiplicity of every pair of vertices is at most $1$. Such a $3$-graph will be called a
 {\em partial Steiner triple system}. 

 Let $G$ be the underlying
graph of a partial Steiner triple system $H=(V,\Delta)$. For an edge $ij\in
E(G)$, the only $k\in V$ such that $ijk\in \Delta$ is denoted $n(ij)$. Recall
that Lemma~\ref{lem.sp trees pm + function} assigns to every spanning tree of
$H$ a pair $(M,f)$, where $M$ is a perfect matching of $G-v$ and
$f:M\rightarrow V$ is such that the triples of $T$ are $\{ijf(ij)\}$. If $H$
is a partial Steiner triple system, the function $f$ is necessarily
$n_{|M}$. In order to describe the perfect matchings that arise we need some
further definitions.

Let $t_1,t_2,\ldots,t_{\ell}$ be the triples of a cycle  spanning $2\ell
$ vertices, that is, there are $2\ell$ different vertices
$a_1,\ldots,a_{\ell}$ and $b_1,\ldots,b_{\ell}$ such that
$t_i=\{a_i,b_i,a_{i+1}\}$  ($a_{\ell+1}=a_1$). The $2\ell$-cycle of the
underlying graph $G$ with edges
$a_1b_1,b_1a_2,\ldots,a_{\ell}b_{\ell},b_{\ell}a_1$ will be called a
\emph{switching cycle}. We say that a perfect matching $M$ of $G$ \emph{alternates
  around a switching cycle} if there is another perfect matching $N$ such that
the symmetric difference $M\bigtriangleup N$ is a switching cycle.

\begin{cor}\label{cor.psts}
Suppose $H=(V,\Delta)$ is a partial Steiner triple system. For any fixed $v\in V$, spanning
trees of $H$ are in bijective correspondence with
perfect matchings of $G-v$ that do not alternate around a switching
cycle.
\end{cor}
\begin{proof}
It follows from Lemma~\ref{lem.sp trees pm + function} that for any
3-graph $H$ spanning trees are in one-one correspondence with pairs
$(M,f)$ where $M$ is a perfect matching of $G-v$ and $f:M\rightarrow
V$ is a function with the property that there are no cycles in
$\{ijf(ij):ij\in M\}$. As noted above, the function $f$ is uniquely 
determined from the matching, since each pair is in at most one triple. The
condition that there are no cycles in $\{ijf(ij):ij\in M\}$ translates directly
to the fact that $M$ does not alternate around a switching cycle.
\hspace*{\fill} \end{proof}

\begin{thm}\label{thm:PSTS_3-Pfaff} If $H=(V,\Delta)$ is a partial Steiner triple system with the
  property that $H-v$ has no cycles for some $v\in V$ then the number
  of spanning trees of $H$ is equal to the number of perfect matchings
  of $G-v$.
Furthermore, $H$ is 3-Pfaffian if and only if $G-v$ is Pfaffian.
\end{thm}
\begin{proof} Let $V=[2n\!+\!1]$.
If $G-v$ has no switching cycles, 
i.e.,\ if
$H-v$ has no cycles, then by Corollary~\ref{cor.psts} perfect matchings of $G-v$ are in
bijective correspondence with spanning trees of $H$. 

To prove the second part, we relate orientations of triples in $H-v$ 
to orientations of edges in $G-v$ so that we can express the
Masbaum-Vaintrob theorem in terms of edge orientations. 

A triangle $abc$ in $G-v$ is called \emph{black} if $abc$ is a triple of
$H-v$. An edge of $G-v$ is \emph{black} if it belongs to a black triangle; it
is \emph{white} otherwise (if $ab$ is white, then $abv$ is a triple of $H$). 
Given an orientation $\omega$ of $H$, we define an orientation of $G-v$ in the
following way. If $abc$ is a black triangle with $a<b<c$ and the corresponding
triple $abc$ is oriented $(a\, b\, c)$, orient
$a\longrightarrow b,b\longrightarrow c,c\longrightarrow a$. Otherwise, if
$abc$ is oriented $(a\, c\, b)$, orient
$a\longrightarrow c,c\longrightarrow b,b\longrightarrow a$. White edges are
arbitrarily oriented. 

The Hirschman-Reiner formulation of the Masbaum-Vaintrob theorem
(Theorem~\ref{thm: HR}) gives

$$\mathcal{P}^{\omega}(H,y)= \mathop{\sum_{\mbox{\rm \tiny perfect
      matchings $M$ of $G-v$}}} {\rm sgn}(M)\mathop{\prod_{ij\in M}}_{i<j}\epsilon_{i,j,n(ij)}\bar{y}_{ijn(ij)},$$  

where $n(ij)$ denotes the only vertex such that $ijn(ij)\in \Delta$ and
$\bar{y}_{ijn(ij)}$ equals $y_{ijn(ij)}$ or $-y_{ijn(ij)}$ according to
whether the orientation of $ijn(ij)$ equals or is opposite to the canonical orientation.

It is straightforward to check that, for $i<j$, 
$$\epsilon_{i,j,n(ij)}\bar{y}_{ijn(ij)}=\left\{\begin{array}{rl}  y_{ijn(ij)}
    & \mbox{ if } i\longrightarrow j;\\ -y_{ijn(ij)} & \mbox{ if }
    j\longrightarrow i. \end{array} \right.  $$
Therefore,
$$\mathcal{P}^{\omega}(H,y)= \mathop{\sum_{\mbox{\rm \tiny perfect
      matchings $M$ of $G-v$}}} {\rm sgn}(M)(-1)^{\#\{i<j: j\rightarrow i\}}\mathop{\prod_{ij\in M}}_{i<j} y_{ijn(ij)}.$$  

Thus if the orientation of $G-v$ is Pfaffian, the orientation $\omega$ is a
$3$-Pfaffian orientation of $H$, and conversely.
Therefore if $G-v$ has a Pfaffian orientation with the property that
each black triangle $ijk$ of $G-v$ is cyclically oriented 
then $H$ is
3-Pfaffian.

We show that any Pfaffian orientation of $G-v$ can be
converted into a Pfaffian orientation cyclic on black triangles of $G-v$. 

Let $abc$ be a black triangle of $G-v$ with some orientation of its
edges.
Suppose this orientation of
$abc$ is not already cyclic. Two of the edges of $abc$ must be in the same
direction when traversing the triangle, say $ab$ and $bc$. Then $abc$ can be
cyclically oriented by reversing the directions of all edges incident with $b$
or by reversing the direction of all edges incident with $a$ and then of those
incident with $c$.   Reversing the direction of all the edges
incident with a given vertex of $G-v$ preserves the property of being
a Pfaffian orientation,
since any even cycle has its parity of forward edges preserved.

We next show how to combine these movements to make all black triangles
cyclic. Since $H-v$ is a forest, there is some ordering
$\tau_1,\ldots,\tau_{\ell}$ of the black triangles such that $|(\cup_{j\leq i}
\tau_j )\cap \tau_{i+1}|\leq 1$. Inductively, suppose that the first $i$ black
triangles are cyclically oriented. Let $a,b$ be two vertices of $\tau_{i+1}$
that do not belong to $\cup_{j\leq i} \tau_j$. Then if $\tau_{i+1}$ is not
cyclically oriented it can be made so by
reversing the orientation of all edges
incident with $a$, or with $b$, or with both. This clearly leaves all
black triangles already processed unaltered, so eventually all black triangles
are cyclically oriented, as needed. \hspace*{\fill} \end{proof}

In particular, in Theorem~\ref{thm:PSTS_3-Pfaff} if $H$ is such that
$H-v$ has no cycles and $G-v$ is planar
then $G-v$ is Pfaffian. In this case the Pfaffian tree polynomial
$P^{\omega}(H;y)$ is up to sign equal to the tree generating polynomial
$P(H;y)$ when $\omega$ is a 3-Pfaffian orientation of $H$. 
Galluccio and Loebl~\cite{GL99} prove a statement first made by
Kasteleyn that the generating function for perfect matchings of a graph embeddable in an orientable
surface of genus $g$ may be written as a linear combination of $4^g$
Pfaffians (with coefficients independent of the graph). Suppose we have a 3-graph
$H=(V,\Delta)$ with the property that there is $v\in V$ such that the graph
$G-v$ underlying $H-v$ is without triple cycles and is of genus
$g$. Then we can use the one-one sign-preserving correspondence
between spanning trees of $H$ and perfect matchings of $G-v$ to deduce
a similar result: there are $4^g$
triple orientations of $H$ such that the tree generating polynomial
$\mathcal{P}(H,y)$ can be expressed as a linear
combination of $4^g$ signed tree generating polynomials
$\mathcal{P}^{\omega}(H,y)$, where $\omega$ ranges over $4^g$ triple orientations.

\subsection{Minimal non-3-Pfaffian 3-graphs}\label{ssec:minimal}

By Theorem~\ref{thm: existence coNP} there is a polynomial-size
certificate witnessing a non-3-Pfaffian 3-graph. Even if the number of spanning trees is exponential in $n$,
there is a polynomial-size subset of spanning trees of
  $H$  whose elements cannot be made all the same
sign.
In view of the fact that a non-minimal non-3-Pfaffian
$3$-graph can be reduced to a minimal non-3-Pfaffian sub-$3$-graph by
deletion and contraction of triples, it is natural to ask whether there is a
finite set of obstructions to being 3-Pfaffian, such as given by
Corollary~\ref{cor:2-susp_obstructions} for 2-suspensions of
graphs. In this subsection we show that this is not the case by giving an
infinite collection of minimal non-3-Pfaffian graphs (see
Theorem~\ref{thm:interlaced}).

\begin{figure}[htb]
\begin{center}
\caption{Some non-3-Pfaffian $3$-graphs $H$ minimal with respect to deletion
  and contraction of triples, given by their underlying graph with
  one vertex deleted. Edges not in shaded triangles are pairs of
  vertices in a triple containicng the removed
  vertex.}\label{fig:minimal_non-3-Pfaffian}
\vspace{0.3cm}
\scalebox{1} 
{\begin{pspicture}(0,-2.25375)(12.485937,2.25375)
\definecolor{color2055b}{rgb}{0.8,0.8,0.8}
\psbezier[linewidth=0.02](4.6746874,-1.488125)(5.0946875,-2.168125)(6.3946877,-1.948125)(6.9146876,-1.508125)(7.4346876,-1.068125)(7.6946874,-0.308125)(7.6346874,0.111875)
\psbezier[linewidth=0.02](1.7346874,1.231875)(3.0946875,1.291875)(3.4546876,0.091875)(3.0746875,-0.688125)
\psbezier[linewidth=0.02](0.3146875,-0.748125)(0.3146875,-1.548125)(3.0746875,-1.508125)(3.0746875,-0.708125)
\psbezier[linewidth=0.02](1.6346875,1.251875)(0.7346875,1.211875)(0.0346875,0.091875)(0.3346875,-0.748125)
\pstriangle[linewidth=0.02,dimen=outer,fillstyle=solid,fillcolor=color2055b](1.6946875,0.251875)(1.4,0.98)
\pstriangle[linewidth=0.02,dimen=outer,fillstyle=solid,fillcolor=color2055b](1.0146875,-0.708125)(1.4,0.98)
\pstriangle[linewidth=0.02,dimen=outer,fillstyle=solid,fillcolor=color2055b](2.3746874,-0.708125)(1.4,0.98)
\psdots[dotsize=0.2,fillstyle=solid,fillcolor=black,dotstyle=o](1.6746875,-0.708125)
\psdots[dotsize=0.2,fillstyle=solid,fillcolor=black,dotstyle=o](0.9946875,0.271875)
\psdots[dotsize=0.2,fillstyle=solid,fillcolor=black,dotstyle=o](2.3746874,0.271875)
\psdots[dotsize=0.2,fillstyle=solid,dotstyle=o](0.3146875,-0.708125)
\psdots[dotsize=0.2,fillstyle=solid,dotstyle=o](3.0546875,-0.728125)
\psdots[dotsize=0.2,fillstyle=solid,dotstyle=o](1.6946875,1.231875)
\psbezier[linewidth=0.02](4.1846876,0.051875)(4.0946875,-0.09947635)(4.4146876,-0.948125)(4.5346875,-1.108125)
\psbezier[linewidth=0.02](5.8746877,1.831875)(5.1546874,1.751875)(3.688464,0.9717801)(3.6746874,0.011875)(3.660911,-0.9480301)(5.3946877,-1.808125)(5.8946877,-1.648125)
\psbezier[linewidth=0.02](10.374687,-1.6494293)(11.674687,-1.508125)(12.134687,-0.608125)(12.134687,0.071875)
\psbezier[linewidth=0.02](10.354688,1.8459886)(9.47878,1.991875)(8.474688,0.671875)(8.654688,0.131875)
\pstriangle[linewidth=0.02,dimen=outer,fillstyle=solid,fillcolor=color2055b](10.404688,0.791875)(1.42,1.06)
\rput{271.53976}(11.190495,11.73909){\pstriangle[linewidth=0.02,dimen=outer,fillstyle=solid,fillcolor=color2055b](11.624687,-0.408125)(1.42,1.06)}
\rput{-179.33186}(20.821007,-1.9948504){\pstriangle[linewidth=0.02,dimen=outer,fillstyle=solid,fillcolor=color2055b](10.404688,-1.588125)(1.42,1.06)}
\rput{90.46375}(9.40106,-9.081525){\pstriangle[linewidth=0.02,dimen=outer,fillstyle=solid,fillcolor=color2055b](9.204687,-0.408125)(1.42,1.06)}
\psdots[dotsize=0.2](9.714687,0.791875)
\psdots[dotsize=0.2](11.074688,0.831875)
\psdots[dotsize=0.2](11.114688,-0.568125)
\psdots[dotsize=0.2](9.734688,-0.568125)
\psdots[dotsize=0.2,fillstyle=solid,dotstyle=o](10.394688,-1.628125)
\psdots[dotsize=0.2,fillstyle=solid,dotstyle=o](12.154688,0.111875)
\psdots[dotsize=0.2,fillstyle=solid,dotstyle=o](10.394688,1.831875)
\psdots[dotsize=0.2,fillstyle=solid,dotstyle=o](8.674687,0.111875)
\pstriangle[linewidth=0.02,dimen=outer,fillstyle=solid,fillcolor=color2055b](5.9046874,0.771875)(1.42,1.06)
\rput{-90.27026}(7.0564203,7.226964){\pstriangle[linewidth=0.02,dimen=outer,fillstyle=solid,fillcolor=color2055b](7.1246877,-0.428125)(1.42,1.06)}
\rput{-179.71352}(11.814792,-2.1667132){\pstriangle[linewidth=0.02,dimen=outer,fillstyle=solid,fillcolor=color2055b](5.9046874,-1.628125)(1.42,1.06)}
\rput{-270.57037}(4.7597246,-4.6035933){\pstriangle[linewidth=0.02,dimen=outer,fillstyle=solid,fillcolor=color2055b](4.7046876,-0.428125)(1.42,1.06)}
\psdots[dotsize=0.2](5.2146873,0.771875)
\psdots[dotsize=0.2](6.5746875,0.811875)
\psdots[dotsize=0.2](6.6146874,-0.588125)
\psdots[dotsize=0.2](5.2346873,-0.588125)
\psdots[dotsize=0.2,fillstyle=solid,dotstyle=o](5.8946877,-1.648125)
\psdots[dotsize=0.2,fillstyle=solid,dotstyle=o](7.6546874,0.091875)
\psdots[dotsize=0.2,fillstyle=solid,dotstyle=o](5.8946877,1.811875)
\psdots[dotsize=0.2,fillstyle=solid,dotstyle=o](4.1746874,0.091875)
\usefont{T1}{ppl}{m}{n}
\rput(1.6714063,1.506875){\footnotesize 1}
\usefont{T1}{ppl}{m}{n}
\rput(0.06671875,-0.873125){\footnotesize 2}
\usefont{T1}{ppl}{m}{n}
\rput(3.305625,-0.873125){\footnotesize 3}
\usefont{T1}{ppl}{m}{it}
\rput(2.0979688,0.146875){\footnotesize a}
\usefont{T1}{ppl}{m}{it}
\rput(1.2771875,0.126875){\footnotesize b}
\usefont{T1}{ppl}{m}{it}
\rput(1.6721874,-0.433125){\footnotesize c}
\usefont{T1}{ppl}{m}{n}
\rput(5.851406,2.066875){\footnotesize 1}
\usefont{T1}{ppl}{m}{n}
\rput(3.9467187,0.126875){\footnotesize 2}
\usefont{T1}{ppl}{m}{n}
\rput(5.885625,-2.093125){\footnotesize 3}
\usefont{T1}{ppl}{m}{n}
\rput(7.869375,0.126875){\footnotesize 4}
\usefont{T1}{ppl}{m}{it}
\rput(6.377969,0.626875){\footnotesize a}
\usefont{T1}{ppl}{m}{it}
\rput(5.3971877,0.626875){\footnotesize b}
\usefont{T1}{ppl}{m}{it}
\rput(5.3921876,-0.373125){\footnotesize c}
\usefont{T1}{ppl}{m}{it}
\rput(6.4082813,-0.373125){\footnotesize d}
\usefont{T1}{ppl}{m}{n}
\rput(10.351406,2.106875){\footnotesize 1}
\usefont{T1}{ppl}{m}{it}
\rput(9.877188,0.626875){\footnotesize b}
\usefont{T1}{ppl}{m}{it}
\rput(10.897968,0.626875){\footnotesize a}
\usefont{T1}{ppl}{m}{n}
\rput(8.406719,0.126875){\footnotesize 2}
\usefont{T1}{ppl}{m}{n}
\rput(12.369375,0.126875){\footnotesize 4}
\usefont{T1}{ppl}{m}{it}
\rput(9.872188,-0.333125){\footnotesize c}
\usefont{T1}{ppl}{m}{it}
\rput(10.928281,-0.333125){\footnotesize d}
\usefont{T1}{ppl}{m}{n}
\rput(10.405625,-1.873125){\footnotesize 3}
\end{pspicture} 
}

\end{center}
\end{figure}
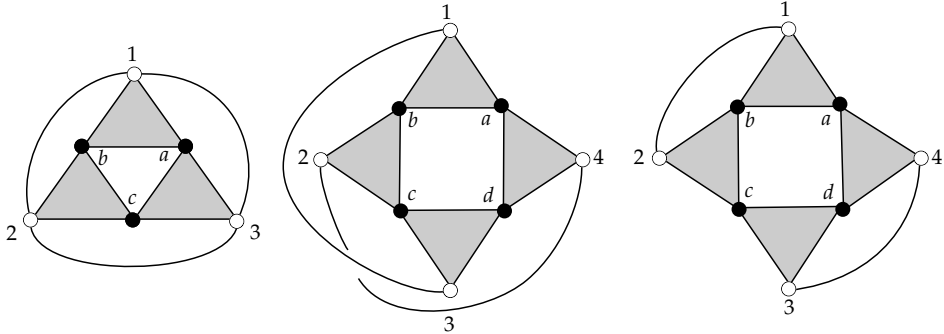

The 3-graphs $H$ in Table~\ref{table:minimal_non-3-Pfaffian} are minimally
non-3-Pfaffian and $H-\{0\}$ has underyling graph $G-\{0\}$ of the
form illustrated in Figure~\ref{fig:minimal_non-3-Pfaffian}.
The orientation of a spanning tree is given as the cyclic permutation
of the vertex set obtained as product of
  3-cycles; to form this product the oriented triples of the spanning
  tree are taken in the order given in the previous column of the
  table. The sign of the orientation is relative
  to the order of vertices given in the first column. In each case
  there are an odd number of negative spanning trees. It is readily
  checked that a given triple belongs to an even number of spanning
  trees, and therefore that it is not possible to change triple
  orientations to obtain spanning trees all of the same sign. 

\begin{table}[htb]
\caption{Three non-3-Pfaffian 3-graphs minimal with respect to
  deletion and contraction of triples.}\label{table:minimal_non-3-Pfaffian}
\begin{center}
\begin{tabular}{c|c|c|c|c}
Vertices & Oriented triples & Spanning tree &
Orientation & Sign\\ [1ex]
\hline
$0,1,2,3,$ & $012, 023, 031,$ & $\{012, 1ab, 3ca\}$ & $(\, 0\;
3\; c\; a\; b\; 1\; 2\,)$ & $+$\\
$a,b,c$ & $1ab, 2bc, 3ca$ & $\{012, 2bc, 3ca\}$ & $(\, 0\; 1\; b\; a\; 3\;
c\; 2\,)$ & $-$\\
 & & $\{023, 2bc, 1ab\}$ & $(\, 0\; 1\; a\; b\; c\; 2\; 3\,)$ & $+$\\
 & & $\{023, 3ca, 1ab\}$ & $(\, 0\; 2\; c\; b\; 1\; a\; 3\,)$ & $-$ \\
 & & $\{031, 3ca, 2bc\}$ & $(\, 0\; 2\; b\; c\; a\; 3\; 1\,)$ & $+$ \\
 & & $\{031, 1ab, 2bc\}$ & $(\, 0\; 3\; a\; c\; 2\; b\; 1\,)$ &
$-$\\[1ex]
\hline 
$0,1,2,3,4,$ & $013, 024,$ & $\{013, 1ab, 2bc, 4da\}$ & $(\, 0\; 4\;
d\; a\; c\; 2\; b\; 1\; 3\,)$ & $-$ \\
$a,b,c,d$ &  $1ab, 2bc, 3cd, 4da$ & $\{013, 3cd, 2bc, 4da\}$ & $(\,
0\; 1\; 2\; b\; c\; a\; 4\; d\; 3\,)$ & $+$\\
 & & $\{024, 2bc, 1ab, 3cd\}$ & $(\, 0\; 1\; a\; b\; d\; 3\; c\; 2\;
4\,)$ & $-$\\
 & & $\{024, 4da, 1ab, 3cd\}$ & $(\, 0\; 2\; 3\; c\; d\; b\; 1\; a\;
4\,)$ & $+$\\
& & $\{013, 024, 2bc, 4da\}$ & $(\, 0\; 1\; 3\; b\; c\; 2\; d\; a\;
4\,)$ & $+$\\
& & $\{013, 024, 1ab, 3cd\}$ & $(\, 0\; a\; b\; 1\; c\; d\; 3\; 2\;
4\,)$ & $-$\\[1ex]

\hline
$0,1,2,3,4,$ & $012, 034,$ & $\{012, 2bc, 3cd, 4da\}$ & $(\, 0\; 1\;
b\; a\; 4\; d\; 3\; c\; 2\,)$ & $+$\\
$a,b,c,d$ & $1ab, 2bc, 3cd, 4da$  & $\{012, 1ab, 4da, 3cd\}$ &  $(\,
0\; 4\; 3\; c\; d\; a\; b\; 1\; 2\,)$ & $-$\\
 & & $\{034, 4da, 1ab, 2bc\}$ &  $(\, 0\; 3\; d\; c\; 2\; b\; 1\; a\;
4\,)$ & $+$ \\
 & & $\{034, 3cd, 2bc, 1ab\}$ &  $(\, 0\; 2\; 1\; a\; b\; c\; d\; 3\;
4\,)$ & $-$\\
& & $\{012, 1ab, 034, 3cd\}$ & $(\, 0\; a\; b\; 1\; 2\; c\; d\; 3\;
4\,)$ & $+$ \\
& & $\{012, 2bc, 034, 4da\}$ & $(\, 0\; 1\; b\; c\; 2\; 3\; d\; a\;
4\,)$ & $-$\\[1ex]
\hline

\end{tabular}
\end{center}
\end{table}

\begin{prop}\label{lem:two_ind_white_edges} Let $H$ be a 3-graph on vertices $0,1,2,\ldots, 2k$ with
  triples 
$$\{2k,1,2\}, \{2,3,4\}, \{4,5,6\},\ldots,\{2k-2,2k-1,2k\} $$

and containing two triples of the form
$$\{0,2x\!-\!1,2y\!-\!1\}, \{0,2z\!-\!1, 2t\!-\!1\}$$

for some distinct $x,y,z,t$. 
Then $H$ is non-3-Pfaffian.
Similarly, a 3-graph with triples 
$$\{2k,1,2\}, \{2,3,4\}, \{4,5,6\},\ldots,\{2k-2,2k-1,2k\} $$
and three triples of the form 
 $$\{0,2x\!-\!1,2y\!-\!1\}, \{0,2y\!-\!1, 2z\!-\!1\}, \{0,2z\!-\!1,
2x\!-\!1\}$$ for some distinct $x,y,z,$ is non-3-Pfaffian.
\end{prop}

\begin{proof}
Since the property of being 3-Pfaffian is preserved by deletion and
contraction of triples we may assume in the first case that $k=4$ and
$\{x,y,z,t\}=\{1,2,3,4\}$ and in the second case that $k=3$ and $\{x,y,z\}=\{1,2,3\}$.
These cases are the non-3-Pfaffian 3-graphs given in Table~\ref{table:minimal_non-3-Pfaffian}.
\hspace*{\fill} \end{proof}

\begin{figure}[htb]
\caption{A 3-Pfaffian and a non-3-Pfaffian 3-graph, given by the
  underlying graph with a vertex deleted. These are first in the
  family of 3-graphs of Theorem~\ref{thm:interlaced} that are 3-Pfaffian or
  non-3-Pfaffian according to the parity of the number of shaded triangles.}\label{fig:interlaced}
\vspace{-0.2cm}
\begin{center}

\scalebox{1} 
{
\begin{pspicture}(0,-2.2192183)(10.089603,2.2292182)
\definecolor{color5628b}{rgb}{0.8,0.8,0.8}
\psbezier[linewidth=0.02](9.432729,0.1092165)(9.47,0.35796347)(9.414918,0.4734006)(9.393051,0.617287)
\psbezier[linewidth=0.02](8.727368,-1.3748988)(9.313247,-1.1680841)(9.961102,-0.7091343)(10.020352,-0.060870465)(10.079603,0.58739334)(8.969313,0.8813464)(8.4902935,0.721122)
\psbezier[linewidth=0.02](7.7900014,-1.6109397)(8.070582,-1.603885)(8.195729,-1.5367408)(8.354991,-1.4943377)
\psbezier[linewidth=0.02](6.3463435,-0.8342227)(6.5419655,-1.4015167)(6.9961953,-2.0390015)(7.6641884,-2.1241097)(8.332182,-2.2092183)(8.670155,-1.1644716)(8.519377,-0.70103914)
\psbezier[linewidth=0.02](6.131365,0.069216594)(6.13,-0.19839682)(6.1955175,-0.32058907)(6.2344184,-0.47418845)
\psbezier[linewidth=0.02](6.72,1.3392181)(6.06,1.0691296)(5.34,0.54682875)(5.3,-0.08980841)(5.26,-0.72644556)(6.54,-0.8807819)(7.08,-0.66856945)
\psbezier[linewidth=0.02](7.78,1.6792182)(7.46,1.6392181)(7.32,1.5592182)(7.14,1.4992181)
\psbezier[linewidth=0.02](0.98,0.51921815)(1.16,0.6792181)(1.58,0.9792181)(1.84,1.0792181)
\psbezier[linewidth=0.02](3.14,-0.8807819)(3.36,-0.6317448)(3.22,-0.24300408)(3.16,-0.06078186)
\psbezier[linewidth=0.02](0.5184615,-0.90078187)(0.58,-1.1207819)(1.02,-1.3207818)(1.04,-1.3407818)
\psbezier[linewidth=0.02](0.7917391,0.35921815)(0.34,-0.04078186)(0.0,-0.60078186)(0.118695654,-1.0282818)(0.23739131,-1.4557818)(1.52,-1.6407819)(1.82,-0.9291747)
\pstriangle[linewidth=0.02,dimen=outer,fillstyle=solid,fillcolor=color5628b](1.86,0.09921814)(1.4,0.98)
\pstriangle[linewidth=0.02,dimen=outer,fillstyle=solid,fillcolor=color5628b](1.18,-0.86078185)(1.4,0.98)
\pstriangle[linewidth=0.02,dimen=outer,fillstyle=solid,fillcolor=color5628b](2.54,-0.86078185)(1.4,0.98)
\psdots[dotsize=0.2,fillstyle=solid,fillcolor=black,dotstyle=o](1.84,-0.86078185)
\psdots[dotsize=0.2,fillstyle=solid,fillcolor=black,dotstyle=o](1.16,0.11921814)
\psdots[dotsize=0.2,fillstyle=solid,fillcolor=black,dotstyle=o](2.54,0.11921814)
\psdots[dotsize=0.2,fillstyle=solid,dotstyle=o](0.48,-0.86078185)
\psdots[dotsize=0.2,fillstyle=solid,dotstyle=o](3.22,-0.8807819)
\psdots[dotsize=0.2,fillstyle=solid,dotstyle=o](1.86,1.0792181)
\psbezier[linewidth=0.02](1.22,-1.4257343)(1.8,-1.5607818)(2.64,-1.6807818)(3.48,-1.1007819)(4.32,-0.5207819)(3.44,0.13921814)(2.5481012,0.13749222)
\psbezier[linewidth=0.02](3.06,0.15921813)(2.74,0.83921814)(2.38,1.3592181)(1.52,1.5592182)(0.66,1.7592181)(0.54,0.89921814)(1.16,0.11921814)
\pstriangle[linewidth=0.02,dimen=outer,fillstyle=solid,fillcolor=color5628b](7.78,0.71921813)(1.4,0.98)
\rput{-89.76206}(8.873573,9.008937){\pstriangle[linewidth=0.02,dimen=outer,fillstyle=solid,fillcolor=color5628b](8.96,-0.44078186)(1.4,0.98)}
\rput{-269.1915}(6.742629,-6.58971){\pstriangle[linewidth=0.02,dimen=outer,fillstyle=solid,fillcolor=color5628b](6.62,-0.46078187)(1.4,0.98)}
\rput{179.33484}(15.586579,-2.3120394){\pstriangle[linewidth=0.02,dimen=outer,fillstyle=solid,fillcolor=color5628b](7.8,-1.6007819)(1.4,0.98)}
\psdots[dotsize=0.2,fillstyle=solid,dotstyle=o](7.78,1.6992182)
\psdots[dotsize=0.2,fillstyle=solid,dotstyle=o](6.14,0.01921814)
\psdots[dotsize=0.2,fillstyle=solid,dotstyle=o](7.8,-1.5807818)
\psdots[dotsize=0.2,fillstyle=solid,dotstyle=o](9.44,0.05921814)
\psdots[dotsize=0.2](7.1,-0.6407819)
\psdots[dotsize=0.2](8.48,-0.62078184)
\psdots[dotsize=0.2](8.48,0.71921813)
\psdots[dotsize=0.2](7.08,0.71921813)
\psbezier[linewidth=0.02](9.34,0.83921814)(9.17542,1.3632543)(8.946499,2.0295882)(7.9964643,2.1244032)(7.046429,2.2192183)(6.42,1.1392181)(7.0907664,0.7392181)
\usefont{T1}{ppl}{m}{n}
\rput(2.0567188,1.1342181){\footnotesize 1}
\usefont{T1}{ppl}{m}{n}
\rput(0.37203124,-0.60578185){\footnotesize 2}
\usefont{T1}{ppl}{m}{n}
\rput(3.0109375,-1.0657818){\footnotesize 3}
\usefont{T1}{ppl}{m}{it}
\rput(1.4175,0.014218139){\footnotesize c}
\usefont{T1}{ppl}{m}{it}
\rput(1.8632812,-0.54578185){\footnotesize a}
\usefont{T1}{ppl}{m}{it}
\rput(2.2225,-0.02578186){\footnotesize b}
\usefont{T1}{ppl}{m}{n}
\rput(8.016719,1.7742181){\footnotesize 1}
\usefont{T1}{ppl}{m}{n}
\rput(6.052031,0.29421815){\footnotesize 2}
\usefont{T1}{ppl}{m}{n}
\rput(7.5109377,-1.6057819){\footnotesize 3}
\usefont{T1}{ppl}{m}{n}
\rput(9.394688,-0.22578186){\footnotesize 4}
\usefont{T1}{ppl}{m}{it}
\rput(8.283281,0.57421815){\footnotesize a}
\usefont{T1}{ppl}{m}{it}
\rput(7.2425,0.5542181){\footnotesize b}
\usefont{T1}{ppl}{m}{it}
\rput(7.2575,-0.38578185){\footnotesize c}
\usefont{T1}{ppl}{m}{it}
\rput(8.273594,-0.38578185){\footnotesize d}
\end{pspicture} 
}

\end{center}
\end{figure}
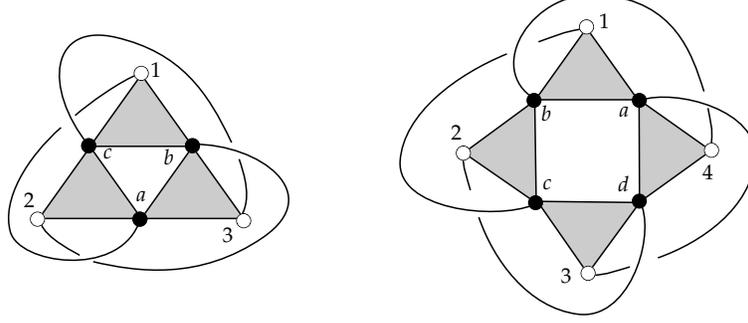

The 3-graphs in Table~\ref{table:interlaced} are illustrated in
Figure~\ref{fig:interlaced}. 

\begin{table}[htb]
\caption{A 3-Pfaffian and a non-3-Pfaffian 3-graph.}\label{table:interlaced}
\begin{center}
\begin{tabular}{c|c|c|c|c}
Vertices & Oriented triples & Spanning tree &
Orientation & Sign\\ [1ex]
\hline
$0,1,2,3,$ & $01a, 02b, 03c,$  & $\{01a, 2ca,
3ab\}$ & $(\,0\;1\;2\;c\;b\;3\;a\,)$ & $-$\\
$a,b,c$ & $1bc, 2ca, 3ab$ & $\{02b, 3ab, 1bc\}$ & $(\, 0\; 2\; 3\; a\; c\; 1\;
b\,)$ & $-$\\
& & $\{03c, 1bc, 2ca\}$ & $(\, 0\; 3\; 1\; b\; a\; 2\; c\,)$ & $-$\\
& & $\{01a, 02b, 03c\}$ &  $(\,0\;1\;a\;2\;b\;3\;c\,)$ & $-$\\[1ex]
\hline
 $0,1,2,3,4,$ & $01c, 02d, 03a, 04b,$ & $\{01c, 2bc, 3cd, 4da\}$ & $(\,0\; 1\; 2\; b\; a\; 4\; d\; 3\;
c\,)$ & $+$\\
$a,b,c,d$ & $1ab, 2bc, 3cd, 4da$ & $\{02d, 3cd, 4da, 1ab\}$ & $(\, 0\; 2\; 3\; c\; b\; 1\; a\; 4\;
d\,)$ & $+$ \\
 & & $\{03a, 4da, 1ab, 2bc\}$ &$(\, 0\; 3\; 4\; d\; c\; 2\; b\; 1\;
a\,)$ & $+$\\[0.25ex]
 & & $\{04b, 1ab, 2bc, 3cd\}$ & $(\, 0\; 4\; 1\; a\; d\; 3\; c\; 2\;
b\,)$ & $+$ \\
&  & $\{01c, 03a, 2bc, 4da\}$ & $(\, 0\;
1\; 2\; b\; c\; 3\; 4\; d\; a\,)$ & $-$\\
 & & $\{02d, 04b, 3cd, 1ab\}$ & $(\, 0\; 2\; 3\; c\; d\; 4\; 1\; a\;
b\,)$ & $-$\\[0.25ex]
&  & $\{01c, 02d, 03a, 04b\}$ & $(\,
0\; 1\; c\; 2\; d\; 3\; a\; 4\; b\,)$ & $+$\\[1ex]
\hline
\end{tabular}
\end{center}
\end{table}

In this case we cannot determine directly from the table whether the second
graph is $3$-Pfaffian or not. To do so, we transform the problem into an
algebraic one. 
The incidence matrix for triples (rows) and spanning trees
(columns)
is as follows (with spanning trees in the same order as in the table and each
column labelled by the sign of the corresponding tree):

$$\begin{array}{c|cccc|cc|c}
 & + & + & + & + & - & - & +\\
\hline
01c & 1 & 0 & 0 & 0  & 1 & 0 & 1\\
02d & 0 & 1 & 0 & 0  & 0 & 1 & 1\\
03a & 0 & 0 & 1 & 0  & 1 & 0 & 1\\
04b & 0 & 0 & 0 & 1  & 0 & 1 & 1\\
\hline 1ab & 0 & 1 & 1 & 1  & 0 & 1 & 0\\
2bc & 1 & 0 & 1 & 1 & 1 & 0 & 0\\
3cd & 1 & 1 & 0 & 1 & 0 & 1 & 0\\
4da & 1 & 1 & 1 & 0 & 1 & 0 & 0
\end{array}
$$

The non-zero positions in the row indexed by triple $t$ correspond to those trees that will change
sign if triple $t$ changes orientation. Therefore, finding a $3$-Pfaffian
orientation is equivalent to finding a subset of rows whose sum (modulo $2$)
is either $\left(\begin{array}{ccccccc}1 &  1 & 1 & 1 & 0 & 0 &
  1\end{array}\right)$ or $\left(\begin{array}{ccccccc} 0 & 0 & 0 & 0
    & 1 & 1 & 0\end{array}\right)$ (in the first case all trees
would be negative and in the second case they would be positive).  In other
words, we need to check whether either of the two vectors belongs to the row
span of the matrix over $\mathbb{F}_2$.

Since row $i$ and row $i+4$ for each $i=1,2,3,4$ sum to the all-one vector
$\left(\begin{array}{ccccccc} 1 & 1 & 1 & 1 & 1 & 1 & 1\end{array}\right)$, the row span is the rank $5$ subspace of $\mathbb{F}_2^7$, generated
by the rows of the matrix
$$\left(\begin{array}{ccccccc}
1 & 0 & 0 & 0 & 1 & 1 & 0\\
 0 & 1 & 0 & 0 & 1 & 0 & 1\\
0 & 0 & 1 & 0 & 1 & 1 & 0\\
 0 & 0 & 0 & 1 & 1 & 0 & 1\\
1 & 1 & 1 & 1 & 1 & 1 & 1\end{array}\right).$$

A 3-Pfaffian orientation exists if and only if the vector $\left(\begin{array}{ccccccc}
0 & 0 & 0 & 0 & 1 & 1 & 0\end{array}\right)$  is spanned by the rows of
this matrix. This is easily seen not to be the case. Hence no
orientation of triples can make all spanning trees have the same sign,
i.e., the second graph in Table~\ref{table:interlaced}  is a non-3-Pfaffian 3-graph.

The two graphs in Table~\ref{table:interlaced} are the first members of an
infinite family. The next member is given in
Table~\ref{table:5-cycle}; it is the $3$-graph $H$ for which the underlying graph $G-\{0\}$ of
$H-\{0\}$ consists of a 5-cycle of triangles $1ab,
2bc, 3cd, 4de, 5ea$ with edges $1c, 2d, 3e, 4a, 5b$.

\begin{table}[htb]
\caption{A 3-Pfaffian 3-graph.}\label{table:5-cycle}
\begin{center} 
\begin{tabular}{c|c|c|c}
 Oriented triples & Spanning tree &
Orientation & Sign\\ [1ex]
\hline
 
$01c, 02d, 03e, 04a, 05b,$ &
$\{01c,2bc,3cd,4de,5ea\}$ & $(\, 0\; 1\; 2\; b\; a\; 5\; e\; 4\; d\;
3\; c\,)$ & $+$\\
 $1ab, 2bc, 3cd, 4de, 5ea$ & $\{02d,3cd,4de,5ea,1ab\}$ & $(\, 0\; 2\;
3\; c\; b\; 1\; a\; 5\; e\; 4\; d\,)$ & $+$\\
& $\{03e,4de,5ea,1ab, 2bc\}$ & $(\, 0\; 3\; 4\; d\; c\; 2\; b\; 1\;
a\; 5\; e\,)$ & $+$\\
& $\{04a, 5ea,1ab, 2bc, 3cd\}$ & $(\, 0\; 4\; 5\; e\; d\; 3\; c\; 2\;
b\; 1\; a\,)$ & $+$\\
& $\{05b, 1ab,2bc, 3cd, 4de\}$ & $(\, 0\; 5\; 1\; a\; e\; 4\; d\; 3\;
c\; 2\; b\,)$ & $+$ \\[0.5ex]
& $\{01c, 04a, 2bc, 3cd, 5ea\}$ & $(\, 0\; 1\; 2\; b\; d\; 3\; c\; 4\;
5\; e\; a\,)$ & $-$\\
& $\{02d, 05b, 3cd, 4de, 1ab\}$ & $(\,0 \; 2\; 3\; c\; e\; 4\; d\; 5\;
1\; a\; b\,)$ & $-$\\
& $\{03e, 01c, 4de, 5ea, 2bc\}$ & $(\, 0\; 3\; 4\; d\; a\; 5\; e\; 1\;
2\; b\; c\,)$ & $-$\\
& $\{04a, 02d, 5ea, 1ab, 3cd\}$ & $(\, 0\; 4\; 5\; e\; b\; 1\; a\; 2\;
3\; c\; d\,)$ & $-$\\
& $\{05b, 03e, 1ab, 2bc, 4de\}$ & $(\, 0\; 5\; 1\; a\; c\; 2\; b\; 3\;
4\; d\; e\,)$ & $-$\\[0.5ex]
& $\{01c, 02d, 03e, 04a, 05b\}$ & $(\, 0\; 1\; 2\; b\; d\; 3\; c\; 4\;
5\; e\; a\,)$ & $+$ \\
\hline
\end{tabular}
\end{center}
\end{table}

The triple--spanning tree incidence matrix is here---taking spanning
trees in the order given in Table~\ref{table:5-cycle}---
given by
$$\begin{array}{c|ccccc|ccccc|c}
& + & + & + & + & + & - & - & - & - & - & +\\
\hline
01c & 1 & 0 & 0 & 0 & 0 & 1 & 0 & 1 & 0 & 0 & 1\\
02d & 0 & 1 & 0 & 0 & 0 & 0 & 1 & 0 & 1 & 0 & 1\\
03e & 0 & 0 & 1 & 0 & 0 & 0 & 0 & 1 & 0 & 1 & 1\\
04a & 0 & 0 & 0 & 1 & 0 & 1 & 0 & 0 & 1 & 0 & 1\\
05b & 0 & 0 & 0 & 0 & 1 & 0 & 1 & 0 & 0 & 1 & 1\\
\hline 
1ab & 0 & 1 & 1 & 1 & 1 & 0 & 1 & 0 & 1 & 1 & 0\\
2bc & 1 & 0 & 1 & 1 & 1 & 1 & 0 & 1 & 0 & 1 & 0\\
3cd & 1 & 1 & 0 & 1 & 1 & 1 & 1 & 0 & 1 & 0 & 0\\
4de & 1 & 1 & 1 & 0 & 1 & 0 & 1 & 1 & 0 & 1 & 0\\
5ea & 1 & 1 & 1 & 1 & 0 & 1 & 0 & 1 & 1 & 0 & 0
\end{array}
$$

Simple inspection shows that the sum of rows 2 to 6 is the vector\newline
$\left(\begin{array}{ccccccccccc} 0 & 0 & 0 & 0 & 0 & 1 & 1 & 1 & 1 &
  1 & 0\end{array}\right)$; therefore by changing the orientation of 
the triples $02d,03e,04a,05b,1ab$ all the trees become positive.

These examples concern the case of a 3-graph $H$ for which the underlying graph $G-\{0\}$ of
$H-\{0\}$ is a cycle of triangles
together with edges each joining an ``inner'' vertex (degree 4) to an
``outer'' vertex (degree 2).  Theorem~\ref{thm:interlaced} below says that
these graphs are non-3-Pfaffian if and only if the cycle of triangles is
even. Moreover, they are all minimal non-3-Pfaffian graphs. 
Recall from  Proposition~\ref{lem:two_ind_white_edges} that if the graph
$G-\{0\}$ underlying $H-\{0\}$ has two independent edges joining pairs of ``outer''
vertices of a cycle of triangles, or $G-\{0\}$ has a 3-cycle of edges joining three
``outer'' vertices of such a cycle of triangles, then the 3-graph $H$ is
non-3-Pfaffian, but in this case it is non-minimal (the minimal examples being
those in Figure~\ref{fig:minimal_non-3-Pfaffian}).

The {\em Lucas numbers} $L_k$ are defined for $k\geq 3$ by
$L_k=L_{k-2}+L_{k-1}$ and $L_1=1, L_2=3$. This sequence is given
explicitly by
$L_k=\left(\frac{1+\sqrt{5}}{2}\right)^k+\left(\frac{1-\sqrt{5}}{2}\right)^k$.

\begin{thm}\label{thm:interlaced}
Let $H$ be the 3-graph on vertices $0,1,1',2,2',\ldots, k,k'$ with triples
$$\{k-1,k, 1'\},
\{k, 1, 2'\}, \{1,2,3'\}, \ldots, \{k-2,k-1,k'\}$$
and $$\{0,1,1'\}, \{0,2,2'\}, \ldots, \{0, k,k'\}.$$

Then $H$ has $L_k$ spanning trees. For odd values of $k$ the 3-graph
$H$ is
3-Pfaffian but for even values of $k\geq 4$ it is
non-3-Pfaffian. Furthermore, when $k\geq 4$ is even $H$ is a minimal
non-3-Pfaffian graph.
\end{thm}
\begin{proof}
The 3-graph $H$ in the case $k=3$ is shown by direct calculation to have a 3-Pfaffian
orientation (see the first entry of
Table~\ref{table:interlaced}) and the cases $k=1$ and $2$ trivially also give 3-Pfaffian 3-graphs. So we assume $k\geq 4$.

Let $s_i=\{0, i,i'\}$ and $t_i=\{i-2,i-1,i'\}$ for
$i=1,\ldots, k$ (in which $t_1=\{k-1, k,
1'\}$, $t_2=\{k, 1,2'\}$). If successive triples $s_i, s_{i+1}$, reading subscripts modulo
$k$, belong to a spanning tree $T$ of $H$ then $s_{i-1}$ must also belong
to $T$. This is because the only triples containing vertex $(i-1)'$ are $s_{i-1}$ and
$t_{i-1}$, and the latter makes a cycle with $s_i$ and $s_{i+1}$. Therefore if there are any
successive triples $s_i$ and $s_{i+1}$ in $T$ then $T$ consists of all
the triples
$s_1,s_2,\ldots, s_k$. For any other spanning tree there are no two
consecutive triples $s_i, s_{i+1}$. On the other hand, given a non-empty subset
$I$ of $\{1,2,\ldots, k\}$ with the property that no two elements are
consecutive (modulo $k$)  the triples $\{s_i:i\in
I\}\cup\{t_j:j\not\in I\}$ form a spanning tree of $H$, which we shall
denote by $T_I$. The singleton subsets $I=\{i\}$ vacuously satisfy the
consecutiveness condition. Since $k\geq 4$ there is at least one such set $I$
with $2$ or more elements. There are
$L_k-1$ such non-empty subsets $I$ uniquely determining spanning trees
$T_I$ in this way. Together with the spanning tree $S$ consisting of triples
$\{s_1,s_2,\ldots, s_k\}$, they account for all $L_k$ spanning trees
of $H$.

Take the vertices of $H$ in the order $0,1',1,2',2, \ldots, k',k$. We shall
choose triple orientations $(0\, i'\, i)$ for the $s_i$ and
$(i-2\, i-1\, i')$ for the $t_i$ and
calculate directly the sign of the spanning tree $T_I$ with set of triples $\{s_i:i\in
I\}\cup\{t_j:j\not\in I\}$ for $I\subseteq\{1,2,\ldots, k\}$ having no
two consecutive elements modulo $k$. Then we shall argue that when $k$
is even no switches of triple orientations can make all the spanning trees the
same sign, whereas the reverse is true when $k$ is odd.  First however
we observe that the spanning tree $S$ with set of triples
$\{s_1,s_2\ldots, s_k\}$ gives the following cyclic permutation of the
vertex set:
$$(\,0\; 1'\; 1\,)\;(\,0\; 2'\; 2'\,)\cdots (\, 0\; k'\; k\,)=(\, 0\; 1'\; 1\; 2'\; 2\; \cdots \; k'\; k\,).$$  
Hence the spanning tree $S$ has positive orientation. 

{\sc Claim.} A spanning tree $T$ of $H$ has sign given by
$$\begin{cases} (-1)^{|I|-1} & T=T_I=\{s_i:i\in I\}\cup\{t_j:j\not\in
  I\},\\
 +1 & T=S=\{s_1,\ldots, s_k\}.\end{cases}$$

We delay the proof of this claim and proceed to determine whether $H$ is
3-Pfaffian or not. If $k$ is odd, switching the orientation of all triples
$s_i$ clearly makes all trees negative, so in this case $H$ is $3$-Pfaffian. 
To treat the case $k$ even it is necessary to look more carefully at the
effect that switching orientations has on the sign of the trees. Since each
tree contains exactly one of $s_i$ and $t_i$,  switching both of them has the
effect of switching the signs of all trees, hence at most one of $s_i$ and
$t_i$ has to be switched. Observe also that if in a given orientation we
switch the triple $t_i$ the signs of the trees are opposite to those obtained
by switching $s_i$. Therefore, we can assume that if $H$ has a $3$-Pfaffian
orientation, then this orientation can be obtained from the initial one by
switching a subset of the triples $s_i$. Now, since we want all the trees
$T_{\{i\}}$ to have the same sign, the only options are to either switch no
$s_i$ or to switch all of them. The first option is clearly not $3$-Pfaffian
and the second one makes all trees but $S$ negative, hence it is also non-$3$-Pfaffian.

We next show that for $k$ even the graph $H$ is minimally
non-3-Pfaffian. By symmetry, it is enough to consider the deletions
$H\backslash s_1$ and $H\backslash t_1$ and the contractions $H/s_1$ and
$H/t_1$. The $3$-graph $H\backslash s_1$   is easily seen to be $3$-Pfaffian
by switching the orientation of $s_2,s_3,\ldots,s_k$ (the spanning tree $S$ of
$H$ is no longer a spanning tree of $H\backslash s_1$). For $H_1=H\backslash t_1$,
we observe that $H_1-0$ has no cycles, and hence satisfies the hypothesis of
Theorem~\ref{thm:PSTS_3-Pfaff}. Since the underlying graph of $H_1-0$ is planar, it
is a Pfaffian graph, hence the theorem implies that $H\backslash t_1$ is
$3$-Pfaffian.  The contraction $H/s_1$ is shown to be $3$-Pfaffian by a
similar argument. Finally, the contraction $H/t_1$ is isomorphic to the
same $3$-graph of the theorem statement corresponding to $k-1$ and with the triple $s_1$ removed, hence it is also $3$-Pfaffian.

It remains to prove the claim about the signs of spanning trees. To calculate the orientation of a spanning
 tree $T_I$ of $H$ we embed $T_I$ in the plane so that positive triple
 orientations correspond to anticlockwise orientations. Traversing $T_I$ in
 an anticlockwise sense and reading off vertices as we encounter them we
 obtain a permutation whose sign is the sign of $T_I$. The $3$-graph $T_I-0$
 has $|I|$ connected components, each of them containing a pair $i,i'$ and
 some triples of the form $t_j$. If $i<j$ are consecutive elements of $I$, the
 component that contains $i,i'$ contains triples $t_{i+1},\ldots, t_{j-1}$;
 all together, the set of vertices is $\{i-1,i',i,(i+1)',i+1,\ldots,
 (j-2)',j-2,(j-1)'\}$, that is, the interval $[i-1,(j-1)']$. These intervals
partition $[1',k]$ (cyclically), so in
 order to determine the orientation of $T_I$ we find the orientation
 on each interval $[i-1,(j-1)']$  separately and then piece these
 together and start the traversal of the whole tree $T_I$ at the vertex $0$.
We can assume that
 $1\in I$, since the spanning trees in the orbit of $T_I$ under the
 permutation $(1,2,\ldots,k)(1',2',\ldots,k')$ have the same sign as $T_I$.
We say that two elements of $I$ are cyclically consecutive if there is no
element of $I$ between them. Hence, the largest element of $I$ is consecutive
with $1$, unless if $I=\{1\}$; in this case, we consider that $k+1$ is
consecutive with $1$. We identify throughout $k+1$ with the vertex $1$  and
$0$ with vertex $k$ in order that
 the following argument works for all $|I|$ consecutive pairs of elements of
 $I$ rather than having to consider two cases separately.

We shall use the following facts about the signs of permutations:
\begin{enumerate}[(i)]
\item reversing the order of $\ell$ elements has sign $(-1)^{\binom{\ell}{2}}$,
\item the permutation 
$$\left(\begin{array}{cccccccc} a_1 & a_2 & \cdots & a_\ell & b_1 & b_2
  & \cdots & b_\ell\\
b_1 & a_1 & b_2 & a_2 & \cdots & \cdots & b_\ell & a_\ell\end{array}\right)$$ 
interleaving  a block of $\ell$ elements with another block of
  $\ell$ elements has sign $(-1)^{\binom{\ell+1}{2}}$.
\end{enumerate}

Traversing the tree $T_I$ in an anticlockwise sense we find that the vertices
$\{i-1, i,i',\ldots, (j-2)',j-2,(j-1)' \}$ appear in the following order up to
even permutation:

\begin{equation}\label{eq:perm}
i', i, (i\!+1)',i\!-1,i\!+1,i\!+2,\ldots,j\!-2,(j\!-1)',(j\!-2)',\ldots, (i\!+3)',(i\!+2)'. 
\end{equation}

It is easily seen that this permutation is of opposite parity as

$$i\!-1,i, i\!+1,i\!+2,\ldots,j\!-2,(j\!-1)',(j\!-2)',\ldots,
(i\!+3)',(i\!+2)', (i\!+1)', i'. $$

Now we reverse the order of the last $j-i$ elements, with a sign change of
$(-1)^{\binom{j-i}{2}}$.

$$i\!-1,i, i\!+1,i\!+2,\ldots,j\!-2,i',(i\!+1)',\ldots, (j\!-2)',(j\!-1)'. $$

Finally, we interleave the block $i,i\!+1,i\!+2,\ldots,j\!-2$ with the block
$i',(i\!+1)',\ldots, (j\!-2)'$, with a sign change of $(-1)^{\binom{j-i}{2}}$ again.

Therefore, permutation~(\ref{eq:perm}) is of opposite parity as
$i-1,i',i,\ldots,(j\!-2)',j\!-2,(j\!-1)'$.

Now, given the tree $T_I$ embedded in the plane, when traversing in
counterclockwise order we encounter the following permutation of the vertex set,
after a change sign of $(-1)^{|I|}$
$$0,k,1',1,2',2,\ldots, k'. $$
Clearly this permutation is of opposite parity as $0,1',1,\ldots,k',k$,
therefore the tree $T_I$ has sign $(-1)^{|I|-1}$, as claimed.\hspace*{\fill} \end{proof}


\section{Some open problems}\label{sec:open_problems}


In Theorem~\ref{thm:Tutte-like} we found a necessary condition for a 3-graph
to have a spanning tree. Is it possible to strengthen this condition to make
it sufficient?

\begin{problem}
Find a characterization of 3-graphs that have a spanning 
tree.
\end{problem}

In Section~\ref{sec:existence_counting} we found a superexponential lower bound of the type asked for in the following problem for the cases
$d=\frac{n-1}{2}$, $m=1$ (Steiner triple systems, Theorem~\ref{thm:lower_bound_STS}) and
$d=\binom{n}{2}, m=n-2$ (the complete 3-graph,
Theorem~\ref{thm:pruefer}).
\begin{problem}\label{prob:lower_bound}
Suppose that $H=([2n\!+\!1],\Delta)$ is a 3-graph such that each vertex is of degree
at least $d$ and each pair of vertices has multiplicity at least $m$. Find a lower bound on the
number of spanning trees of $H$ (as a function of $n, d$ and $m$).  
\end{problem}


In Section~\ref{sec:suspensions} we considered suspensions of graphs. A suspension of a bipartite graph $G=(A\cup B, E)$ is a
  special form of tripartite 3-graph $H=(A\cup B\cup C,\Delta)$, where all triples are of the form
  $abc$ for $a\in A, b\in B, c\in C$. For $k$-regular bipartite graphs Schrijver~\cite{Schrijver98}
  established exponential lower bounds on the number of perfect
  matchings. It may be easier to solve
  Problem~\ref{prob:lower_bound} when restricted to the case
  when $H$ is a tripartite 3-graph.

Little~\cite{Little75} gave a forbidden subgraph characterization of
bipartite Pfaffian graphs: if there is an even subdivision of $K_{3,3}$ whose
  complement has a perfect matching then the graph is non-Pfaffian.

\begin{problem}\label{prob:3-Pfaffian_tripartite} Is there a forbidden subgraph characterization
  for 3-Pfaffian tripartite 3-graphs?\end{problem}
  Corollary~\ref{cor:2-susp_obstructions} provides an affirmative
  answer to the question raised in
  Problem~\ref{prob:3-Pfaffian_tripartite} for
  2-suspensions and likewise Theorem~\ref{thm:1-susp} for 1-suspensions.

Theorem~\ref{thm:interlaced} provides an example of an infinite
set of non-3-Pfaffian partial Steiner triple systems, no two of which can be obtained from the
other by deletion or contraction of triples. 
Each partial Steiner triple system $H$ belonging to this set has the
property that if $v$ is a vertex such that $H-v$ has no cycles, then the underlying
graph $G-v$ is non-planar. 
\begin{problem}\label{prob:minimalSTS}
Let $\mathcal{H}$ be the set of non-3-Pfaffian partial Steiner triple systems
$H=(V,\Delta)$ with the property that there is $v\in V$ such that
$H-v$ has no cycles and the underlying
graph $G-v$ of $H-v$ is planar. Is there an infinite set of 3-graphs in $\mathcal{H}$
that are minimal with respect to deletion and contraction of triples?
\end{problem}


\section*{Acknowledgements}

We would like to thank Marc Noy for his generosity and the impetus he
gave to our work on the subject of this paper -- many of the results would not have been
obtained without him. Also, we thank Martin Loebl for the insights he
shared on the topic of Pfaffian orientations and the many ideas he
produced that inspired our work.


\bibliographystyle{elsarticle-num}

\begin{thebibliography}{10}
\expandafter\ifx\csname url\endcsname\relax
  \def\url#1{\texttt{#1}}\fi
\expandafter\ifx\csname urlprefix\endcsname\relax\def\urlprefix{URL }\fi
\expandafter\ifx\csname href\endcsname\relax
  \def\href#1#2{#2} \def\path#1{#1}\fi

\bibitem{Abd04}
A.~Abdesselam, Grassmann--Berezin calculus and theorems of the matrix-tree type, Adv. in Appl. Math. 33 (2004), 51--70.


\bibitem{AKS97}
N.~Alon, J.-H. Kim, J.~Spencer, Nearly perfect matchings in regular simple
  hypergraphs, Israel J. Math. 100~(1) (1997) 171--187.


\bibitem{AF95}
L.~Andersen, H.~Fleischner, The {NP}-completeness of finding {A}-trails in
  {Eulerian} graphs and of finding spanning trees in hypergraphs, Discrete
  Appl. Math. 59 (1995) 203--214.


\bibitem{Brouwer81}
A.~E. Brouwer, On the size of a maximum transversal in a {Steiner} triple
  system, Canadian J. Math. 33 (1981) 1202--1204.


\bibitem{CMSS08}
S.~Caracciolo, G.~Masbaum, A.~Sokal, A.~Sportiello, A randomized
  polynomial-time algorithm for the spanning hypertree problem on 3-uniform
  hypergraphs, arXiv:math/0812.3593v1 [cs.CC] (2008).


\bibitem{E65}
J. Edmonds, Paths, trees, and flowers, Canad. J. Math. 17 (1965), 449--467.

\bibitem{GS85}
H.~Gabow, M.~Stallmann, Efficient algorithms for graphic matroid
intersection and parity, Automata, Languages and Programming; Lecture
Notes in Computer Science 194, (ed. W. Brauer), Springer-Verlag, New
York, 1985, 210--220.
 
\bibitem{GS86}
H.~Gabow, M.~Stallmann, Augmenting path algorithm for linear matroid parity,
  Combinatorica 6~(2) (1986) 123--150.

\bibitem{GL99}
A.~Galluccio, M.~Loebl, On the theory of {Pfaffian} orientations. {I.}
  {Perfect} matchings and permanents, Electron. J. Combin. 6 (1),
  1999, R6.


\bibitem{HR04}
S.~Hirschman, V.~Reiner, Note on the {Pfaffian} matrix-tree theorem, Graphs
  Combin. 20~(1) (2004) 59--63.


\bibitem{JK82}
P.~Jensen, B.~Korte, Complexity of matroid property algorithms, SIAM J. Comput.
  11 (1982) 184--190.


\bibitem{K67} P.W.~Kasteleyn, Graph theory and crystal physics. In: Graph Theory and Theoretical
Physics (ed. F. Harary), Academic Press, New York, 1967, 43--110.


\bibitem{Little75}
C.~Little, A characterization of convertible (0, 1)-matrices, J. Combin. Theory
  B 18 (1975) 187--208.



\bibitem{L78}
L.~Lov\'asz, The matroid matching problem, in: Algebraic Methods in Graph
  Theory, Vol. I, II, Colloquia Mathematica Societatis J\'anos Bolyai, Szeged,
  Hungary, 1978, pp. 495--517.


\bibitem{L80}
L.~Lov\'{a}sz, Matroid matching and some applications, J. Combin. Theory Ser. B
  28~(2) (1980) 208--236.


\bibitem{LP86}
L.~Lov\'asz, M.~Plummer, Matching Theory, Vol.~29 of Annals of Discrete
  Mathematics, North-Holland, 1986, also at: Akad\'emia Kiad\'o, Budapest,
  1986.

\bibitem{Lovasz87}
L.~Lov\'{a}sz, Matching structure and the matching lattice, J. Combin. Theory
  Ser. B 43~(2) (1987) 187--222.

\bibitem{MV02}
G.~Masbaum, A.~Vaintrob, A new matrix-tree theorem, Internat. Math. Res.
  Notices 27 (2002) 1397--1426.

\bibitem{Schrijver98}
A.~Schrijver, Counting 1-factors in regular bipartite graphs, J. Combin. Theory
  B 72~(1) (1998) 122--135.

\bibitem{SS06}
S.~Sivasubramanian, Spanning trees in complete uniform hypergraphs and a
  connection to extended {$r$-Shi} hyperplane arrangements,
  arXiv:math/0605083v2 [math.CO] (2006).


\bibitem{Stem90}
J.~Stembridge, Nonintersecting paths, {Pfaffians} and plane partitions, Adv.
  Math. 83 (1990) 96--131.



\bibitem{RT06}
R. Thomas, A survey of {Pfaffian} orientations of graphs, in:
Proceedings of the International Congress of Mathematicians, Madrid,
Spain, European Mathematical Society, Z\"{u}rich, 2006.

\bibitem{Val79}
L.~Valiant, The complexity of enumeration and reliability problems, SIAM J.
  Comput. 8~(3) (1979) 410--421.


\bibitem{VY89}
V.~Vazirani, M.~Yannakakis, Pfaffian orientations, 0-1 permanents, and even
  cycles in directed graphs, Discrete Appl. Math. 25~(1--2) (1989) 179--190.



\bibitem{Wilson74}
R.~Wilson, Nonisomorphic {Steiner} triple systems, Math. Z. 135~(4) (1974)
  303--313.


\end{thebibliography}

\end{document}